\magnification=1200

\hsize=11.25cm    
\vsize=18cm       
\parindent=12pt   \parskip=5pt     

\hoffset=.5cm   
\voffset=.8cm   

\pretolerance=500 \tolerance=1000  \brokenpenalty=5000

\catcode`\@=11

\font\eightrm=cmr8         \font\eighti=cmmi8
\font\eightsy=cmsy8        \font\eightbf=cmbx8
\font\eighttt=cmtt8        \font\eightit=cmti8
\font\eightsl=cmsl8        \font\sixrm=cmr6
\font\sixi=cmmi6           \font\sixsy=cmsy6
\font\sixbf=cmbx6

\font\tengoth=eufm10 
\font\eightgoth=eufm8  
\font\sevengoth=eufm7      
\font\sixgoth=eufm6        \font\fivegoth=eufm5

\skewchar\eighti='177 \skewchar\sixi='177
\skewchar\eightsy='60 \skewchar\sixsy='60

\newfam\gothfam           \newfam\bboardfam

\def\tenpoint{
  \textfont0=\tenrm \scriptfont0=\sevenrm \scriptscriptfont0=\fiverm
  \def\rm{\fam\z@\tenrm}
  \textfont1=\teni  \scriptfont1=\seveni  \scriptscriptfont1=\fivei
  \def\oldstyle{\fam\@ne\teni}\let\old=\oldstyle
  \textfont2=\tensy \scriptfont2=\sevensy \scriptscriptfont2=\fivesy
  \textfont\gothfam=\tengoth \scriptfont\gothfam=\sevengoth
  \scriptscriptfont\gothfam=\fivegoth
  \def\goth{\fam\gothfam\tengoth}
  
  \textfont\itfam=\tenit
  \def\it{\fam\itfam\tenit}
  \textfont\slfam=\tensl
  \def\sl{\fam\slfam\tensl}
  \textfont\bffam=\tenbf \scriptfont\bffam=\sevenbf
  \scriptscriptfont\bffam=\fivebf
  \def\bf{\fam\bffam\tenbf}
  \textfont\ttfam=\tentt
  \def\tt{\fam\ttfam\tentt}
  \abovedisplayskip=12pt plus 3pt minus 9pt
  \belowdisplayskip=\abovedisplayskip
  \abovedisplayshortskip=0pt plus 3pt
  \belowdisplayshortskip=4pt plus 3pt 
  \smallskipamount=3pt plus 1pt minus 1pt
  \medskipamount=6pt plus 2pt minus 2pt
  \bigskipamount=12pt plus 4pt minus 4pt
  \normalbaselineskip=12pt
  \setbox\strutbox=\hbox{\vrule height8.5pt depth3.5pt width0pt}
  \let\bigf@nt=\tenrm       \let\smallf@nt=\sevenrm
  \normalbaselines\rm}

\def\eightpoint{
  \textfont0=\eightrm \scriptfont0=\sixrm \scriptscriptfont0=\fiverm
  \def\rm{\fam\z@\eightrm}
  \textfont1=\eighti  \scriptfont1=\sixi  \scriptscriptfont1=\fivei
  \def\oldstyle{\fam\@ne\eighti}\let\old=\oldstyle
  \textfont2=\eightsy \scriptfont2=\sixsy \scriptscriptfont2=\fivesy
  \textfont\gothfam=\eightgoth \scriptfont\gothfam=\sixgoth
  \scriptscriptfont\gothfam=\fivegoth
  \def\goth{\fam\gothfam\eightgoth}
  
  \textfont\itfam=\eightit
  \def\it{\fam\itfam\eightit}
  \textfont\slfam=\eightsl
  \def\sl{\fam\slfam\eightsl}
  \textfont\bffam=\eightbf \scriptfont\bffam=\sixbf
  \scriptscriptfont\bffam=\fivebf
  \def\bf{\fam\bffam\eightbf}
  \textfont\ttfam=\eighttt
  \def\tt{\fam\ttfam\eighttt}
  \abovedisplayskip=9pt plus 3pt minus 9pt
  \belowdisplayskip=\abovedisplayskip
  \abovedisplayshortskip=0pt plus 3pt
  \belowdisplayshortskip=3pt plus 3pt 
  \smallskipamount=2pt plus 1pt minus 1pt
  \medskipamount=4pt plus 2pt minus 1pt
  \bigskipamount=9pt plus 3pt minus 3pt
  \normalbaselineskip=9pt
  \setbox\strutbox=\hbox{\vrule height7pt depth2pt width0pt}
  \let\bigf@nt=\eightrm     \let\smallf@nt=\sixrm
  \normalbaselines\rm}

\tenpoint

\def\pc#1{\bigf@nt#1\smallf@nt}         \def\pd#1 {{\pc#1} }

\catcode`\;=\active
\def;{\relax\ifhmode\ifdim\lastskip>\z@\unskip\fi
\kern\fontdimen2  -1.2 \fontdimen3 \string;}

\catcode`\:=\active
\def:{\relax\ifhmode\ifdim\lastskip>\z@\unskip\fi\penalty\@M\ \fi\string:}

\catcode`\!=\active
\def!{\relax\ifhmode\ifdim\lastskip>\z@
\unskip\fi\kern\fontdimen2  -1.1 \fontdimen3 \string!}

\catcode`\?=\active
\def?{\relax\ifhmode\ifdim\lastskip>\z@
\unskip\fi\kern\fontdimen2  -1.1 \fontdimen3 \string?}

\frenchspacing

\def\raggedbottom{\topskip 10pt plus 36pt\r@ggedbottomtrue}

\def\pointir{\unskip . --- \ignorespaces}

\def\Medbreak{\vskip-\lastskip\medbreak}

\long\def\th#1 #2\enonce#3\endth{
   \Medbreak\noindent
   {\pc#1} {#2\unskip}\pointir{\it #3}\smallskip}

\def\proof{\vskip-\lastskip\smallskip\noindent
 {\it Proof} : }

\def\decale#1{\smallbreak\hskip 28pt\llap{#1}\kern 5pt}
\def\decaledecale#1{\smallbreak\hskip 34pt\llap{#1}\kern 5pt}
\def\puce{\smallbreak\hskip 6pt{$\scriptstyle\bullet$}\kern 5pt}

\def\eqalign#1{\null\,\vcenter{\openup\jot\m@th\ialign{
\strut\hfil$\displaystyle{##}$&$\displaystyle{{}##}$\hfil
&&\quad\strut\hfil$\displaystyle{##}$&$\displaystyle{{}##}$\hfil
\crcr#1\crcr}}\,}

\catcode`\@=12

\showboxbreadth=-1  \showboxdepth=-1

\newcount\numerodesection \numerodesection=1
\def\section#1{\bigbreak
 {\bf\number\numerodesection.\ \ #1}\nobreak\medskip
 \advance\numerodesection by1}

\mathcode`A="7041 \mathcode`B="7042 \mathcode`C="7043 \mathcode`D="7044
\mathcode`E="7045 \mathcode`F="7046 \mathcode`G="7047 \mathcode`H="7048
\mathcode`I="7049 \mathcode`J="704A \mathcode`K="704B \mathcode`L="704C
\mathcode`M="704D \mathcode`N="704E \mathcode`O="704F \mathcode`P="7050
\mathcode`Q="7051 \mathcode`R="7052 \mathcode`S="7053 \mathcode`T="7054
\mathcode`U="7055 \mathcode`V="7056 \mathcode`W="7057 \mathcode`X="7058
\mathcode`Y="7059 \mathcode`Z="705A


\def\diagram#1{\def\normalbaselines{\baselineskip=0pt\lineskip=5pt}
\matrix{#1}}

\def\vfl#1#2#3{\llap{$\textstyle #1$}
\left\downarrow\vbox to#3{}\right.\rlap{$\textstyle #2$}}

\def\ufl#1#2#3{\llap{$\textstyle #1$}
\left\uparrow\vbox to#3{}\right.\rlap{$\textstyle #2$}}

\def\hfl#1#2#3{\smash{\mathop{\hbox to#3{\rightarrowfill}}\limits
^{\textstyle#1}_{\textstyle#2}}}

\def\agoth{{\goth a}}
\def\bgoth{{\goth b}}
\def\dgoth{{\goth d}}
\def\ogoth{{\goth o}}
\def\units{{\ogoth^\times}}

\def\pgoth{{\goth p}}
\def\lgoth{{\goth l}}

\def\P{{\bf P}}
\def\Q{{\bf Q}}
\def\Qp{\Q_p}
\def\Ql{\Q_l}
\def\R{{\bf R}}
\def\C{{\bf C}}
\def\N{{\bf N}}

\def\Z{{\bf Z}}
\def\Zp{\Z_p}
\def\F{{\bf F}}
\def\Fp{{\F_{\!p}}}

\def\Id{\mathop{\rm Id}\nolimits}

\def\Card{\mathop{\rm Card}\nolimits}
\def\Gal{\mathop{\rm Gal}\nolimits}
\def\Ker{\mathop{\rm Ker}\nolimits}
\def\Coker{\mathop{\rm Coker}\nolimits}

\def\droite#1{\,\hfl{#1}{}{8mm}\,}
\def\series#1{(\!(#1)\!)}
\def\zero{\{0\}}
\def\one{\{1\}}

\def\to{\rightarrow}

\def\normressym(#1,#2)_#3{\displaystyle\left({#1,#2\over#3}\right)}

\def\mod{\mathop{\rm mod.}\nolimits}
\def\pmod#1{\;(\mod#1)}

\newcount\refno 
\long\def\ref#1:#2<#3>{                                        
\global\advance\refno by1\par\noindent                              
\llap{[{\bf\number\refno}]\ }{#1} \pointir{\it #2} #3\goodbreak }

\def\citer#1(#2){[{\bf\number#1}\if#2\empty\relax\else,\ {#2}\fi]}

\newbox\bibbox
\setbox\bibbox\vbox{\bigbreak
\centerline{{\pc BIBLIOGRAPHIC} {\pc REFERENCES}}

\ref{\pc ANONYME}:
Rapport sur l'article 37-08,
<Papiers secrets des Annales de la Facult{\'e} des Sciences de Toulouse, 20
janvier 2009.> 
\newcount\anonyme \global\anonyme=\refno

\ref{\pc ARTIN} (E.) and {\pc SCHREIER} (O.):
{\"U}ber eine Kennzeichnung der reell abgeschlossenen K{\"o}rper,
<Abh.\ Math.\ Sem.\ Hamburg {\bf 5}, 1927, pp.~225--232.>
\newcount\artinschreier \global\artinschreier=\refno

\ref{\pc BERLEKAMP} (E.):
An analog[ue] to the discriminant over fields of characteristic two,  
<J.\ Algebra {\bf 38} (1976), no.~2, pp.~325--327.>
\newcount\berlekamp \global\berlekamp=\refno

\ref{\pc BOURBAKI} (N.):
Alg{\`e}bre, Chapitres 4 {\`a} 7,
<Masson, Paris, 1981, 422~pp.>
\newcount\bourbaki \global\bourbaki=\refno

\ref{\pc BRILL} (A.):
{U}eber die Discriminante,
<Math.\ Annalen, {\bf 12}, 1877, pp.~87--89.>
\newcount\brill \global\brill=\refno

\ref{\pc CASSELS} (J.):
Local Fields,
<Cambridge University Press, 1986, 360 pp.>
\newcount\cassels \global\cassels=\refno

\ref{\pc CHILDS} (L.):
Taming wild extensions,
<American Mathematical Society, 2000. viii+215 pp.>
\newcount\childs \global\childs=\refno

\ref{\pc COHEN} (H.):
Number theory. Vol. I. 
<Springer, 2007, xxiv+650 pp.>
\newcount\cohen \global\cohen=\refno

\ref{\pc DALAWAT} (C.):
Minimal bases, minimal cubics,
<in preparation.>
\newcount\minimal \global\minimal=\refno

\ref{\pc DALAWAT} (C.):
Wilson's theorem,
<J. de Th{\'e}orie de nombres de
 Bordeaux, {\bf 21} (2009) 3, pp.~517--521. Cf.~arXiv\string:0711.3879v1.> 
\newcount\wilson \global\wilson=\refno

\ref{\pc DALEN} (K.):
On a theorem of Stickelberger,
<Math.\ Scand.\ {\bf 3}, 1955, pp.~124--126.>
\newcount\dalen \global\dalen=\refno

\ref{\pc DEDEKIND} (R.):
{\"U}ber die Diskriminanten endlicher K{\"o}rper,
<Abh.\ der K{\"o}nigl.\ Ges.\ d.\ Wiss.\ zu G{\"o}ttingen,
{\bf 29} (1882), pp.~1--56.>
\newcount\dedekind \global\dedekind=\refno

\ref{\pc DEL \pc CORSO} (I.) and {\pc DVORNICICH} (R.):
The compositum of wild extensions of local fields of prime degree,
<Monatsh.\ Math.\ {\bf 150} (2007), no.~4, pp.~271--288.>
\newcount\delcorso \global\delcorso=\refno

\ref{\pc EISENSTEIN} (G.): 
{\"U}ber ein einfaches Mittel zur Auffindung der h{\"o}heren
Reciprocit{\"a}tsgesetze und der mit ihnen zu verbindenen 
Erg{\"a}nzungss{\"a}tze, 
<J.\ f.\ d.\ reine und angewandte Math., {\bf 39}, 1850, pp.~351-364. $=$
Math.\ Werke II, pp.~623--636.>
\newcount\eisenstein \global\eisenstein=\refno

\ref{\pc FESENKO} (I.) and {\pc VOSTOKOV} (S.):
Local fields and their extensions,
<American Mathematical Society, 2002. xii+345 pp.>
\newcount\fesvost \global\fesvost=\refno

\ref{\pc FREI} (G.) and {\pc ROQUETTE} (P.):
Emil Artin und Helmut Hasse, Die Korrespondenz 1923--1934,
<Universit{\"a}tsverlag G{\"o}ttingen, 2009, 499 pp.>
\newcount\freiroquette \global\freiroquette=\refno

\ref{\pc FR{\"O}HLICH} (A.):
Discriminants of algebraic number fields,
<Math.\ Zeitschrift, {\bf 74}, 1960, pp.~18--28.>
\newcount\frohlich \global\frohlich=\refno

\ref{\pc FURTW{\"A}NGLER} (Ph.):
{\"U}ber die Reziprozit{\"a}tsgesetze zwischen $l$-ten Potenzresten in
algebraischen Zahlk{\"o}rpern, wenn $l$ eine ungerade Primzahl bedeutet,
<Math.\ Ann.\ {\bf 58}, 1903, pp.~1--50.>
\newcount\furtwangler \global\furtwangler=\refno

\ref{\pc GRAS} (G.):
Class field theory,
<Springer, Berlin, 2003, xiv+491 pp.>
\newcount\gras \global\gras=\refno

\ref{\pc HASSE} (H.):
Bericht {\"u}ber neuere Untersuchungen und Probleme aus der Theorie de
algebraischen Zahlk{\"o}rper, Teil Ia,
<Jahresbericht der Deutschen Mathematiker-Vereinigung, {\bf 436}, 1927,
pp.~233--311.>
\newcount\klassenkorperbericht \global\klassenkorperbericht=\refno

\ref{\pc HASSE} (H.):
{\"U}ber das Reziprozit{\"a}tsgesetz der $m$-ten Potenzreste,
<J.\ f.\ d.\ reine und angewandte Math., {\bf 158}, 
1927, pp.~228--259.>
\newcount\hassepotenz \global\hassepotenz=\refno


\ref{\pc HASSE} (H.):
Zahlentheorie,
<Akademie-Verlag, Berlin, 1969, 611~pp.>
\newcount\hasse \global\hasse=\refno

\ref{\pc HECKE} (E.):
Vorlesungen {\"u}ber die Theorie der algebraischen Zahlen,
<Akad.\ Verlagsges., Leipzig, 1923, pp.~265.> 
\newcount\hecke \global\hecke=\refno

\ref{\pc HENSEL} (K.):
{\"U}ber die zu einem algebraischen K{\"o}rper geh{\"o}rigen Invarianten,
<J.\ f.\ d.\ reine und angewandte Math., {\bf 129}, 1905,
pp.~68--85.> 
\newcount\hensel \global\hensel=\refno

\ref{\pc HENSEL} (K.):
Die multiplikative Darstellung der algebraischen Zahlen f{\"u}r den Bereich
eines beliebigen Primteilers,
<J.\ f.\ d.\ reine und angewandte Math., {\bf 146}, 1916,
pp.~189--215.> 
\newcount\henselmult \global\henselmult=\refno

\ref{\pc HENSEL} (K.):
{\"U}ber die Zerlegung der Primteiler in relativ zyklischen K{\"o}rpern, nebst
einer Anwendung auf die\/ {\rm Kummer}schen K{\"o}rper,
<J.\ f.\ d.\ reine und angewandte Math., {\bf 151}, 1921, pp.~112--120.>
\newcount\henselkummer \global\henselkummer=\refno

\ref{\pc HENSEL} (K.):
Die Zerlegung der Primteiler eines beliebigen Zahlk{\"o}rpers in einem
aufl{\"o}sbaren Oberk{\"o}rper, 
<J.\ f.\ d.\ reine und angewandte Math., {\bf 151}, 1921, pp.~200--209.>
\newcount\henselsolvable \global\henselsolvable=\refno

\ref{\pc HENSEL} (K.):
Zur multiplikativen Darstellung der algebraischen Zahlen f{\"u}r den Bereich
eines Primteilers,
<J.\ f.\ d.\ reine und angewandte Math., {\bf 151}, 1921, pp.~210--212.>
\newcount\henselmultii \global\henselmultii=\refno

\ref{\pc  HERBRAND} (J.):
Sur la th{\'e}orie des groupes de d{\'e}composition, d'inertie et de
ramification, 
<J.\ de math. pures et  appl., {\bf 10}, 1931, pp.~481--498.>
\newcount\herram \global\herram=\refno

\ref{\pc HERBRAND} (J.):
Une propri{\'e}t{\'e} du discriminant des corps alg{\'e}brique,
<Ann.\ sci.\ {\'E}cole Norm.\ Sup.\ {\bf 49}, 1932, pp.~105--112.>
\newcount\herbrand \global\herbrand=\refno

\ref{\pc HILBERT} (D.):
Grundz{\"u}ge einer Theorie des Galoisschen Zahlk{\"o}rpers,
<Nachrichten d.\ Ges.\ d.\ Wiss.\ zu G{\"o}ttingen, 1894, pp.~224--236.  $=$
Ges.\ Abhandlungen I, pp.~13--23.>
\newcount\hilbert \global\hilbert=\refno

\ref{\pc HILBERT} (D.):
Die Theorie der algebraischen Zahlk{\"o}rper,
<Jahresbericht der Deutschen Mathematiker-Vereinigung, {\bf 4}, 1897,
pp.~175--546. $=$
Ges.\ Abhandlungen I, pp.~63--363.>
\newcount\zahlbericht \global\zahlbericht=\refno

\ref{\pc KRAUS} (A.):
Quelques remarques {\`a} propos des invariants $c_4$, $c_6$ et $\Delta$ d'une
courbe elliptique. 
<Acta Arith.\ {\bf 54} (1989), no.~1, pp.~75--80.>
\newcount\krausinv \global\krausinv=\refno

\ref{\pc KRAUS} (A.):
Sur le d{\'e}faut de semi-stabilit{\'e} des courbes elliptiques {\`a}
r{\'e}duction additive, 
<Manuscripta Math.\ {\bf 69} (1990), no.~4, pp.~353--385.>
\newcount\kraus \global\kraus=\refno

\ref{\pc MARTINET} (J.):
Les discriminants quadratiques et la congruence de Stickelberger,
<S{\'e}m.\ Th{\'e}or.\ Nombres Bordeaux (2) {\bf 1} (1989),  no.~1,
pp.~197--204.> 
\newcount\martinet \global\martinet=\refno

\ref{\pc MOVAHHEDI} (A.) and {\pc ZAHIDI} (M.):
Symboles des restes quadratiques des discriminants dans les extensions
mod{\'e}r{\'e}ment ramifi{\'e}es,
<Acta Arith.\ {\bf 92}, 2000, no.~3, pp.~239--250.>
\newcount\zahidi \global\zahidi=\refno

\ref{\pc NAIPAUL} (V.):
A way in the world,
<William Heinemann, London, 1994.>
\newcount\naipaul \global\naipaul=\refno

\ref{\pc NAKAGOSHI} (N.):
The structure of the multiplicative group of residue classes modulo
${\pgoth}\sp{N+1}$,  
<Nagoya Math.\ J.\ {\bf 73} (1979), pp.~41--60.>
\newcount\nakagoshi \global\nakagoshi=\refno 

\ref{\pc NARKIEWICZ} (W.):
Elementary and analytic theory of algebraic numbers,
<$2^{\rm nd}$ ed., Springer-Verlag, 
1990, 746~pp.> 
\newcount\narkiewicz \global\narkiewicz=\refno 

\ref{\pc NEUKIRCH} (J.):
Class Field Theory,
<Springer-Verlag, Berlin, 1986, 140~pp.>
\newcount\neukirch \global\neukirch=\refno

\ref{\pc PELLET} (A.):
Sur la d{\'e}composition d'une fonction enti{\`e}re en facteurs
irr{\'e}ducibles suivant un module permier, 
<Comptes Rendus {\bf 86}, 1878, pp.~1071--1072.>
\newcount\pellet \global\pellet=\refno

\ref{\pc PISOLKAR} (S.):
Absolute norms of\/ $p$-primary numbers,
<J. de Th{\'e}orie de nombres de Bordeaux, {\bf 21} (2009) 3,
pp.~733--740. Cf.~arXiv\string:0807.1174v1.>  
\newcount\supriya \global\supriya=\refno

\ref{\pc ROQUETTE} (P.):
Letter to the author,
<8 July 2007.>
\newcount\roquette \global\roquette=\refno 

\ref{\pc SCHUR} (I.):
Elementarer Beweis eines Satzes von {\pc L.\ \pc STICKELBERGER},
<Math.\ Zeitschrift, {\bf 29}, 1929, pp.~464--465.>
\newcount\schur \global\schur=\refno

\ref{\pc SEKIGUCHI} (T.), {\pc OORT} (F.) and {\pc SUWA} (N.):
On the deformation of Artin-Schreier to Kummer,
<Ann.\ sci.\ {\'E}cole Norm.\ Sup.\ (4) {\bf 22}, 1989, no.~3, pp.~345--375.>
\newcount\sekiguchi \global\sekiguchi=\refno

\ref{\pc SERRE} (J.-P.):
Corps locaux,
<Publications de l'Universit{\'e} de Nancago, No.~{\sevenrm VIII}, Hermann,
Paris, 1968, 245 pp.>
\newcount\serre \global\serre=\refno

\ref{\pc SERRE} (J.-P.):
Une ``formule de masse" pour les extensions totalement ramifi{\'e}es de
degr{\'e} donn{\'e} d'un corps local, 
<Comptes Rendus {\bf 286}, 1978, pp.~1031--1036.>
\newcount\serremass \global\serremass=\refno

\ref{\pc SILVERMAN} (J.):
Weierstrass equations and the minimal discriminant of an elliptic curve,
<Mathematika {\bf 31}, 1984, no.~2, pp.~245--251.>
\newcount\silverman \global\silverman=\refno

\ref{\pc SKOLEM} (Th.):
On a certain connection between the discriminant of a polynomial and the
number of its irreducible factors\/ $\pmod p$,
<Norsk Mat.\ Tidsskr.\ {\bf 34}, 1952, pp.~81--85.>
\newcount\skolem \global\skolem=\refno

\ref{\pc STICKELBERGER} (L.):
{\"U}ber eine neue Eigenschaft der Diskriminanten algebraischer
Zahlk{\"o}rper, 
<Verh.\ I. int.\ Math.\ Kongr., 1897, pp.~182--193.>
\newcount\stick \global\stick=\refno

\ref{\pc SWAN} (R.):
Factorization of polynomials over finite fields, 
<Pacific J.\ Math.\ {\bf 12}, 1962, pp.~1099--1106.>
\newcount\swan \global\swan=\refno

\ref {\pc TATE} (J.) :
The arithmetic of elliptic curves,
<Invent.\ Math.\ {\bf 23}, 1974, pp.~179--206.>
\newcount\tate \global\tate=\refno

\ref{\pc VORONO{\"I}} (G.):
Sur une propri{\'e}t{\'e} du discriminant des fonctions enti{\`e}res,
<Verh.\ des 3.~int.\ Math.\ Kongr., 1905, pp.~187--189.>
\newcount\voronoi \global\voronoi=\refno

\ref{\pc VOSTOKOV} (S.):
Explicit formulas for the Hilbert symbol. 
< in Invitation to higher local fields,
Geom.\ Topol.\ Monogr., 3, Geom.\ Topol.\ Publ., Coventry, 2000, pp.~81--89.>
\newcount\vostokov \global\vostokov=\refno

\ref{\pc WADSWORTH} (A.):
Discriminants in characteristic two,
<Linear and Multilinear Algebra {\bf 17} (1985), no.~3-4, pp.~235--263.>
\newcount\wadsworth \global\wadsworth=\refno

\ref{\pc WATERHOUSE} (W.):
A unified Kummer-Artin-Schreier sequence,
<Math.\ Ann.\ {\bf 277}, 1987, no.~3, pp.~447--451.>
\newcount\waterhouse \global\waterhouse=\refno

\ref{\pc WATERHOUSE} (W.):
Discriminants of {\'e}tale algebras and related structures,
<J.\ f.\ d.\ reine und angewandte Math.\ {\bf 379}, 1987, pp.~209--220.>
\newcount\waterdisc \global\waterdisc=\refno

\ref{\pc WITT} (E.):
Konstruktion der galoischen K{\"o}rpern der Charakteristik~$p$ zu gegebener
Gruppe der Ordnung~$p^f$,
<J.\ f.\ d.\ reine und angewandte Math. {\bf 174}, 1936,
pp.~237--245.> 
\newcount\witt \global\witt=\refno

} 

\centerline{\bf Local discriminants, kummerian extensions, and elliptic curves}
\bigskip\bigskip 
\leftline{Chandan Singh Dalawat} 
\leftline{\it Harish-Chandra Research Institute}
\leftline{\it Chhatnag Road, Jhunsi,   Allahabad 211019, India} 
\leftline{\it Email~:  dalawat@gmail.com}

\bigskip\bigskip

{\eightpoint {\bf Abstract}.  Starting from Stickelberger's congruence for the
  absolute discriminant of a number field, we ask a series of natural
  questions which ultimately lead to an orthogonality relation for the
  ramification filtration on $K(\root p\of{K^\times})$, where $K$ is any
  finite extension of $\Q_p$ containing a primitive $p$-th root of~$1$. An
  extensive historical survey of discriminants and primary numbers is
  included.  Among other things, we give a direct proof of Serre's mass
  formula in the case of quadratic extensions.  Incidentally, it is shown that
  every unit in a local field is the discriminant of some elliptic
  curve.\footnote{}{Keywords~: Stickelberger's congruence, discriminants,
    primary numbers, local fields, ramification filtration, cyclic extensions,
    Serre's mass formula, Kummer pairing, elliptic curves.}}

\bigskip\bigskip\bigskip

{{\it Die hier charakterisierte neue Theorie der algebraischen
Zahlen [...] scheint mir auch aus dem Grunde ein brauchbares Hilfsmittel
f{\"u}r arithmetische Untersuchungen zu sein, weil mit ihrer Hilfe Fragen der
Zahlentheorie vollst{\"a}ndig und einfach gel{\"o}st werden k{\"o}nnen, deren
Beantwortung mit den bisherigen Methoden entweder {\"u}berhaupt nicht gelang,
oder doch bedeutende Schwierigkeiten bereitete.\/}\hfill --- Kurt Hensel
\citer\hensel(p.~70).} 

\bigskip

Stickelberger showed that if the discriminant~$D$ of a degree-$n$ number field
$\Omega$ is not divisible by an odd prime~$p$, then $D$ becomes a square in
$\F_p^\times$ if and only if $n-m_p$ is even, where $m_p$ is the number of
places of $\Omega$ above~$p$.  If $2$ does not divide $D$, then $D$ is
$\equiv1\pmod8$ if $n-m_2$ is even, $\equiv5\pmod8$ if $n-m_2$ is odd.

Hensel showed that these global results are immediate consequences of purely
local ones.  We give the relative versions of Hensel's local results~: the
base field is no longer $\Qp$, but a finite extension $K$ thereof
(prop.~{15}).  Along the way, we also specify the $\F_l$-line in
$K^\times\!/K^{\times l}$ which corresponds, via Kummer theory, to the
unramified degree-$l$ extension of $K$ when $K^\times$ has an element of prime
order~$l$ (prop.~{16}).  Prop.~{15} turns out to be a somewhat sharper version
of a theorem of Fr{\"o}hlich (th.~{14}), which we had indeed set out to
sharpen.  In the presence of prop.~{30} and~{45}, prop.~{16} becomes the local
version of theorems of Hecke (th.~{53} and~{54}), which are immediate
corollaries, and which generalise a part of Hilbert's results (cf.~th.~{56}
and~{58}).


The other part of Hilbert's results deals with the valuation of the
discriminant~; a generalisation of this part can be deduced from results in
Hasse's {\it Klassenk{\"o}rperbericht\/} (cf.~th.~{60} and~{63}).
Stickelberger, Hilbert, Hecke and Hasse all four deal with number fields, but
the questions are purely local, and deserve a purely local proof, a proof
which Hensel could have given and which he did indeed give for Stickelberger's
theorem.  (Hensel did more~; read on to find out.)  This is what we do for the
other theorems.

An interesting local consequence of prop.~{16} is an explicit formula
(prop.~{17}) for the pairing $G\times D\rightarrow{}_p\mu$, where $D\subset
K^\times\!/K^{\times p}$ is the $\Fp$-line which corresponds to the unramified
degree-$p$ extension $L$ of $K$, and $G=\Gal(L|K)$.

An interesting global consequence of these theorems, apart from the
decomposition law in prime-degree kummerian extensions of number fields, is a
theoretical procedure for computing the relative discriminant of {\it any\/}
extension of number fields.  This also provides a test for an order in a
number field to be the maximal order.  Everything boils down to the
computation of the relative discriminant of a local kummerian extension of
degree equal to the residual characteristic, which is achieved in terms of the
ramification filtration and its relation to the natural filtration on the
multiplicative group~; see the final remark in Part~VII.

Our proofs require no more than a study of the filtration on the $\Z_p$-module
$U_1$ of principal units or {\it Einseinheiten\/} of $K$, especially of the
endomorphism $(\ )^p$ of raising to the exponent~$p$, which goes back to
Hensel.  Thus they have the appearance of a piece of late-nineteenth- or
early-twentieth-century arithmetic which fell into the twenty-first.

The paper consists of nine Parts, two of them of an historical nature.  Part~I
is a brief chronology of the work of Pellet, Brill, Stickelberger,
Vorono{\"\i}, Hensel, Schur, Herbrand and Fr{\"o}hlich on discriminants,
interspersed with a series of questions which lead to later developments.

\medbreak 

Part~II contains the statements of our results about discriminants of
unramified extensions, about unramified kummerian extensions of prime degree,
about the explicit $p$-tic (quadra{\it tic}, cu{\it bic}, quin{\it tic},
$\ldots$) character, about the filtration on $K^\times\!/K^{\times p}$, about
rings of integers, about discriminants of elliptic curves and about the
orthogonality relation.

\medbreak 

Part~VI is a brief review of the work of Fr{\"o}hlich, Hasse, Hecke and
Hilbert on primary numbers.  After the proof of prop.~{16}, it was natural to
try to extend Fr{\"o}hlich's results about rings of integers from quadratic
extensions to prime-degree kummerian extensions.  It was in seeking to do so
that we became aware of Hecke's global results, and, somewhat later, of those
of Hilbert.  Hecke's {\it S{\"a}tze\/}~118 (th.~{53}) and~119 (th.~{54}),
which constitute a generalisation of Hilbert's {\it Satz\/}~96 (th.~{58}) and
a part of {\it Satz\/}~148 (th.~{56}), follow from our local prop.~{16}.
Hensel makes a reentry on the scene at the end of this Part, closely followed
by Eisenstein.

The mathematical Parts can be read independently of the historical ones~; a
more detailed listing of the contents can be found in Part~II.

\medbreak

Part~III determines the structure of the multiplicative group $K^\times$ of a
local field $K$, following chapter~15 of Hasse's {\it Zahlentheorie}, itself
based upon Hensel's results.  Our treatment is more intrinsic, and some of our
proofs differ from theirs.  Using this, general properties of the
discriminant, and the compatibility with Artin-Schreier theory, we prove our
results about discriminants and $p$-primary numbers in Part~IV.

\medbreak 

In Part~V, we take a closer look at the filtration on
$K^\times\!/K^{\times p}$, which leads to an understanding of the precise
relationship of our results with the theorems of Fr{\"o}hlich and Hecke.  We
also give a few examples to show that some disparate results in the literature
follow from a systematic theory.

\medbreak

Part~VII deals with the computation of the discriminant of ramified
prime-degree cyclic extensions of local fields and the determination of their
rings of integers (cor.~{61}).  We determine the number of such extensions in
the kummerian case and, as a consequence, give a direct elementary proof, in
the case of quadratic extensions, of Serre's mass formula (lemma~67).

\medbreak

In Part~VIII, we introduce the discriminant of an elliptic curve~$E$ over a
local field $K$ as an element of $K^\times\!/\ogoth^{\times12}$~; this
definition was anticipated by J.~Silverman.  We show that --- in contrast to
the Stickelberger-Hensel condition (th.~{6}) and to our prop.~{15} --- every
element of $\ogoth^\times\!/\ogoth^{\times12}$ occurs as the discriminant of a
good-reduction elliptic curve (cor.~{71}).

\medbreak

Part~IX contains a few words about the genesis of these notes, and determines
the ramification filtration on $\Gal(M|K)$, where $M=K(\root p\of{K^\times})$
(the maximal abelian extension of exponent dividing~$p$) and $K|\Qp$ is a
finite extension containing a primitive $p$-th root of~$1$, in terms of the
filtration on $K^\times\!/K^{\times p}$.  The question was natural in the
light of prop.~{16}, which can be interpreted as saying that $\Gal(M|K)^n=\bar
U_{pe_1-n+1}^\perp$ for $n=1$, where the orthogonal is with respect to the
Kummer pairing, and where $\bar U_{pe_1}\subset K^\times\!/K^{\times p}$ is
the ``deepest'' $\Fp$-line~; the statement of orthogonality for $n>1$ is taken
from \citer\delcorso(), where the special case
$K=F(\!\root{p-1}\of{F^\times})$ for some (finite) extension $F|\Q_p$ is
treated.  Our proof in the general case is simpler and more conceptual.


\bigskip \centerline{***} 

\bigbreak
\centerline{\bf I. A brief history of discriminants}

\bigskip
{\bf Stickelberger (1897).} At the first international
Congress of mathematicians, held at Z{\"u}rich, L.~Stickelberger
proved some new properties of the discriminant ({\it Grundzahl\/}) of an
algebraic number field $\Omega$ of degree~$n$.  We quote the two theorems
relevant to us.

\th THEOREM {1} (\citer\stick(p.~186))
\enonce 
Die Diskriminante des K{\"o}rpers $\Omega$ ist durch die Primzahl $p$ nicht
teilber, wenn $p$ ein Produkt von lauter verschiedenen Primidealen in $\Omega$
ist~; zugleich ist sie, wenn $p$ ungerade, quadratischer Rest oder Nichtrest
von $p$, je nachdem die Anzahl der in $p$ aufgehenden Primideale von geradem
Grade eine gerade oder ungerade ist, oder je nachdem die Anzahl aller
Primfaktoren von $p$ dem Grade des K{\"o}rpers kongruent ist nach dem Modul
Zwei oder nicht.
\endth
(If a prime $p$ is unramified in $\Omega$, then it does not divide the
discriminant $D$ of $\Omega$~; if such a $p$ is odd, then $D$ is a quadratic
residue $\pmod p$ or not according as the number of even-degree places above
$p$ is even or odd.)

The {\it etwas m{\"u}hsamer\/} prime $2$ is treated in the last result of
the paper. 

\th THEOREM {2} (\citer\stick(p.~192))
\enonce
Die Diskriminante des Zahlk{\"o}rpers $\Omega$, in dem $2$ ein Produkt von $m$ 
verschiedenen Primidealen $\pgoth_1$, $\pgoth_2$, $\ldots$, $\pgoth_m$ der
Grade $f_1$, $f_2$, $\ldots$, $f_m$ ist, ist von der Form $8q+1$ oder $8q+5$,
je nachdem unter jenen Primfaktoren solche geraden Grades in gerader oder
ungerader Zahl vorkommen~; in Zeichen ist
$$
{D-1\over4}\equiv\sum(f-1)=n-m\pmod2
$$
oder
$$
(-1)^{D-1\over4}=(-1)^{n-m}.
$$
\endth
(If the prime~$2$ is unramified in $\Omega$, then the discriminant of $\Omega$
is $\equiv1\pmod8$ or $\equiv5\pmod8$ according as the number of even-degree
places above $2$ is even or odd.)

\bigbreak {\bf Vorono{\"\i} (1905).} G.~Vorono{\"\i}, unaware of
Stickelberger's results, rediscovers th.~{1} and gives a different proof at the
third international Congress at Heidelberg in 1904.  He takes an irreducible
polynomial $F(x)\in\Z[x]$ of degree~$n$ whose discriminant is not divisible by
a prime $p$ and factors it as $F=\varphi_1\varphi_2\cdots\varphi_\nu$ into
irreducible polynomials in $\F_p[x]$.  Implicitly assuming that $p$ is odd,
his version of th.~{1} is the following.

\th THEOREM {3} (\citer\voronoi(p.~186))
\enonce
Le nombre $\nu$ des facteurs irr{\'e}ductibles de la fonction $F(x)$ par
rapport au module $p$ v{\'e}rifie l'{\'e}quation
$$
\left({D\over p}\right)=(-1)^{n-\nu},
$$
o{\`u} $\displaystyle\left({D\over p}\right)$ est le symbole de Legendre.
\endth
(The number $\nu$ of irreducible factors of $F$ modulo $p$ satisfies the
displayed equation involving the Legendre symbol.)

Th.~1 was subsequently rediscovered by Th.~Skolem \citer\skolem() and R.~Swan
\citer\swan().  Indeed, it had been anticipated by A.~Pellet~:

\th THEOREM {4} (\citer\pellet(p.~1071))
\enonce
Soit\/ $\Delta$ le produit des carr{\'e}s des diff{\'e}rences des racines
d'une congruence\/ $f(x)\equiv0\pmod p$ n'ayant pas de racines {\'e}gales~;
$\Delta$ est non-r{\'e}sidu quadratique $\pmod p$, si $f(x)$ admet un nombre
impair de facteurs irr{\'e}ductibles de degr{\'e} pair~; $\Delta$ est, au
contraire, r{\'e}sidu quadratique, si $f(x)$ n'admet pas de facteurs
irr{\'e}ductibles de degr{\'e} pair ou en admet un nombre pair.
\endth
(The discriminant $\Delta$ of a separable polynomial $f\in\Fp[x]$ is in
$\F_p^{\times2}$ if and only if the number of even-degree irreducible factors
of $f$ is even.)

\medskip\centerline{***}\medskip

With these beginnings, we shall ask a series of questions which lead naturally
to many twentieth-century results about discriminants and related topics.  The
idea is thus to see the new results in an old light.

Our first questions is~: {\it Aren't these theorems global manifestations
  of purely local results~?}

\bigbreak

{\bf Hensel (1905).}  This was first realised by K.~Hensel, who proved the
local result at the finite places.

\th THEOREM {5} (\citer\hensel(p.~78))
\enonce
Ist $f(x)$ f{\"u}r den Bereich von $p$ irreduktibel, und ist $p$ kein
Diskriminantenteiler des zugeh{\"o}rigen K{\"o}rpers $K(x)$, so ist die
Diskriminante einer jeden Gleichung dieses K{\"o}rpers quadratischer Rest oder
Nichtrest zu $p$, je nachdem $n$ ungerade oder gerade ist~; est ist also
stets~:
$$
\left({D\over p}\right)=(-1)^{n-1}.
$$
\endth
(If $f\in\Qp[x]$ is irreducible, and if the odd prime $p$ does not divide the
discriminant $D$ of $\Qp[x]/f$, then $D$ is a quadratic residue $\pmod p$ or
not according as $n$ is odd or even.)

He says that this formula remains valid for $p=2$ if we define the Legendre
symbol in this case as
$$
\left({D\over 2}\right)=i^{D-1\over2}=(-1)^{D-1\over4}.
$$
\vskip-\lastskip
\th THEOREM {6} (\citer\hensel(p.~79))
\enonce
Ist $f(x)$ f{\"u}r den Bereich der Primzahl $2$ irreduktibel und ist $2$ nicht
in der K{\"o}rperdiskriminante enthalten, so ist jede Diskriminante dieses
K{\"o}rpers von der Form $8\nu+1$ oder $8\nu+5$, je nachdem der Grad jenes
K{\"o}rpers ungerade oder gerade ist.  Est gibt keine Diskriminante von der
Form $8\nu+3$ oder $8\nu+7$.
\endth
(If $f\in\Q_2[x]$ is irreducible, and if $2$ does not divide the discriminant
$D$ of $\Q_2[x]/f$, then $D$ is $\equiv1$ or $\equiv5\pmod8$ according
as $n$ is odd or even.  No discriminant is $\equiv3$ or $\equiv7\pmod8$.)

Pellet, Hensel-Mirimanoff, Laskar, Swan, and Barrucand-Laubie have used these
theorems to give a new proof of the law of quadratic reciprocity~; they
continue to inspire current research~: see, for example, \citer\zahidi().  

To these results about the reduction of the discriminant (modulo an odd
prime~$p$, or modulo~$8$) must be added information about its sign, which goes
back to A.~Brill.  He is working with polynomials with real coefficients and
finds that (cf.~prop.~{9})~:

\th THEOREM {7} (\citer\brill(p.~87))
\enonce
Das Vorzeichen der Discriminante einer Gleichung --- lauter verschiedene
Wurzeln vorausgesetzt --- ist negativ, wenn die Anzahl der complexen
Wurzelpaaren eine Ungerade ist, positiv, wenn diese Zahl gerade ist.
\endth
(The sign of the discriminant of a separable real polynomial is negative if
the number of pairs of complex conjugate roots is odd, positive if this number
is even.)

Unaware of Brill's local th.~{7}, Hensel considers the number field $K(x)$
obtained by adjoining a root $x$ of an irreducible polynomial $f\in\Q[x]$ of
degree~$n$, with factorisation $f=f_1f_2\cdots f_h$ into irreducible
polynomials in $\R[x]$ (necessarily linear or quadratic by a theorem of
C.~Gauss, as he reminds us~: {\it Gau\ss\ hat zuerst streng
bewiesen,}$\ldots$).

\th THEOREM {8} (\citer\hensel(p.~70))
\enonce
Die Basisdiskriminanten eines K{\"o}rpers $K(x)$ sind s{\"a}mtlich positiv
oder s{\"a}mtlich negativ, je nachdem $n-h$ gerade oder ungerade ist, d.~h.~es
ist 
$$
\mathop{\rm sgn}D=(-1)^{n-h}.
$$
\endth 
(The sign of the discriminant $D$ of $K(x)$ is given by the above equation.)

By multiplicativity, Brill's th.~{7} comes down to the following prop., which
also implies Hensel's global th.~{8}~:

\th PROPOSITION {9}
\enonce
The discriminant of a finite extension\/ $K\,|\,\R$ is\/ $(-1)^{n-1}$ in\/ 
$\R^\times\!/\R^{\times2}$, where $n=[K:\R]$. 
\endth
\proof  
Only the case $K=\C$ needs to be considered. Computing the discriminant
$d_{\C|\R}$ using the standard basis $1,i$, we get 
$d_{\C|\R}=-4$, which is the same as $-1$ in $\R^\times\!/\R^{\times2}$.
We could have equally well computed the discriminant (``$b^2-4ac$'') of
$T^2+1$.

Our next question is~: {\it What is the valuation of the discriminant when the
  extension $K\,|\,\Q_p$ is ramified~?}

This leads to a host of results due to R.~Dedekind, D.~Hilbert, E.~Artin,
J.~Herbrand, etc., some of which are taken up in Part~VII.  Let us content
ourselves here with a modest corollary of a result \footnote{$({}^1)$}{{\it
Ist\/ $\pgoth$ ein beliebieges Primideal, $p$ der durch\/ $\pgoth$ theilbare
rationale Primzahl, und\/ $\pgoth^e$ die h{\"o}chste in $p$ aufgehende Potenz
von\/ $\pgoth$, so ist das Grundideal $\dgoth$ allemal theilbar durch\/
$\pgoth^{e-1}$~; ist ferner der Exponent\/ $e$ nicht theilbar durch\/ $p$, so
ist $\dgoth$ nicht theilbar durch $\pgoth^e$~; ist aber $e$ theilbar
durch\/~$p$, so ist $\dgoth$ theilbar durch $\pgoth^e$ und vielleicht durch
noch h{\"o}heren Potenzen von~$\pgoth$.}  \citer\dedekind(p.~52) (If $\pgoth$
is a prime ideal, $p$ the rational prime divisible by $\pgoth$ and $\pgoth^e$
the highest power of $\pgoth$ dividing $p$, then the different $\dgoth$ is
divisible by $\pgoth^{e-1}$~; if $p$ does not divide $e$, then $\pgoth^e$ does
not divide $\dgoth$~; if $p$ divides $e$, then $\pgoth^e$ and possibly a
higher power of $\pgoth$ divides $\dgoth$.)} of Dedekind.

\th THEOREM {10} (\citer\dedekind(p.~54))
\enonce
Ist aber\/ $p=2$, also der Exponent\/ $e$ theilbar durch\/ $p$, so ist\/
$\dgoth$ mindestens durch\/ $\pgoth^2$, und folglich\/ $D$ mindestens durch\/
$4$ theilber.
\endth
(If $p=2$ and the ramification index $e$ is even, then $\pgoth^2$ divides
the different $\dgoth$ and consequently $4$ divides the discriminant $D$.)

\bigbreak 
{\bf Schur (1929).}  A simple proof of the conjunction of this fact
and of a somewhat less precise form of th.~{2} was given by I.~Schur~;
regrettably, this is the version most commonly cited these days.

\th THEOREM {11} (\citer\schur(p.~29))
\enonce
Die Diskriminante\/ $D$ eines algebraischen Zahlk{\"o}rpers is stets
kongruent\/ $0$ oder\/ $1$ nach dem Modul\/ $4$.
\endth
(The discriminant $D$ of a number field is always $\equiv0$ or
$\equiv1\pmod4$.)  

Our next question is~: {\it What are the relative versions of these results~?} 

\bigbreak

{\bf Herbrand (1932).}
Indeed, J.~Herbrand asked himself this question before his tragic death in a
mountaineering accident at the young age of~23.  In a paper which he could
write only partially, and which was completed by C.~Chevalley based upon his
rough notes ({\it brouillons\/}), he proves the following theorem.

\th THEOREM {12} (\citer\herbrand(p.~105))
\enonce
$\vartheta$ {\'e}tant le discriminant par rapport au corps\/ $k$ d'un surcorps
relativement metacyclique $K$ de degr{\'e} relatif\/ $N$, on a
$\vartheta=\agoth^2(\alpha)$, o{\`u} $\alpha$ est un nombre tel que~:

a. $\alpha\equiv1\pmod\bgoth$, $\bgoth$ {\'e}tant le plus grand id{\'e}al
divisant $4$ et premier {\`a} $\vartheta$.

b. Le conjugu{\'e} de $\alpha$ dans un corps r{\'e}el conjugu{\'e} de $k$
n'est n{\'e}gatif que si le conjugu{\'e} correspondant de\/ $K$ est imaginaire
et si\/ $N\equiv2\pmod4$.
\endth
($\vartheta$ being the discriminant of a metacyclic extension $K|k$ of
degree~$N$, one has $\vartheta=\agoth^2(\alpha)$, where $\alpha$ is a number
such that $\alpha\equiv1\pmod\bgoth$, $\bgoth$ being the greatest ideal
dividing $4$ and prime to $\vartheta$.  Moreover, the conjugate of $\alpha$ in
a real field is not negative unless the corresponding conjugate of $K$ is
imaginary and $N\equiv2\pmod4$.)

\noindent
As a corollary, he gets  a theorem of Hecke~:

{\it Le discriminant d'un corps alg{\'e}brique, par rapport {\`a} un
  sous-corps, est le dans carr{\'e} d'une classe de ce sous-corps} (The
  discriminant of a number field, with respect to a subfield, is the square of
  a class in that subfield),

\noindent
and, taking $k=\Q$, the theorem of Stickelberger-Schur~:

{\it Le discriminant d'un corps alg{\'e}brique est congru {\`a} $0$ ou {\`a}
  $1\pmod4$} (The discriminant of a number field is always $\equiv0$ or
  $\equiv1\pmod4$),

\noindent
in the case of metacyclic extensions (resp.~number fields).  He says that as
these two theorems are true without restriction to the metacyclic --- or even
to the galoisian --- case, it is probable that his theorem is also true
without the hypothesis that the extension $K|k$ be metacyclic.

\bigbreak {\bf Fr{\"o}hlich (1960).}  This conjecture was proved some thirty
years later by A.~Fr{\"o}hlich when he introduced the id{\'e}lic discriminant.
For an extension $\Lambda|K$ of number fields, the discriminant
$\dgoth_{\Lambda|K}$ is an element of the restricted direct product of the
various $K_{p}^\times\!/\ogoth_{p}^{\times2}$, where $p$ runs through the
places of $K$.  The classical discriminant of $\Lambda|K$ becomes the
(integral) ideal $(\dgoth_{\Lambda|K})$ associated to the id{\'e}lic
discriminant.

\th THEOREM {13} (\citer\frohlich(p.~28))
\enonce
The ideal $(\dgoth_{\Lambda|K})$ ($\Lambda|K$ normal) can be written in the
form $\agoth^2(\alpha)$, $\alpha\in K^\times$ where

{\rm (i)} $\alpha\equiv1\pmod\bgoth$, $\bgoth$ the greatest divisor of $4$
prime to $(\dgoth_{\Lambda|K})$.

{\rm (ii)} $\alpha_p>0$ for each real prime divisor $p$, except when
$\Lambda_p$ is a direct sum of copies of the field of complex numbers and 
$(\Lambda:K)\equiv2\pmod4$. 
\endth

He remarks that (i) is also true when $\Lambda|K$ is not normal, while in the
place of the congruence in (ii) we have $r_2\equiv1\pmod2$.  His relative
version of the Stickelberger-Schur theorem (th.~{11}) is the following.

\th THEOREM {14} (\citer\frohlich(p.~23))
\enonce
Every discriminant is a quadratic residue $\mathop{\rm mod} 4$.
\endth

The proof is by reduction to the local case, where Schur's proof applies
almost without change.  (This result was anticipated by K. Dalen
\citer\dalen(p.~125) and rediscovered by R. Swan \citer\swan(p.~1100).)
  
But in the formulation ``{\it Every discriminant is a quadratic residue
$\mathop{\rm mod} 4$}'', the local version is contentless for odd residual
characteristics, because $\ogoth/4\ogoth=0$~: one does not recover Hensel's
th.~{6}.  When the residual characteristic is $2$, it leaves open the question
as to which squares in $\ogoth/4\ogoth$ are discriminants, and which units the
discriminants of even-degree unramified extensions~; in the absolute case
$K=\Q_2$, th.~{6} gives precise answers to these questions.  Moreover, it is
not very aesthetic to lift elements (of positive valuation) of the
multiplicative group $K^\times\!/\ogoth^{\times2}$ to the additive group
$\ogoth$ with the purpose of going modulo~$4\ogoth$.

In the next section, we give our formulation of the relative version and show
how it can be interpreted as a characterisation of ``$2$-primary'' numbers,
leading to the question of a characterisation of ``$p$-primary'' numbers ---
to use a piece of terminology we learnt much later --- for every prime $p$,
answered in prop.~{16}.

\bigbreak
\centerline{\bf II.  The main propositions}
\bigskip

{\bf The relative version.}
The correct relative version in the local case is not far to seek.  The
discriminant $d_{L|K}$ is an element of $K^\times\!/\ogoth^{\times2}$, a group
which comes equipped with a filtration
$$
\cdots\subset\bar U_n\subset\cdots\subset\bar
U_1\subset\ogoth^\times\!/\ogoth^{\times2}\subset 
K^\times\!/\ogoth^{\times2}\leqno{(1)}
$$
by vector $\F_2$-spaces, deduced from the filtration
$$
\cdots\subset U_n\subset\cdots\subset
U_1\subset\ogoth^\times\subset 
K^\times,\leqno{(2)}
$$
where
$U_n=\Ker(\ogoth^\times\rightarrow(\ogoth/\pgoth^n)^\times)$ ($n>0$).  
The filtration $(1)$ is finite~; indeed $\bar U_{2e+1}=\one$, and $\bar
U_{2e}=\{1,u\}$ is an $\F_2$-line, where $e$ is the (absolute) ramification
index of $K$ if $p=2$, and $e=0$ if $p\neq2$ (cf.~prop.~{33}).  We have adopted
the convention that $U_0=\ogoth^\times$.

\th PROPOSITION {15} 
\enonce
Let\/ $K$ be a finite extension of\/ $\Q_p$ ($p$ prime) and let\/ $L$ be a
finite unramified extension of\/ $K$, of degree\/ $r=[L:K]$.  Then the
discriminant\/ $d_{L|K}$ belongs to the\/ $\F_2$-line\/ $\bar
U_{2e}=\{1,u\}$~; one has\/ $d_{L|K}=1$ if\/ $r$ is odd, $d_{L|K}=u$ if\/ $r$
is even.  When\/ $p=2$, none of the other\/ $2^{d+1}-2$ ($d=[K:\Qp]$) elements
of\/ $\ogoth^\times\!/\ogoth^{\times2}$ is a discriminant.
\endth

The proof is given in Part~III.  Hensel could have easily proved prop.~{15},
because all we need for the proof are Hensel's results on the structure of the
multiplicative group of $K$, and general properties of the discriminant in a
tower of extensions, which go back to Hilbert.  We wish to argue that if he
had done it, and if this version had become established instead of the
Schur-Fr{\"o}hlich version (th.~{11}, th.~{14}), the history of mathematics
would have been different in a few respects.

The first thing to notice is that, when $p\neq2$, the discriminant of the
residual extension of $L|K$ is the reduction of the discriminant via the
isomorphism $\bar U_{2e}\rightarrow k^\times\!/k^{\times2}$.  What about the
case $p=2$~?  In this case, there is an isomorphism $\bar U_{2e}\rightarrow
k/\wp(k)$, where $\wp$ is the $\F_2$-linear endomorphism $x\mapsto x^2-x$ of
$k$.  

Our next question therefore is~: {\it Is there a way of defining discriminants
  --- with values in $k/\wp(k)$ --- 
  when $k$ is an  extension of\/ $\F_2$~?}

\bigbreak 

{\bf Discriminants in characteristic 2}.  There is indeed one~; it was found
by E.~Berlekamp \citer\berlekamp() in 1976, who seems to have been motivated
by coding theory.  He does notice (p.~326) the analogy between his definition
and Stickelberger's th.~{1}, but fails to mention the even stronger analogy
with th.~{2}, or its substantial identity with Hensel's th.~{6}, which can be
construed as implying that the discriminant of a finite extension of $\F_2$ is
$0$ (trivial) if the degree is odd, $1$ (not trivial) if the degree is even.
When the base field $\F_2$ is replaced by any finite extension $k$ thereof,
prop.~{15} implies the analogous result~: the discriminant of an odd-degree
(resp.~even-degree) extension of $k$ is $0$ (resp.~$1$) in the additive
$2$-element group $k/\wp(k)=\F_2$.

The interpretation of the charcteristic-$2$ discriminant as the reduction of a
\hbox{characteristic-$0$} discriminant was found by A.~Wardsworth
\citer\wadsworth() in 1985, eighty years after Hensel's absolute local
versions~(th.~{5}, th.~{6}).

\bigbreak 

{\bf Unramified kummerian extensions.}  There is {\it another reading\/} of
prop.~{15}.  It can be viewed as specifying the $\F_2$-line in
$K^\times\!/K^{\times2}$ which gives us the unramified quadratic extension of
$K$ upon adjoining square roots.

If the local field $K$ contains a primitive $l$-th root of $1$ for some prime
$l$, then degree-$l$ cyclic extensions of $K$ correspond to $\F_l$-lines in
$K^\times\!/K^{\times l}$ (``Kummer theory'').

Our next question is~: {\it Which line in the\/ $\F_l$-space\/
  $K^\times\!/K^{\times l}$ gives us the unramified $(\Z/l\Z)$-extension of
  $K$~?} 

If $l=p$, let ${e_1}$ stand for the ramification index of $K|\Qp$ divided by
$p-1$ (the quotient ${e_1}$ is an integer~; cf.~prop.~{25})~; if $l\neq p$, put
${e_1}=0$.

It can be shown that the induced filtration on
$K^\times\!/K^{\times l}$ by $\F_l$-spaces
$$
\cdots\subset\bar U_n\subset\cdots\subset\bar
U_1\subset\ogoth^\times\!/\ogoth^{\times l}\subset 
K^\times\!/K^{\times l}\leqno{(3)}
$$
has $\bar U_{p{e_1}+1}=\one$, and that $\dim_{\F_l}\bar U_{p{e_1}}=1$
(see Part~III).

\th PROPOSITION {16}
\enonce
The\/ $\F_l$-line in\/ $K^\times\!/K^{\times l}$ which gives the unramified\/ 
$(\Z/l\Z)$-extension of\/ $K$ upon adjoining\/ $l$-th roots is\/ $\bar
U_{p{e_1}}$. 
\endth

The proof is to be found in Part~III.  Notice that $\bar U_{pe_1}=\bar
U_{le_1}$ even when $l\neq p$, for then $e_1=0$.

\bigbreak 
{\bf Degree-$p$ cyclic extensions in characteristic $p$.} When $l=p$,
the choice of an element $\zeta\in K^\times$ of order~$p$ gives us an
isomorphism $\bar U_{p{e_1}}\rightarrow k/\wp(k)$, where $\wp:k\rightarrow k$
is the $\F_p$-linear map $x\mapsto x^p-x$ (cf.~the discussion after
prop.~{33}).  This leads us to the next question.

We ask~: {\it Do $(\Z/p\Z)$-extensions of a characteristic-$p$ field\/ $k$
  correspond to $\F_p$-lines in $k/\wp(k)$~?}

They indeed do, as was discovered by E.~Artin and O.~Schreier
\citer\artinschreier(), who were led to their result by an entirely different
route (maximally ordered fields).  Neither they, nor E.~Witt \citer\witt(),
make any connection\footnote{$({}^2)$}{According to Prof.~Peter Roquette
  \citer\roquette(), Artin made this connection in a letter to Hasse in 1927,
  with the words {\it Ich entdeckte, da\ss\ hier ein alter Bekannter von mir
    vorlag}$\ldots$~; see \citer\freiroquette().  This connection also appears
  independently and implicitly in Hasse \citer\hassepotenz(p.~234) as
$$
\left(\alpha\over{\goth l}\right)_{\!l}=
\zeta_0^{S_{\goth l}({\alpha-1\over l\lambda_0})},\leqno{(13)}
$$
which is essentially the case $a=1$ of our (purely local) prop.~{17}.  See
also his {\it Klassenk{\"o}rperbericht}, Teil~II, \S~17, 3, IV, to which we
don't have access.  Hasse derives $(13)$ from Artin's general reciprocity
law~; not even its local version is needed for the proof of prop.~{17}.} with
Hensel's determination of the filtration on $K^\times\!/K^{\times p}$, which
preceded their results by more than twenty years.  The connection between the
two theories was fully understood only in the eighties, in the works of
T. Sekiguchi, F. Oort \& N. Suwa \citer\sekiguchi() and W. Waterhouse
\citer\waterhouse().  It is no accident that last-named author published a
paper \citer\waterdisc() on the cohomological interpretation of the
discriminant in the same year~; as we have seen, the two topics are not
unrelated.  As it happens, we shall make a modest use of Artin-Schreier theory
in our proofs.

\bigbreak 
{\bf The Kummer pairing.}  
Recall that Kummer theory provides not only a bijection between the set of
$\Fp$-lines\/ $D\subset K^\times\!/K^{\times p}$ and the set of degree-$p$
cyclic extensions\/ $L$ of\/ $K$, but also a pairing $\Gal(L|K)\times
D\rightarrow{}_p\mu$ when $D$ and $L$ correspond to each other~: 
$$
L=K(\root p\of D),\quad 
D=\Ker(K^\times\!/K^{\times p}\rightarrow L^\times\!/L^{\times p}).
$$ 

Our next question is~: {\it Can we give an explicit description of the
pairing\/ $\langle\ ,\ \rangle:G\times\bar U_{pe_1}\rightarrow{}_p\mu$,
where\/ $G=\Gal(L|K)$ and $L$ is the degree-$p$ unramified extension of\/
$K$~?}

The group $G$ has a canonical generator $\varphi$ (``Frobenius'')~: the unique
element $\varphi$ such that $\varphi(\alpha)\equiv \alpha^q\pmod{\pgoth_L}$
for every $\alpha\in\ogoth_L$, where $q=\Card k$ and $k$ is the residue field
of $K$.  We shall see in Part~III that the choice of a generator
$\zeta\in{}_p\mu$ leads to a specific isomorphism $\bar\eta\mapsto\hat c:\bar
U_{pe_1}\rightarrow k/\wp(k)$, where $c=(\eta-1)/p(\zeta-1)$, which belongs to
$\ogoth$, and $\hat c$ is its image in $k/\wp(k)$.  This can be composed with
the trace map $S_{k|\Fp}:k/\wp(k)\rightarrow\Fp$, which is also an
isomorphism.

\th PROPOSITION 17
\enonce
Choose and fix a generator $\zeta\in{}_p\mu$.  For\/ $a\in\Z/p\Z$ and\/
$\eta\in U_{pe_1}$, we have
$\langle\varphi^a,\bar\eta\rangle=\zeta^{a.S_{k|\Fp}(\hat c)}$.
\endth
The proof is to be found in Part~IV (prop.~{38}), where the implication
\citer\serre(p.~230) for the pairing $K^\times\times\bar
U_{pe_1}\rightarrow{}_p\mu$ coming from the reciprocity map
$K^\times\rightarrow G$ is also mentioned.

Another interesting consequence of the characterisation of the unramified
degree-$p$ kummerian extension (prop.~{16}) is a generalisation of the
Pellet-Vorono{\"\i} theorem (th.~{3} and~{4}) to arbitrary finite fields~$k$:
the discriminant of a separable polynomial $f\in k[T]$ is trivial precisely
when the number of even-degree irreducible factors of $f$ is even (cor.~{41}).

\bigbreak 

{\bf The filtration on $K^\times\!/K^{\times p}$.}  In Part~V, we will give
three different characterisations of the filtration $(\bar U_n)_{n>0}$ on
$K^\times\!/K^{\times p}$ for a finite extension $K$ of $\Qp$ and compute the
$\Fp$-dimension of the quotients $\bar U_n/\bar U_{n+1}$.  The main  results
are the following two propositions.

\th PROPOSITION {18}
\enonce
Let\/ $\zeta\in\bar K^\times$ be an element of order\/ $p$.  Then
$$
\dim_\Fp\bar U_{n}/\bar U_{n+1}=\cases{
0&if\/ $n>p{e_1}$,\cr
1&if\/ $n=p{e_1}$ and\/ $\zeta\in K^\times$,\cr
0&if\/ $n=p{e_1}$ and\/ $\zeta\notin K^\times$,\cr
0&if\/ $n<p{e_1}$ and\/ $p|n$,\cr
f&otherwise.\cr}
$$
\endth

\th PROPOSITION {19}
\enonce
For\/ $x\in\units$ and\/ $n>0$, let $\bar x$ be its image in\/
  $\units\!/\ogoth^{\times p}$ and\/ $\hat x$ the image in\/
  $(\ogoth/\pgoth^{n})^\times$.  Then
$$
\bar x\in\bar U_{n}
\ \Longleftrightarrow\ 
\hat x\in(\ogoth/\pgoth^{n})^{\times p}.
$$
\endth 
It is this result which allows us to deduce the theorems of Hecke and Hilbert 
from our local results.  In this Part, we also work out a number of explicit
examples.  

\bigbreak 

{\bf Rings of integers in kummerian extensions.}  In Part~VII, we determine,
following Hasse, the valuation $v_{K}(d_{L|K})$ of the discriminant $d_{L|K}$
of cyclic degree-$p$ extensions $L|K$ of local fields (prop.~{60}).  The proof
allows us to determine the ring of integers $\ogoth_L$ explicitly (prop.~{61}).
One can also compute the number of such $L$ with a given $v_{K}(d_{L|K})$
(prop.~{66}).

We also explain how this solves the global problem of determining the relative
discriminant of an extension of number fields.

\bigbreak 
{\bf Discriminants of elliptic curves over local fields.}  In view of the
definition of the discriminant $d_{L|K}$ of a finite extension $L|K$ of local
fields as an element of $K^\times\!/\ogoth^{\times 2}$, it is natural to
define the discriminant of an elliptic curve $E|K$ as an element 
$d_{E|K}\in K^\times\!/\ogoth^{\times 12}$ (see Part VIII).  

We may ask : {\it Suppose that\/ $E|K$ has good reduction.  Does\/
  $d_{E|K}\in\ogoth^\times\!/\ogoth^{\times 12}$ have to satisfy some
  congruence~?}

In contrast to th.~{6} and prop.~{15}, every element of
$\ogoth^\times\!/\ogoth^{\times 12}$ occurs as $d_{E|K}$ for some
good-reduction elliptic curve $E|K$ (prop.~{71}).  As corollary, every element
of $k^\times\!/k^{\times12}$ occurs as the discriminant of some elliptic
$k$-curve over any finite field $k$ (cor.~{72}).

\bigbreak 

{\bf The orthogonality relation.}  In Part~IX, we derive the orthogonality
relation $\Gal(M|K)^n=\bar U_{pe_1-n+1}^\perp$ which can be thought of as a
generalisation of prop.~{16}.  It determines the ramification filtration on
the maximal exponent-$p$ kummerian extension $M=K(\!\root p\of{K^\times})$ of
a finite extension $K|\Q_p$ (containing a primitive $p$-th root $\zeta$
of~$1$).  Upon taking $K=F(\zeta)$, it can be used to derive the ramification
filtration on the maximal exponent-$p$ abelian extension of any finite
extension $F$ of $\Q_p$.

\bigbreak
\centerline{\bf III.  The multiplicative group of a local field}
\bigskip

In chapter 15 of his {\it Zahlentheorie} \citer\hasse(), H.~Hasse studies the
multiplicative group $K^\times$ of a finite extension $K$ of $\Q_p$ ($p$
prime), whose determination goes back to Hensel \citer\henselmult().  In this
Part, we give a brief account of the results, not all of which are needed for
what follows.  Some of our proofs (for example in \S1 and \S4) are different
from those of Hensel and Hasse.  A part of \S3 can be also be found in
\citer\fesvost(I,\S5).

We first study the analogous local field $k\series{T}$, where $k$ is a finite
extension of $\F_p$.  In both cases, the multiplicative group comes equipped
with a decreasing sequence of subgroups $(U_n)_{n}$ which are $\Z_p$-modules
for $n>0$.  The result for $k\series{T}$ states that the $\Z_p$-module $U_1$
of {\it Einseinheiten\/} is not finitely generated, so the filtration on
$U_1/U_1^p$ is not of finite length.

By contrast, the group of $1$-units $U_1$ in $K$ is finitely generated as a
$\Z_p$-module (cor.~{32}) and the filtration on $K^\times\!/K^{\times p}$ is
of finite length.  There is a criterion for $K^\times$ to contain an element
of order~$p$ (prop.~{25}).  For an element $\zeta$ of order $p^\alpha$, the
precise {\it level} --- the integer $n$ such that $\zeta\in U_n$ but
$\zeta\notin U_{n+1}$ --- is known (prop.~{26}).  We study the
raising-to-the-exponent-$p$ map $(\ )^p$ and show that $U_n^p\subset
U_{\lambda(n)}$, where $\lambda(n)=\inf(pn, n+e)$, and
$e=(v(K^\times):v(\Qp^\times))$ is the ramification index (prop.~{27}).  Next,
we show that the induced map $\rho_n:U_n/U_{n+1}\rightarrow
U_{\lambda(n)}/U_{\lambda(n)+1}$ is always an isomorphism except when
$K^\times$ has an element of order~$p$ and $n=e/(p-1)$, in which case both
$\Ker\rho_n$ and $\Coker\rho_n$ are cyclic of order~$p$~; we determine these
groups explicitly (prop.~{33}).  Finally, we determine the structure of the
$p$-groups $U_1/U_n$ for sufficiently large $n$ (prop.~{34}).

\bigbreak 
\centerline{\bf 1. The multiplicative group $k\series{T}^\times$} 
\medskip 

Let $k$ be a finite extension of $\Fp$ and let $K=k\series{T}$.  For every
$n>0$, let $U_n$ be the kernel of the reduction map from $k[[T]]^\times$ to
$(k[[T]]/T^nk[[T]])^\times$~; in particular, $U_1=\Ker(k[[T]]^\times\to
k^\times)$.  The $U_n$ are $\Z_p$-modules, because they are commutative
pro-$p$-groups.

\th PROPOSITION {20}
\enonce
The $\Z_p$-module $U_1=\Ker(k[[T]]^\times\to k^\times)$ is not
finitely generated.
\endth

It is sufficient to show that $(U_1:U_1^p)$ is not finite.  Supposing that it
is, we shall get a contradiction.

We have $U_n^p\subset U_{pn}$ for every $n$, because $(1+aT^n)^p=1+a^pT^{pn}$
for every $a\in k[[T]]$.  The inclusions $U_n^p\subset U_{pn}\subset U_n$
imply that
$$
(U_n^{\phantom{p}}:U_n^p)\ge(U_n:U_{pn})=q^{pn-n}\qquad (q=\Card k),
$$
because $(U_i:U_{i+1})=q$, as $(1+aT^i)\mapsto a$ induces an isomorphism
$U_i/U_{i+1}\rightarrow k$ for every $i>0$.

We also have the inclusions $U_n^p\subset U_n\subset U_1$ and
$U_n^p\subset U_1^p\subset U_1$ which allow us to compute the index
$(U_1:U_n^p)$ in two different ways.  Comparing them, and using the
fact that $(U_1:U_n)=(U_1^p:U_n^p)$ (for which the absence of torsion
--- the only $p$-th root of $1$ in a field of characteristic $p$ is
$1$ --- is needed), we get the equality $(U_1:U_1^p)=(U_n:U_n^p)$,
which is larger than $q^{pn-n}$.  But we can take $n$ as large as we
please, so $(U_1:U_1^p)$ cannot be finite~!

\th COROLLARY {21}
\enonce
The\/ $\F_p$-space\/ $K^\times\!/K^{\times p}$ is infinite.
\endth

\medskip\centerline{***}\medskip

We fix the notation for the rest of Part~III~: $K$ is a finite extension of
$\Q_p$, $v:K^\times\rightarrow\Z$ is its surjective valuation,
$\ogoth=v^{-1}([0,+\infty])$ is its ring of integers, with unique maximal
ideal $\pgoth=v^{-1}(]0,+\infty])$ and residue field $k=\ogoth/\pgoth$.  The
units $\ogoth^\times$ will also be denoted $U_0$~; for $n>0$, we put
$U_n=1+\pgoth^n$.  Denote by $d=[K:\Q_p]$ the degree of $K$, by
$e=(v(K^\times):v(\Qp^\times))$ the ramification index, and by $f=[k:\F_p]$
the residual degree~; we have $d=ef$ and $q=p^f$, where $q=\Card k$.

We put ${e_1}=e/(p-1)$.  For $n>0$, define $\lambda(n)=\inf(pn, n+e)$, so that
$\lambda(n)=pn$ if $n\le{e_1}$ and $\lambda(n)=n+e$ if $n\ge{e_1}$.

We identify $k^\times$ with a subgroup of $\units$ by the section
$x\mapsto\lim\limits_{n\rightarrow+\infty} y^{q^n}$ (``Teichm{\"u}ller''),
where $y\in\units$ is any preimage of $x\in k^\times$~; it is the subgroup of
solutions of $z^{q-1}=1$.  We have $\units=k^\times U_1$ and $k^\times\cap
U_1=\one$.

Fix an algebraic closure $\bar K$ of $K$.

\bigbreak 
\centerline{\bf 2. Roots of $1$} 
\medskip 

In this section, we study cyclotomic extensions of $K$, give a criterion
for $K^\times$ to contain an element of order~$p$, and determine the level of
an element of order~$p^\alpha$.

\th PROPOSITION {22}
\enonce
Let\/ $\zeta\in\bar K^\times$ be an element of order\/ $n$ prime to\/ $p$.
Then the extension\/ $K(\zeta)\,|\,K$ is unramified of
degree\/~$g$, where $g$ is the order of\/ $q$ in
$(\Z/n\Z)^\times$.
\endth
Recall that, for each $m>0$, $K$ has a unique unramified extension $K_m$ in
$\bar K$ of degree~$m$, that $K_m^\times$ contains an element of order
$(q^m-1)$, and that ``$m'|m$'' is equivalent to $K_{m'}\subset K_m$.  If
$K_m^\times$ has an element of order $l$ prime to $p$, then $l|q^m-1$ (and
conversely).

As $n|q^g-1$, we have $K(\zeta)\subset K_g$.  Therefore $K(\zeta)$ is
unramified over $K$, and hence $K(\zeta)=K_{g'}$ for some $g'|g$.  As
$K_{g'}^\times$ contains an element --- $\zeta$, for example --- of order $n$
prime to $p$, we have $n|q^{g'}-1$, which means that $q^{g'}\equiv1\pmod n$
and $g|g'$, because $g$ is the order of $q$ in $(\Z/n\Z)^\times$.  Thus
$g'=g$.

From now on, let $K_0$ be the maximal unramified subextension of $K$.
\th PROPOSITION {23}
\enonce
Let\/ $\zeta\in\bar K^\times$ be an element of order\/ $p^n$ ($n>0$).
Then the extension\/ $K_0(\zeta)\,|\,K_0$ is totally ramified of
degree\/~$\varphi(p^n)=p^{n-1}(p-1)$, and\/ $1-\zeta$ is a uniformiser
of\/ $K_0(\zeta)$.
\endth

The proof proceeds by induction on $n$.  Put $\xi_n=\zeta$,
$\xi_{n-1}=\xi_n^p$, $\ldots$, $\xi_1=\xi_2^p$.  Also put $K_i=K_0(\xi_i)$ and
$\pi_i=1-\xi_i$.

Let us show that $\pi_1$ is a uniformiser of $K_1$, which is totally ramified
of degree $p-1$ over $K_0$.  As $\xi_1\in K_1^\times$ is an element of order
$p$, we have
$$
{1-\xi_1^p\over1-\xi_1}=1+\xi_1+\xi_1^2+\cdots+\xi_1^{p-1}=0,
$$
which, in terms of $\pi_1=1-\xi_1$, means that
$$
{1-(1-\pi_1)^p\over\pi_1}=p-{p\choose 2}\pi_1+\cdots+(-1)^{p-1}\pi_1^{p-1}=0. 
$$
Thus $\pi_1$ is a root of a degree-$(p-1)$ Eisenstein polynomial over
$K_0$, and hence $K_1\,|\,K_0$ is totally ramified of degree~$p-1$, and
$\pi_1$ is a uniformiser of $K_1$.

Assume now that $K_{n-1}\,|\,K_0$ is totally ramified of
degree~$\varphi(p^{n-1})$ with $\pi_{n-1}=1-\xi_{n-1}$ as a uniformiser.  From
$\xi_n^p-\xi_{n-1}=0$ it follows that
$$
(1-\pi_n)^p-(1-\pi_{n-1})=\pi_{n-1}-
{p\choose 1}\pi_n+
{p\choose 2}\pi_n^2-+\cdots 
+(-1)^p\pi_n^p=0,
$$
which means that $\pi_n$ is a root of a degree-$p$ Eisenstein polynomial
with coefficients in $K_{n-1}$.  Therefore $\pi_n$ is a uniformiser of $K_n$,
which is totally ramified of degree~$p$ over $K_{n-1}$~; in other words,
$K_n\,|\,K_0$ is totally ramified of degree~$\varphi(p^n)$.  This completes
the proof by induction.

(By contrast, for an arbitrary finite extension $K$ of $\Qp$, the extension
$K(\zeta)|K$ may be unramified. Hasse notes that if $p=2$ and $n=2$, then
$K(\sqrt{-1})$ is an {\it unramified\/} quadratic extension of the {\it
  ramified\/} extension $K=\Q_2(\sqrt3)$.  Cf.~ex.~{50}.)

\th PROPOSITION {24}
\enonce The two extensions $\Q_p(\omega)$ ($\omega^{p-1}+p=0$)
and $\Q_p(\zeta)$ ($\zeta^{p-1}+\zeta^{p-2}+\cdots+1=0)$ are isomorphic to
each other.  
\endth

We know that $\pi=1-\zeta$ is a uniformiser of $\Qp(\zeta)$ and that $p$ is
the norm of $\pi$~:
$$
p=(1-\zeta)(1-\zeta^2)\cdots(1-\zeta^{p-1}).
$$
Putting
$$
u_r={1-\zeta^r\over1-\zeta}=1+\zeta+\cdots+\zeta^{r-1}
\equiv r\pmod{\pi}\qquad(0<r<p)
$$ 
we have $p=(1-\zeta)^{p-1}u_1u_2\cdots u_{p-1}$.  But we all know that
$u_1u_2\cdots u_{p-1}\equiv1.2\ldots(p-1)\equiv-1$ in $\F_p^\times$
(``Wilson's theorem'').  So $-p=u\pi^{p-1}$ with $u\in U_1$ depending on
$\zeta$.

But $U_1$ is a $\Z_p$-module, and $p-1$ is invertible in $\Z_p$, so there is a
unique $\eta\in U_1$ such that $\eta^{p-1}=u$.  Therefore we have
$-p=(\eta\pi)^{p-1}$, and thus $\eta\pi$ is a root of $T^{p-1}+p$ in
$\Q_p(\zeta_p)$.

Hence, there is an embedding $\Q_p(\omega)\to\Q_p(\zeta_p)$~; it is an
isomorphism because the two extensions have the same degree over\/ $\Qp$.

(Prop.~{24} answers a question similar to the one answered by prop.~{16}.
Knowing that $\Qp^\times$ has an element of order $p-1$, and that the
extension $\Qp(\zeta)$ is cyclic of degree~$p-1$, which cyclic subgroup of
order~$p-1$ in $\Q_p^\times\!/\Q_p^{\times(p-1)}$ does it correspond to~?
Answer~: the subgroup generated by the image of $-p$.)

\th PROPOSITION {25}
\enonce
For the group\/ $K^\times$ to contain an element of order\/~$p$, it is
necessary and sufficient that\/ $p-1$ divide $e$ and, upon writing\/
$-p=u\pi^e$ ($u\in\ogoth^\times$ and\/ $\pi$ uniformiser of\/ $K$), to
have\/ $\bar u\in k^{\times(p-1)}$.
\endth

The first thing to ask is : Is this independent of the choice of the
uniformiser~?  If $\pi'$ is another uniformiser and we write
$-p=u'\pi'^e$, does one also have $\bar{u'}\in k^{\times(p-1)}$~?
Well, we then have 
$$
u'=u\pi^e\pi'^{-e}=u((\pi/\pi')^{e_1})^{p-1}
$$
and therefore $\bar{u'}\in k^{\times(p-1)}$, because $\bar{u}\in
k^{\times(p-1)}$.

Let us prove the proposition.  Suppose first that $K^\times$ has an element
$\zeta$ of order~$p$ and consider the tower of extensions
$K\,|\,\Qp(\zeta)\,|\,\Qp$.  The absolute ramification index~$e$ is divisible
by the absolute ramification index $p-1$ of $\Qp(\zeta)$, so ${e_1}=e/(p-1)$
is an integer.  The maximal unramified subextension $K_0$ of $K\,|\,\Qp$ is
linearly disjoint from $\Qp(\zeta)$ because the latter is totally ramified;
the extension $K|K_0(\zeta)$ is totally ramified of degree ${e_1}$.

There is a root $\omega$ of $T^{p-1}+p$ in $K_0(\zeta)$ (prop.~{24}).  Further,
$\omega$ is a uniformiser of $K_0(\zeta)$ and, if $\pi$ is a uniformiser of
$K$, then $\omega=\varepsilon\pi^{e_1}$, for some
$\varepsilon\in\ogoth^\times$.  Now observe that
$-p=\omega^{p-1}=\varepsilon^{p-1}\pi^e$, and the unit $\varepsilon^{p-1}$
clearly reduces modulo~$\pi$ to an element of $k^{\times(p-1)}$.

Conversely, suppose that $p-1|e$ and that, writing $-p=u\pi^e$, we have $\bar
u\in k^{\times(p-1)}$.  Write $u=\varepsilon u_1$, where $\varepsilon\in
k^\times$ and $u_1\in U_1$.  By hypothesis, $\varepsilon=\eta^{p-1}$ for some
$\eta\in k^\times$ and, as we have observed earlier, $u_1=\eta_1^{p-1}$ for a
unique $\eta_1\in U_1$.  Then $\eta\eta_1\pi^{e_1}$ is a root of $T^{p-1}+p$
in $K$ so, by prop.~{24}, $K$ contains $\Qp(\zeta)$.

\medskip
{\it 
\rightline{La proposition~25 est grossi{\`e}rement fausse.}
\rightline{\rm --- Anonyme \citer\anonyme(p.~1).}
}

\bigbreak

\th PROPOSITION {26}
\enonce
Suppose that\/ $K^\times$ has an element $\zeta$ of order $p^n$ ($n>0$).
Then\/ $\zeta\in U_a$ but\/ $\zeta\notin U_{a+1}$, where $a=e/\varphi(p^n)$
and $\varphi(p^n)=p^{n-1}(p-1)$.  
\endth

Note that $a$ is an integer by prop.~{23}.  By the same prop., $1-\zeta$ is a
uniformiser of $K_0(\zeta)$, of which  $K$ is a totally ramified extension of
degree~$a$.  For any uniformiser $\pi$ of $K$, writing $1-\zeta=u\pi^a$
($u\in\units)$, we see that $\zeta\in U_a$ but\/ $\zeta\notin U_{a+1}$.

\bigbreak 
\centerline{\bf 3. Raising to the exponent $p$} 
\medskip 

In this section we study the ``homothetie of ratio $p$'' in the $\Z_p$-modules
$U_n$.  We show that this raising to the exponent~$p$ maps $U_n$ into
$U_{\lambda(n)}$ (recall that $\lambda(n)=\inf(pn,n+e)$) for every $n>0$.  The
induced map $\rho_n:U_n/U_{n+1}\to U_{\lambda(n)}/U_{\lambda(n)+1}$ is an
isomorphism in all cases except when\/ $K^\times$ has an element of order\/
$p$ and\/ $n={e_1}$, in which case\/ $\Ker(\rho_n)$ and\/ $\Coker(\rho_n)$
are cyclic of order $p$.  We also specify these two groups.

\th PROPOSITION {27}
\enonce
For every $\eta\in U_n$ ($n>0$), one has $\eta^p\in U_{\lambda(n)}$.
\endth
Let $\pi$ be a uniformiser of $K$ and write 
$\eta=1+a\pi^n$ ($a\in\ogoth$).  We have
$$
(1+a\pi^n)^p=1+pa\pi^n+\cdots+a^p\pi^{pn}.\leqno{(4)}
$$
It is sufficient to restrict to $a\in\units$~; then the valuation of the
second (resp.~last) term on the right is $n+e$ (resp.~$pn$) and the terms not
displayed have valuation $>\lambda(n)$, so
$(1+a\pi^n)^p\equiv1+h(a)\pi^{\lambda(n)}\pmod{\pgoth^{\lambda(n)+1}}$,
with
$$
h(a)=\cases{a^p&if\/ $n<{e_1}$,\cr
a^p-\varepsilon a&if\/ $n={e_1}$,\cr
\phantom{a^p}-\varepsilon a&if\/ $n>{e_1}$.\cr}
\leqno{(5)}
$$
where $\varepsilon\in\units$ is such that $-p=\varepsilon\pi^e$.  Hence
$\eta^p\in U_{\lambda(n)}$, as claimed.  

\th COROLLARY {28}
\enonce
Let\/ $\pi$ be a uniformiser of\/ $K$ and write $-p=\varepsilon\pi^e$
($\varepsilon\in\units$).  The following diagram commutes~:
$$
\diagram{
U_n/U_{n+1}&\droite{(\ )^p}&U_{\lambda(n)}/U_{\lambda(n)+1}\cr
\vfl{}{}{5mm}&&\vfl{}{}{5mm}\cr
k&\droite{h}&k,\cr
} 
$$
where the vertical maps are the isomorphisms\/
$U_i/U_{i+1}\rightarrow k$, $\overline{1+a\pi^i}\mapsto\bar a$ ($a\in\ogoth$)
and the the bottom arrow $h$ is given by $(5)$.
\endth

\th PROPOSITION {29}
\enonce
The map $\rho_n:U_n/U_{n+1}\to U_{\lambda(n)}/U_{\lambda(n)+1}$ is an
isomorphism for all $n>0$ except when\/ $K^\times$ has an element of order\/
$p$ and\/ $n={e_1}$, in which case\/ $\Ker\rho_n$ and\/ $\Coker\rho_n$
are cyclic of order $p$.
\endth

The $\F_p$-linear map $(5)$ from $k$ to $k$ is an isomorphism in all cases
except when $n={e_1}$ and $\bar\varepsilon\in k^{(p-1)}$, which happens
precisely when $K^\times$ has an element of order $p$ (prop.~{25}). In this
case, the kernel has at least two (for if $\bar\varepsilon=b^{p-1}$ for some
$b\in k^\times$, then $0,b\in\Ker\rho_n$) and at most $p$ elements (because
the polynomial $T^p-\bar\varepsilon T$ can have at most $p$ roots), hence it
has exactly $p$ elements.  Consequently, $\Coker\rho_n$ is also cyclic of
order~$p$.
 
\th PROPOSITION {30}
\enonce
For every $n>{e_1}$, the map $(\ )^p:U_n\to U_{n+e}$ is bijective.   
\endth
Let $\pi$ be a uniformiser of $K$ and write $-p=\varepsilon\pi^e$
($\varepsilon\in\units$).  Let $y\in U_{n+e}$ and write $y=1+b\pi^{n+e}$
(with $b\in\ogoth$).  We seek a root of $x^p=y$ such that $x=1+a\pi^n$ for
some $a\in\ogoth$.  This leads to the equation
$$
1+b\pi^{n+e}=1+pa\pi^n+\cdots+a^p\pi^{np}.
$$
As we have seen, all the terms on the right, except the first two, have
valuation $>n+e$, in view of $n>{e_1}$.  The equation can therefore be
rewritten
$$
b=-\varepsilon a+\pi f(a)
$$
for some polynomial $f\in\ogoth[T]$.  Reducing modulo $\pi$ yields
$\bar b=-\bar\varepsilon\bar a$, and since $\bar\varepsilon\neq0$, this
equation has a unique root.  By Hensel's lemma, the same holds for
$x^p=y$. (The injectivity also follows from prop.~{26}).
Cf.~\citer\serre() (pp.~212--3). 

\th PROPOSITION {31}
\enonce
For every $n>{e_1}$, the $\Z_p$-module $U_n$ is free of rank $d=[K:\Qp]$.   
\endth
We have seen that $U_n$ has no element of order~$p$ (prop.~{26}), so it is
sufficient to show that the $\F_p$-space $U_n/U_n^p$ is of dimension $d$.
This is the case because $U_n^p=U_{n+e}$ (prop.~{30}) is of index
$q^e=p^{fe}=p^d$ in $U_n$.  

\th COROLLARY {32}
\enonce
The $\Z_p$-module $U_1$ is finitely generated of rank $d$.   
\endth

Suppose that $K^\times$ has an element of order~$p$, and let ${}_p\mu$ be the
$p$-torsion of $K^\times$.  We have seen (prop.~{26}) that
$$
{}_p\mu\subset U_{e_1},\qquad 
{}_p\mu\cap U_{{e_1}+1}=\one.
$$
We also know (prop.~{30}) that $U_{p{e_1}+1}\subset (U_{p{e_1}}\cap U_1^p)$.
So  we get a sequence 
$$
1\rightarrow
{}_p\mu\rightarrow
U_{{e_1}}/U_{{{e_1}}+1}\hfl{(\ )^p}{}{8mm}
U_{p{{e_1}}}/U_{p{{e_1}}+1}\rightarrow
\bar U_{p{e_1}}\rightarrow1.
\leqno{(6)}$$ 
in which $\bar U_{p{{e_1}}}=U_{p{e_1}}/(U_{p{e_1}}\cap U_1^p)$.  

\th PROPOSITION {33}
\enonce
Suppose that $K^\times$ has an element of order $p$.  The sequence $(6)$
is then exact. 
\endth
It is clear that ${}_p\mu\subset\Ker(\ )^p$~; as the latter group has $p$ 
elements (prop.~{29}), the inclusion is an equality.  It is also clear that
$U_{e_1}^p\subset (U_{p{e_1}}\cap U_1^p)$~; let us show equality here too. If
$\eta\in U_r$ but $\eta\notin U_{r+1}$ for some $r<{e_1}$, then $\eta^p\in
U_{pr}$ and $\eta^p\notin U_{pr+1}$ (cor.~{28}).  As $pr+1<p{e_1}$,
we have $\eta^p\notin U_{p{e_1}}$, which was to be shown.  Remark that the
proof also gives the exactness of the sequence   
$$
1\rightarrow
{}_p\mu\rightarrow
U_{{e_1}}\hfl{(\ )^p}{}{8mm}
U_{p{{e_1}}}\rightarrow
\bar U_{p{e_1}}\rightarrow1.
$$ 
Therefore the kernel and cokernel of $(\ )^p:U_{e_1}\rightarrow U_{pe_1}$ are
trivial when $K^\times$ doesn't have an element of order~$p$~; when it does,
$\Ker(\ )^p$ and $\Coker(\ )^p$ are both of order~$p$.

{\it Upon choosing a uniformiser\/ $\pi$.}  Let $\pi$ be an $\ogoth$-basis of
$\pgoth$ and write $-p=\varepsilon\pi^e$ ($\varepsilon\in\units$).  Consider
the following diagram, in which the two middle vertical arrows are induced by
the $\ogoth$-bases $\pi^{e_1}$, $\pi^{pe_1}$ of $\pgoth^{e_1}$,
$\pgoth^{pe_1}$, in which $\wp_\varepsilon(a)=a^p-\bar\varepsilon a$, and
where we have put $t=e_1$ to save space.  We have seen (cor.~{28}, prop.~{33})
that with these choices, it is commutative.
$$
\def\\{\mskip-2\thickmuskip}
\def\droite#1{\\\hfl{#1}{}{8mm}\\}
\diagram{
1&\rightarrow&{}_p\mu&\droite{}&
U_{t}/U_{{t}+1}&\droite{(\ )^p}
&U_{p{t}}/U_{p{t}+1}&\droite{}&\bar U_{pt}&\rightarrow&1\cr
&&\vfl{}{}{5mm}&&\vfl{}{}{5mm}&&\vfl{}{}{5mm}&&\vfl{}{}{5mm}\cr
0&\rightarrow&\Ker\wp_\varepsilon&\droite{}&k&\droite{\wp_\varepsilon}
&k&\droite{}&k/\wp_\varepsilon(k)&\rightarrow&0.\cr
}
$$

{\it Upon choosing a primitive $p$-th root of\/ $1$.}  If we choose such a
root $\zeta$, then $\pi_1=1-\zeta$ is an $\ogoth$-basis of $\pgoth^{e_1}$
(prop.~{23}), and $-p\pi_1=p(\zeta-1)$ a basis of $\pgoth^{p{e_1}}$.  These
bases lead to the commutative diagram
$$
\def\\{\mskip-2\thickmuskip}
\def\droite#1{\\\hfl{#1}{}{8mm}\\}
\diagram{
1&\rightarrow&{}_p\mu&\droite{}&
U_{{e_1}}/U_{{{e_1}}+1}&\droite{(\ )^p}
&U_{p{{e_1}}}/U_{p{{e_1}}+1}&\droite{}&\bar U_{p{e_1}}&\rightarrow&1\cr
&&\vfl{}{}{5mm}&&\vfl{}{}{5mm}&&\vfl{}{}{5mm}&&\vfl{}{}{5mm}\cr
0&\rightarrow&\F_p&\droite{}&k&\droite{\wp}
&k&\droite{S_{k|\F_p}}&\F_p&\rightarrow&0,\cr
}
$$ 
in which the vertical maps are isomorphisms, with $\wp(x)=x^p-x$ and
$S_{k|\F_p}$ the trace map.  Here we have $\wp$ instead of $\wp_\varepsilon$
because $-p\pi_1/\pi_1^p$ is a $1$-unit (``Wilson's theorem'', cf.~proof of
prop.~{24}). 

Explicitly, the isomorphism $\bar U_{pe_1}\rightarrow\Fp$ --- induced by the
choice $\zeta\in{}_p\mu$ of a generator --- sends $\bar\eta$ to
$S_{k|\Fp}(\hat c)$, where $\eta\in U_{pe_1}$, $c=(1-\eta)/p(1-\zeta)$ and
$\hat c$ denotes its image in $k/\wp(k)$.  \footnote{$({}^3)$}{Hasse's
convention in \citer\hassepotenz(p.~233) amounts to taking $p\pi_1$ as an
$\ogoth$-basis for $\pgoth^{pe_1}$~; one would then have to replace $\wp$ by
$-\wp$ in the displayed diagram.  Compare the exponent $S_{\goth
l}({\alpha-1\over l\lambda_0})$ in formula $(13)$ in footnote $({}^2)$ with
our $S_{k|\Fp}\left({1-\eta\over p(1-\zeta)}\right)$.}  In particular, when
$k=\Fp$, we have $\wp=0$ (Fermat's ``little'' theorem), and $k/\wp(k)=k=\Fp$.
The isomorphism is then simply
$\bar\eta\mapsto(1-\eta)/p(1-\zeta)\pmod\pgoth$.

{\it Upon choosing a $(p-1)$-th root of\/ $-p$.}  Equivalently, as we
discovered recently in \citer\henselmultii(p.~211), one can choose a
$(p-1)$-th root $\Pi$ of $-p$ in $K$ (cf.~prop.~{24})~; then $\Pi$ is an
$\ogoth$-basis of $\pgoth^{e_1}$ and $\Pi^p$ a basis of $\pgoth^{p{e_1}}$~;
with these choices for the two middle vertical arrows, the above diagram is
commutative.

Notice that if we fix $\Pi$ and $\zeta$, there is a unique ``natural''
bijection, sending $\Pi$ to $\zeta$, between the set $R$ of $(p-1)$-th roots
of $-p$ and the set $P$ of $p$-th roots of~$1$.  Indeed, $R$ is a
$({}_{p-1}\mu)$-torsor, and $P$ a $(\Z/p\Z)^\times$-torsor.  But we have a
natural isomorphism $\xi\mapsto\bar\xi:{}_{p-1}\mu\rightarrow(\Z/p\Z)^\times$
of groups (over $\Zp$, so to speak), induced by the passage to the quotient
$\Zp\rightarrow\Fp$.  The bijection $R\rightarrow P$ in question is
$\xi\Pi\mapsto\zeta^{\bar\xi}$ (for every $\xi\in{}_{p-1}\mu$).

All this would be true for any two torsors under the ``same'' group
${}_{p-1}\mu=(\Z/p\Z)^\times$.  But more is true here~: there is a unique
bijection $\zeta\mapsto\Pi_\zeta:P\rightarrow R$ such that
$\Pi_\zeta/(1-\zeta)$ is a $1$-unit for every $\zeta\in P$.  Indeed, as we saw
during the proof of prop.~{24}, $u=-p/(1-\zeta)^{p-1}$ is a $1$-unit for every
$\zeta\in P$~; denoting by $\eta\in U_1$ the unique $(p-1)$-th root of $u$,
take $\Pi_\zeta=\eta(1-\zeta)$.  Moreover, this bijection is ``equivariant''~:
$\Pi_{\zeta^r}=\chi(r)\Pi_\zeta$ for every $r\in(\Z/p\Z)^\times$, where
$\chi(r)\in{}_{p-1}\mu$ is the ``Teichm{\"u}ller'' lift of $r$~:
$\overline{\chi(r)}=r$. 


{\it The case\/ $p=2$.}  The choice $\zeta=-1$ (or $\Pi=-2$) is forced upon
us.  Consequently, the isomorphisms are canonical.  Thus, when $k=\F_2$, the
isomorphism $\bar U_{2e}\rightarrow\F_2$ is
$\bar\eta\mapsto(1-\eta)/4\pmod\pgoth$, which is the same as the more familiar
$\bar\eta\mapsto(\eta-1)/4\pmod\pgoth$ because $\pgoth$ is even.

\bigbreak 
\centerline{\bf 4. The multiplicative group\/ $(\ogoth/\pgoth^n)^\times$}  
\medskip

\th PROPOSITION {34}
\enonce
The group\/ $(\ogoth/\pgoth^n)^\times$ is the direct product of its
subgroups\/ $k^\times$ and $U_1/U_n$.  For\/ $n>{e_1}$, the restriction of\/
$U_1\to U_1/U_n$ to the torsion subgroup\/ $W\subset U_1$ is injective, and
the image of\/ $W$ is a direct factor of\/ $U_1/U_n$.  For\/ $n>{e_1}+e$,
the group\/ $U_1/U_n$ is the direct product of the image of\/ $W$ with\/ $d$
cyclic\/ $p$-groups (of order\/ $>1$).
\endth
The exact sequence
$
1\rightarrow U_1/U_n\rightarrow\units\!/U_n\rightarrow
                                     k^\times\rightarrow1
$
has a canonical splitting for every $n>0$, because $U_1/U_n$ is a
$p$-group and the order of $k^\times$ is prime to~$p$.  

For $n>{e_1}$, we have $W\cap U_n=\one$ (prop.~{26}), so the restriction of the
projection $U_1\to U_1/U_n$ is injective on $W$ and the restriction of $U_1\to
U_1/W$ injective on $U_n$.  Choosing a section $s:U_1/W\to U_1$, which is
possible because the $\Zp$-module $U_1/W$ is free (of rank $d$, prop.~{31}),
write $U_1=W\times (U_1/W)$; then, quotienting modulo the sub-$\Zp$-module
$U_n\subset U_1/W$, we still get a direct product decomposition
$U_1/U_n=W\times (U_1/WU_n)$.

Let us show that for $n>{e_1}+e$, the (finite) $\Zp$-module $M=U_1/WU_n$ is a
direct product of $d$ cyclic $p$-groups (of order $>1$).  As it can be
generated by $d$ elements (take the image of any $\Zp$-basis of $U_1/W$), it
is sufficient to show that it cannot be generated by $<d$ elements.  It is in
fact sufficient to exhibit one subgroup which cannot be generated by $<d$
elements.  Now, the $\Zp$-module $U_{n-e}$ is free of rank~$d$ (prop.~{31}),
and $U_{n-e}^p=U_n$.  Thus the subgroup $U_{n-e}/U_n\subset U_1/WU_n$, a
$d$-dimensional vector $\Fp$-space, cannot be generated by $<d$ elements.
Note that $U_{n-e}/U_n$ is a subgroup of $U_1/WU_n$ because $W\cap
U_{n-e}=\one$ (prop.~{26}).

Information about the groups $U_1/U_n$ for small $n$ can be found in
\citer\nakagoshi().  

The methods of ideal theory had not succeeded in determining the structure of
the groups $(A/{\goth a})^\times$, where $A$ is the ring of integers in a
number field and ${\goth a}\subset A$ an ideal.  As a curiosity, let us
mention that Wilson's theorem, which we have had to invoke a certain number of
times, and which was generalised by Gauss in his {\it Disquisitiones\/} (\S78)
from $(\Z/p\Z)^\times$ to $(\Z/a\Z)^\times$ for any $a>0$, can be further
generalised by local means to all $(A/{\goth a})^\times$.  There are four
possibilities for the product of all elements in this group, and precise
conditions for each of these possibilities to occur have been given
\citer\wilson().

\bigbreak
\centerline{\bf IV.   Unramified kummerian extensions}
\bigskip

This Part and the next contain the proofs of our main propositions about
discriminants of unramified extensions, unramified kummerian extensions, their
rings of integers, the $p$-tic character (prop.~{38}), and the filtration on
$K^\times\!/K^{\times l}$.  But before proving prop.~{16} and deriving a few
corollaries, among them prop.~{15}, let us compute the greatest $n$ such that
the $\F_l$-dimension of $\bar U_{n}$ is $\neq0$.

Recall the notation in vigour~: $K$ is a finite extension of $\Q_p$, $\ogoth$
is its ring of integers, with residue field $k$ having $q$ elements.  The
units $\ogoth^\times$ will also be denoted $U_0$~; we denote by $U_{n+1}$ the
$1$-units of level~$>n$ (cf.~Part~II).

Let $l$ be a prime number such that $K^\times$ has an element of order~$l$.
We define ${e_1}=e/(p-1)$, where $e$ is the absolute ramification index of $K$
if $l=p$, and $e=0$ if $l\neq p$.  The filtration $(U_n)_n$ on $K^\times$
induces the filtration $(\bar U_n)_n$ on $K^\times\!/K^{\times l}$.

Thus, in both cases ($l\neq p$ and $l=p$), we have $\dim_{\F_l}\bar
U_{p{e_1}}=1$ and $\dim_{\F_l}\bar U_{n}=0$ for $n>p{e_1}$~: the case $l\neq
p$ follows from the fact that $U_1$ is a $\Z_p$-module and $k^\times$ a cyclic
group of order divisible by~$l$~; for the case $l=p$, see prop.~{29} and~{33}.

\medskip\centerline{***}\medskip

Prop.~{16} says that the $\F_l$-line in $K^\times\!/K^{\times l}$ which
corresponds to the unramified degree-$l$ extension of $K$ is the one which
lies the deepest in the induced filtration, namely the line $\bar
U_{p{e_1}}$. 

{\bf Proof of prop.~{16}.}  Suppose first that $l\neq p$~; then $l|q-1$, as
$k^\times$ has an element of order~$l$.  Let $M$ be the unramified degree-$l$
extension of $K$~; its residue field $m$ is the degree-$l$ extension of $k$.
To show that $M$ is associated to the $\F_l$-line $\bar U_{p{e_1}}=\bar
U_0=\units\!/\ogoth^{\times l}=k^\times\!/k^{\times l}$, it is sufficient to
show that $k^\times\subset m^{\times l}$.

As the group $m^\times$ is cyclic of order $q^l-1$, and as the subgroup
$k^\times$ is of order $q-1$, we have $k^\times=m^{\times a}$, where
$
a={(q^l-1)/(q-1)}
$
But, as $q\equiv1\pmod l$, we have $1+q+q^2+\cdots+q^{l-1}\equiv0\pmod l$,
so $a$ is a multiple of $l$, and hence $k^\times\subset m^{\times l}$.  It
would have also sufficed to remark that the map $k^\times\!/k^{\times
  l}\rightarrow m^\times\!/m^{\times l}$ is trivial, by Kummer theory.

Let us remark that it is not so much the primality of $l$, but the fact that
$l$ divides $q-1$ (and hence is prime to $p$) which has been used in the proof.
Thus, {\it for every divisor\/ $s\,|\,q-1$, the degree-$s$ unramified
extension of\/ $K$ is the kummerian extension obtained by adjoining\/ $\root
s\of{\ogoth^\times}$ to\/ $K$.} 

\smallskip

Let us come to the case $l=p$.  Let $M$ now be the degree-$p$ unramified
extension of $K$~; its absolute ramification index is the same as that of $K$,
namely $e$.  We have to show that every element of $U_{p{e_1}}$ has a $p$-th
root in $M$.  Denoting by $(V_n)_{n>0}$ the filtration on the $1$-units of
$M$, this amounts to showing that the map $\bar U_{p{e_1}}\rightarrow \bar
V_{p{e_1}}$ (induced by the inclusions $U_i\subset V_i$) is trivial.

The map $\bar U_{p{e_1}}\rightarrow \bar V_{p{e_1}}$, whose triviality is in
question, and which is therefore denoted by $1_?$ in the multiplicative
notation in use, is part of the following commutative diagram~:
$$
\def\\{\mskip-2\thickmuskip}
\def\droite#1{\\\hfl{#1}{}{8mm}\\}
\diagram{
1&\rightarrow&{}_p\mu&\rightarrow&
V_{e_1}/V_{{e_1}+1}&\rightarrow
&V_{p{e_1}}/V_{p{e_1}+1}&\rightarrow&
\bar V_{p{e_1}}&\rightarrow&1\phantom{.}\cr
&&\ufl{=}{}{5mm}&&\ufl{}{}{5mm}&&\ufl{}{}{5mm}&&\ufl{}{1_?}{5mm}\cr
1&\rightarrow&{}_p\mu&\rightarrow&
U_{e_1}/U_{{e_1}+1}&\rightarrow
&U_{p{e_1}}/U_{p{e_1}+1}&\rightarrow&\bar U_{p{e_1}}&\rightarrow&1.\cr}
$$ 
Upon choosing a generator $\zeta$ of ${}_p\mu$, or,
equivalently, a $(p-1)$-th root of $-p$, --- cf.~the discussion after
prop.~{33} --- the above diagram gets identified with 
$$
\def\\{\mskip-2\thickmuskip}
\def\droite#1{\\\hfl{#1}{}{8mm}\\}
\diagram{
0&\rightarrow&\F_p&\droite{}&m&\droite{\wp}
&m&\droite{S_{m|\F_p}}&\F_p&\rightarrow&0,\cr
&&\ufl{=}{}{5mm}&&\ufl{}{}{5mm}&&\ufl{}{}{5mm}&&\ufl{}{0_?}{5mm}\cr
0&\rightarrow&\F_p&\droite{}&k&\droite{\wp}
&k&\droite{S_{k|\F_p}}&\F_p&\rightarrow&0\cr}
$$
in which $m$ is the residue field of $M$ and the first three vertical arrows
are inclusions.  We have to show that the last arrow $0_?:\F_p\rightarrow\F_p$
is indeed~$0$.  But this is the case because $m$ is the degree-$p$ extension
of $k$~: we have $k\subset\wp(m)$.

The following corollaries are immediate.

\th COROLLARY {35}
\enonce 
Let\/ $\tilde K$ be the maximal unramified extension of\/ $K$.  The kernel of
the map\/ $K^\times\!/K^{\times l}\rightarrow \tilde K^\times\!/\tilde
K^{\times l}$ is\/ $\bar U_{p{e_1}}$.
\endth

\th COROLLARY {36}
\enonce 
Let\/ $\tilde K$ be the maximal tamely ramified extension of\/ $K$.  The map\/
$K^\times\!/K^{\times l}\rightarrow \tilde K^\times\!/\tilde K^{\times l}$ is
trivial if\/ $l\neq p$~; it has the kernel\/ $\bar U_{p{e_1}}$ if\/ $l=p$.
\endth

When $l=2$ and $K$ is the extension of $\Q_2$ obtained by adjoining $1$, or
$\root3\of2$, or $\sqrt3$, or $\sqrt{-1}$, A.~Kraus \citer\kraus(p.~376) does
explicit calculations in each case to determine the units which become squares
in the maximal unramified extension~; cf.~cor.~{46}.  Cor.~{35} and~{36} give
a criterion for an element of $K$ to become an $l$-th power in $\tilde K$ for
any finite extension $K$ of $\Qp$ having a primitive $l$-th root of $1$, where
the primes $p$ and $l$ are otherwise arbitrary.

In view of the discussion after prop.~{33}, things can be made more explicit.
Let $u\in\ogoth$ be such that its image generates\/ $k/\wp(k)$.  Then the
image of $\eta=1-up(1-\zeta)$ generates $\bar U_{pe_1}$.  Let $\root
p\of\eta$ be a root of $T^p-\eta$~; then $\root p\of\eta-1$ is a root of
$$
T^p+pT^{p-1}+\cdots+pT+up(1-\zeta),
$$ 
where the coefficients of the suppressed terms are all divisible by $p$.
Dividing  throughout by $(1-\zeta)^p$ and setting $S=T/(1-\zeta)$, we see that 
$\rho=(\root p\of\eta-1)/(1-\zeta)$ is a root of
$$
h(S)=S^p+\cdots+{p\over (1-\zeta)^{p-1}}S+{up\over(1-\zeta)^{p-1}},
$$ 
where the coefficients of the suppressed terms are divisible by $\pi$ ---
any uniformiser of $K$.  Denoting by $\hat h$, $\hat u$ the reductions modulo
$\pi$ and recalling that $-p/(1-\zeta)^{p-1}$ is a $1$-unit (``Wilson's
theorem''), we see that 
$$
\hat h(S)=S^p-S-\hat u.
$$ 
As we had chosen $u\in\ogoth$ such that the image of $\hat u\in k$ generates
$k/\wp(k)$, this shows that $\ogoth[\rho]/\pi$ is the degree-$p$ extension of
$k$. (Cf.~\citer\gras(p.~60)).


\th COROLLARY {37}
\enonce
For\/ $\eta=1-up(1-\zeta)$, with\/ $u\in\ogoth$ such that\/
$S_{k|\Fp}(\hat u)\neq0$, the ring of integers of\/ $K(\root p\of\eta)$ is\/
$\ogoth[(\root p\of\eta-1)/(1-\zeta)]$.
\endth

(If we had worked more generally with a characteristic-$0$ field $K$ complete
with respect to a discrete valuation whose residue field $k$ is perfect of
prime characteristic $p$ and which contains a primitive $p$-th root of $1$,
the choice of such a root $\zeta$ would still lead to an isomorphism $\bar
U_{pe_1}\rightarrow k/\wp(k)$, but these $\Fp$-spaces need no longer be
$1$-dimensional.)

Abbreviate $D=\bar U_{pe_1}$, $L=K(\root p\of D)$ and $G=\Gal(L|K)$.  We have
seen that the choice of a generator $\zeta\in{}_p\mu$ allows us to identify
$D$ and ${}_p\mu$ with $\Z/p\Z$.  On the other hand, the group $G$ has a
canonical generator $\varphi$~: the unique element such that
$\varphi(\alpha)\equiv \alpha^q\pmod{\pgoth_L}$ for every $\alpha\in\ogoth_L$,
where $q=\Card k$~; it can be used to identify $G$ with $\Z/p\Z$.  Can the
pairing
$$
G\times D\rightarrow{}_p\mu,\qquad
\langle\sigma,\bar\eta\rangle={\sigma(\xi)\over\xi} \quad
(\xi^p=\eta),\leqno{(7)}
$$
which comes from Kummer theory, be made explicit in terms of these
identifications~?  

\th PROPOSITION {38} (``Poor man's explicit reciprocity law'')
\enonce
Choose and fix a generator $\zeta\in{}_p\mu$.  For\/ $a\in\Z/p\Z$ and\/
$\eta\in U_{pe_1}$, we have
$\langle\varphi^a,\bar\eta\rangle=\zeta^{a.S_{k|\Fp}(\hat c)}$, where
$\displaystyle c={1-\eta\over p(1-\zeta)}$ and $\hat c$ is its image in
$k/\wp(k)$.  
\endth
(More prosaically, if we identify $G$ (resp.~$\bar U_{pe_1}$ and ${}_p\mu$)
with $\Z/p\Z$ using $\varphi$ (resp.~$\zeta$), then the pairing (7) is just
$\langle a,b\rangle=ab$.  Recall (cf.~the discussion after prop.~{33}) that
the isomorphism $\bar U_{pe_1}\rightarrow\Fp$ is given by $\bar\eta\mapsto
S_{k|\Fp}(\hat c)$, with $c=(1-\eta)/p(1-\zeta)$.)  Note, in particular, that
$\langle\varphi,\bar\eta\rangle=\zeta$ if $\bar\eta\in\bar U_{pe_1}$
corresponds to $1\in\Fp$ under the isomorphism $\bar U_{pe_1}\rightarrow\Fp$
induced by a generator $\zeta\in{}_p\mu$.)

It is sufficient to show that the choice of $\zeta$ leads to an identification
of the Kummer pairing $G\times\bar U_{pe_1}\rightarrow{}_p\mu$ with the
Artin-Schreier pairing $\Gal(l|k)\times k/\wp(k)\rightarrow\Z/p\Z$, where $l$
denotes --- not a prime number as hitherto, but momentarily the residue
field of $L$~; the prop.\ will follow from this because the Artin-Schreier
pairing is just $\langle\hat\varphi^a,b\rangle=aS_{k|\Fp}(b)$, where
$\hat\varphi$ --- the image of $\varphi$ --- is the canonical generator of
$\Gal(l|k)$.

(Recall that there is a reciprocity map $K^\times\rightarrow G$, and hence a
pairing $(\ ,\ )_K:K^\times\times\bar U_{pe_1}\rightarrow{}_p\mu$ deduced from
the Kummer pairing.  Remarking that the reciprocity map sends uniformisers of
$K$ to $\varphi$, prop.~{38} allows us to retrieve the last prop.\ of
\citer\serre(p.~230).)

\medskip\centerline{***}\medskip

Let us come to the proof of prop.~{15}.  For the time being, let $l$ be any
prime for which the finite extension $K$ of $\Qp$ has a primitive $l$-th root
of~$1$~; for the application to discriminants, the case $l=2$ is sufficient.

Let $M$ be an unramified extension of $K$~; it has the same $e$ as $K$~:
$e=v_K(l)=v_M(l)$ ($=0$ if $p\neq l$).  We denote by $(\bar V_n)_{n>0}$ the
filtration induced on $M^\times\!/M^{\times l}$ by the canonical filtration of
$M^\times$.  The residue field $m$ of $M$ is a finite extension of~$k$.

\th PROPOSITION {39}
\enonce
The norm map $N_{M|K}:M^\times\rightarrow K^\times$ induces an isomorphism 
$\bar V_{p{e_1}}\rightarrow\bar U_{p{e_1}}$.
\endth

If $p\neq l$, then $\bar V_{p{e_1}}=m^\times\!/m^{\times l}$ and $\bar
U_{p{e_1}}=k^\times\!/k^{\times l}$, and the map induced by $N_{M|K}$ on these
two spaces is the same as the isomorphism induced by $N_{m|k}$ (which is
surjective $m^\times\rightarrow k^\times$~; cf.~also the proof of prop.~{16}).

Let us come to the  case $p=l$.  We have a commutative diagram of horizontal
isomorphisms 
$$
\def\\{\mskip-2\thickmuskip}
\def\droite#1{\\\hfl{#1}{}{8mm}\\}
\diagram{
\bar V_{p{e_1}}&\droite{}&m/\wp(m)\cr
\vfl{}{}{5mm}&&\vfl{}{S_{m|k}}{5mm}\cr
\bar U_{p{e_1}}&\droite{}&k/\wp(k)\cr
}
$$ 
in which $S_{m|k}$ is the trace map~; it suffices to show that it is an
isomorphism.  This follows from the fact that $S_{m|\Fp}=S_{k|\Fp}\circ
S_{m|k}$, where $S_{m|\Fp}:m/\wp(m)\rightarrow\Fp$ and
$S_{k|\Fp}:k/\wp(k)\rightarrow\Fp$ are the trace maps, which are isomorphisms.
Therefore the norm map $\bar V_{p{e_1}}\rightarrow\bar U_{p{e_1}}$ is an
isomorphism, as claimed. 

Remark that the isomorphism $S_{m|k}$ sits in the commutative diagram of
$\Fp$-linear maps
$$
\def\\{\mskip-2\thickmuskip}
\def\droite#1{\\\hfl{#1}{}{8mm}\\}
\diagram{
0&\rightarrow&\F_p&\droite{}&m&\droite{\wp}
&m&\droite{}&m/\wp(m)&\rightarrow&0\cr
&&\vfl{0}{}{5mm}&&\vfl{}{}{5mm}&&\vfl{}{}{5mm}&&\vfl{}{S_{m|k}}{5mm}\cr
0&\rightarrow&\F_p&\droite{}&k&\droite{\wp}
&k&\droite{}&k/\wp(k)&\rightarrow&0\cr
}
$$
which should be contrasted with the diagram appearing before cor.~{35}.  The
above diagram results from
$$
\def\\{\mskip-2\thickmuskip}
\def\droite#1{\\\hfl{#1}{}{8mm}\\}
\diagram{
1&\rightarrow&{}_p\mu&\rightarrow&
V_{e_1}/V_{{e_1}+1}&\rightarrow
&V_{p{e_1}}/V_{p{e_1}+1}&\rightarrow&
\bar V_{p{e_1}}&\rightarrow&1\cr
&&\vfl{1}{}{5mm}&&\vfl{}{}{5mm}&&\vfl{}{}{5mm}&&\vfl{}{}{5mm}\cr
1&\rightarrow&{}_p\mu&\rightarrow&
U_{e_1}/U_{{e_1}+1}&\rightarrow
&U_{p{e_1}}/U_{p{e_1}+1}&\rightarrow&\bar U_{p{e_1}}&\rightarrow&1\cr
}
$$ 
in which the arrows are the norm maps, upon choosing a $p$-th root of unity.

Till the end of this Part~III, take $l=2$, which is allowed because $K^\times$
has an element of order~$2$, namely~$-1$.

\th PROPOSITION {40}
\enonce
Let $M|K$ be the unramified quadratic extension.  Its discriminant $d_{M|K}=u$
is the unique element\/ $\neq1$ in $\bar U_{p{e_1}}$.   
\endth
This is clear if $p\neq2$, because the quadratic extension $m$ of $k$ is
obtained by adjoining $\sqrt u$.  The discriminant of $T^2-u$ is $4u$,
which is the same as $u$ modulo squares (of units).

Suppose that $p=2$.  The quadratic extension $m$ of $k$ is obtained by
adjoining a root of the polynomial $T^2-T-\alpha$, for some
$\alpha\notin\wp(k)$.  Let $\omega\in\ogoth$ be a lift of $\alpha$~; then $M$
is obtained by adjoining a root $t$ of $T^2-T-\omega$.  As the ring of
integers of $M$ is $\ogoth[t]$, the discriminant $d_{M|K}$ equals
$d=1+4\omega$, modulo $\ogoth^{\times2}$.  Because $4\ogoth=\pgoth^{2e}$, one
has $d\in U_{2e}$~; one also has $d\notin\ogoth^{\times 2}$, because $K(\sqrt
d)=M$.  Thus $\bar d=u$ in $\bar U_{2e}$, which was to be proved.

Prop.~{15} says that the discriminant $d_{L|K}$ of an unramified extension
$L|K$ lies in the $\F_2$-line $\bar U_{p{e_1}}=\{1,u\}$, and equals $1$ if the
degree $[L:K]$ is odd, $u$ if the degree in question is even.  Prop.~{40} was
the case $[L:K]=2$.  (Recall that $e_1=0$ if $p\neq2$.)

{\bf Proof of prop.~{15}}.  As $L|K$ is galoisian, it contains
$K(\sqrt{d_{L|K}})$, which is therefore an unramified extension of $K$, and
hence (prop.~{16}) $d_{L|K}\in\bar U_{p{e_1}}$.  If $[L:K]$ is odd,
necessarily $K(\sqrt{d_{L|K}})=K$, and therefore $d_{L|K}=1$.  If $[L:K]=2g$
is even, $L$ contains the quadratic unramified extension (prop.~{16})
$M=K(\sqrt u)$ of $K$, and
$$\let\\\strut
d\\_{L|K}=d\\_{M|K}^g.N\\_{M|K}(d\\_{L|M})\quad
(\hbox{\it Schachtelungssatz\/}) \leqno{(8)}
$$ 
according to \citer\frohlich(p.~19) or \citer\cassels(p.~143) (or
\citer\hasse(Kap.~25), or \citer\waterdisc(p.~213) or even Hilbert ---
$D=d^rn({\goth d})$ --- \citer\hilbert(p.~22)), where $N_{M|K}$ is the norm
map $M^\times\rightarrow K^\times$~; it induces the isomorphism $\bar
V_{p{e_1}}\rightarrow\bar U_{p{e_1}}$ (prop.~{39}).  The proof is now
complete~: if $g$ is odd, the first factor on the right in $(8)$ is $1$ and
the second $u$ (prop.~{40}), whereas if $g$ is even, the first factor is $u$
(by the induction hypothesis and prop.~{39}), and the second is $1$
(prop.~{40}).  Irrespective of the parity of~$g$, the product is always $u$,
which proves the proposition.  Finally, it is clear that when\/ $p=2$, none of
the other\/ $2^{d+1}-2$ ($d=[K:\Qp]$) elements of\/
$\ogoth^\times\!/\ogoth^{\times2}$ can be a discriminant, for $v(d_{L|K})>0$
if $L|K$ is a ramified extension, so
$d_{L|K}\notin\ogoth^\times\!/\ogoth^{\times2}$.

Notice that in the case $K=\Q_p$, prop.~{15} reduces to the
Stikcelberger-Hensel th.~{5} if $p\neq2$, to their th.~{6} if $p=2$.  The
analogy between prop.~{15}, which determines relative discriminants at the
$p$-adic places, and prop.~{9}, which determines them at the archimedean
places, is striking.

The Pellet-Vorono{\"\i} theorem (th.~{3} and~{4}) can now be extended from
$\Fp$ ($p\neq2$) to any finite field $k$.

\th COROLLARY {41} 
\enonce
Let\/ $p$ be an odd prime (resp.~$p=2$), $k$ a finite extension of\/ $\Fp$
and\/ $f\in k[T]$ a separable polynomial.  The discriminant of\/ $k[T]/fk[T]$
is trivial in the\/ $2$-element multiplicative group\/
$k^\times\!/k^{\times2}$ (resp.~additive group\/ $k/\wp(k)$) if and only if
the number of even-degree irreducible factors of\/ $f$ is even.
\endth

\medskip\centerline{***}\medskip

{\it Relative discriminants of ramified extensions.}  Having found the correct
relative analogue (prop.~{15}) of the Stickelberger-Hensel congruence (th.~{5}
and~{6}), what about the relative version of Dedekind's th.~{10}, which says
that the discriminant of a ramified extension of $\Q_2$ is $\equiv0\pmod4$~?
Part~VII is devoted to this problem, where the valuation of the relative
discriminant of a prime-degree cyclic extension of local fields (of arbritrary
residual characteristic) is computed.

\bigbreak
\centerline{\bf V.  The filtration on $K^\times\!/K^{\times p}$}
\bigskip

In this Part, we make a finer study of the filtration on $K^\times\!/K^{\times
  p}$, which helps us see the precise relationship between th.~{14} and
prop.~{15}.  In fact, we shall give three equivalent ways of looking at this
filtration.  Then we compute a number of examples which illustrate some of our
results and explain the theoretical underpinnings of some results in the
literature.

When the $\Z_p$-module $U_1$ of {\it Einseinheiten\/} is free, a basis can be
found in \citer\fesvost(I,\S6).  One can also find there and in $VI,\S5$ a
basis for the $(\Z/p^r\Z)$-module $U_1/U_1^{p^r}$ when the torsion subgroup of
$U_1$ has order $p^r$.

\bigbreak
{\it The first definition of the filtration on\/
  $K^\times\!/K^{\times p}$.}  Let us put $U_0=\units$ and $\bar
U_n=U_n/(U_n\cap U_0^p)$.  We thus get a filtration on
$K^\times\!/K^{\times p}$ by sub-$\F_p$-spaces
$$
\cdots\subset\bar U_n\subset\cdots\subset\bar U_1\subset\bar U_0\subset
K^\times\!/K^{\times p}.
$$

We have seen that $\bar U_{n}=\{\bar1\}$ if $n>p{e_1}$ (prop.~{30}~; note that
$p{e_1}={e_1}+e$) and that $\bar U_{p{e_1}}$ is an $\Fp$-line if $K^\times$
has an element of order~$p$ (prop.~{33}), otherwise $\bar
U_{p{e_1}}=\{\bar1\}$.  Let us show that the other inclusions $\bar
U_{n+1}\subset\bar U_n$ are an {\it equality\/} if $p|n$ and of {\it
codimension\/} $f$ otherwise.

\th PROPOSITION {42}
\enonce
Let\/ $\zeta\in\bar K^\times$ be an element of order\/ $p$.  Then
$$
\dim_\Fp\bar U_{n}/\bar U_{n+1}=\cases{
0&if\/ $n>p{e_1}$,\cr
1&if\/ $n=p{e_1}$ and\/ $\zeta\in K^\times$,\cr
0&if\/ $n=p{e_1}$ and\/ $\zeta\notin K^\times$,\cr
0&if\/ $n<p{e_1}$ and\/ $p|n$,\cr
f&otherwise.\cr}
$$
(\/Note that the condition ``\thinspace\thinspace$\zeta\in K^\times$'' just
means that\/ $K^\times$ has an element of order~$p$).
\endth 
We have already seen the cases $n>p{e_1}$ (prop.~{30}) and $n=p{e_1}$
(prop.~{29}).  Suppose next that $n=ps$ for some $s<{e_1}$.  Then, $U_{ps}\cap
U_0^p=U_s^p$, whereas $U_{ps+1}\cap U_0^p=U_{s+1}^p$, as follows from
prop.~{29} (cf.~the proof of prop.~{33}).  Thus we have a commutative diagram 
\vskip-20pt
$$
\def\\{\mskip-2\thickmuskip}
\def\droite#1{\\\hfl{#1}{}{8mm}\\}
\diagram{
&&&&&&1\cr
&&&&&&\vfl{}{}{5mm}\cr
1&\rightarrow&U_{s+1}&\droite{}&
U_s&\droite{}&U_s/U_{s+1}&\rightarrow&1\cr
&&\vfl{(\ )^p}{}{5mm}&&\vfl{(\ )^p}{}{5mm}&&\vfl{(\ )^p}{}{5mm}\cr
1&\rightarrow&U_{ps+1}&\droite{}&
U_{ps}&\droite{}&U_{ps}/U_{ps+1}&\rightarrow&1\cr
&&\vfl{}{}{5mm}&&\vfl{}{}{5mm}&&\vfl{}{}{5mm}\cr
&&\bar U_{ps+1}&\droite{}&
\bar U_{ps}&&1\cr
}
$$
both whose rows are exact, as are the columns (prop.~{29}).  Hence the
inclusion $\bar U_{ps+1}\rightarrow\bar U_{ps}$ is an isomorphism by the
serpent lemma.  Consider finally the case when $n=ps+g$ ($s<{e_1}$, $0<g<p$)
is $<p{e_1}$ and prime to~$p$.  Then $U_n\cap U_0^p=U_{s+1}^p$ and
$U_{n+1}\cap U_0^p=U_{s+1}^p$, as follows from prop.~{29} by considerations
similar to the ones above.  The inclusions $U_{s+1}^p\subset U_{n+1}\subset
U_n$ induce the short exact sequence
$$
1\to\bar U_{n+1}\to\bar U_{n}\to U_{n}/U_{n+1}\to1
$$
which proves the proposition because we know that the last 
$\F_p$-space is of dimension $f$.

As a check, the dimension of $\bar U_0=\bar U_1$, according to prop.~{42},
turns out to be ${e_1}.(p-1)f=ef=d$ if $\zeta\notin K^\times$, which is what
it is according to cor.~{32}.

\th COROLLARY {43}
\enonce
Suppose that\/ $K^\times$ has an element\/ $\zeta$ of order\/ $p$.  Then
$\dim_\Fp K^\times\!/K^{\times p}=2+d=2+ef$ and, for $m=pe_1-t$ with
$t\in[0,pe_1[$, 
$$
\dim_\Fp\bar U_{m}=1+\left(t-\left[t\over p\right]\right)\!f.
$$
\endth
When $K^\times$ has no element of order\/ $p$, subtract $1$ from these
dimensions. 

Concretely, let $pe_1-1=m_1>\cdots>m_e=1$ be the $e=(p-1)e_1$ numbers in
$[1,pe_1]$ which are prime to $p$.  When $\zeta\in K$, we have
$$
\dim_\Fp\bar U_{m_i}=1+if\quad (i\in[1,e]).
$$

\bigbreak

{\it The second definition of the filtration on\/ $K^\times\!/K^{\times p}$.}
Because $U_{p{e_1}+1}\subset U_0^p$ (prop.~{30}), the $\Fp$-spaces $\bar U_n$
can be described entirely in terms of the ring $\ogoth/\pgoth^{p{e_1}+1}$, or
rather its group of units $W_0=(\ogoth/\pgoth^{p{e_1}+1})^\times$, which comes
equipped with the filtration
$$
W_n=U_n/U_{p{e_1}+1}=
\Ker((\ogoth/\pgoth^{p{e_1}+1})^\times\rightarrow(\ogoth/\pgoth^{n})^\times)
\quad (n\in[0,p{e_1}+1]).
$$
\vskip-\lastskip
\th PROPOSITION {44}
\enonce
We have\/ $\bar U_n=W_n/(W_n\cap W_0^p)$ for\/ $n\in[0,p{e_1}+1]$.
\endth

This follows from the fact that the image of $(U_n\cap U_0^p)\subset U_n$,
modulo $U_{p{e_1}+1}$, is $(W_n\cap W_0^p)\subset W_n$.  Indeed, the inclusion
$(U_n\cap U_0^p)/U_{p{e_1}+1}\subset (W_n\cap W_0^p)$ is clear.  Conversely,
suppose that $\bar y\in W_0^p$ for some $y\in U_n$~; we have to show that
$y\in U_0^p$.  Let $x\in U_0$ be such that $\bar y=\bar x^p$.  Then $y/x^p\in
U_{p{e_1}+1}$, hence there is a (unique) $\gamma\in U_{{e_1}+1}$ such that
$y/x^p=\gamma^p$.  Consequently, $y\in U_0^p$.  An equivalent way of saying
it~: $W_0^p=U_0^p/U_{p{e_1}+1}$ (prop.~{30}).

\bigbreak

{\it The third definition of the filtration on\/ $K^\times\!/K^{\times p}$.}
Let us come to the definition most commonly used in the literature.  The
following proposition is also true for $l\neq p$, but it is then trivial.

\th PROPOSITION {45}
\enonce
For\/ $x\in\units$ and $n>0$, let $\bar x$ be its image in\/
  $\units\!/\ogoth^{\times p}$ and\/ $\hat x$ the image in\/
  $(\ogoth/\pgoth^{n})^\times$.  Then
$$
\bar x\in\bar U_{n}
\ \Longleftrightarrow\ 
\hat x\in(\ogoth/\pgoth^{n})^{\times p}.
$$
\endth 
Suppose first that $\bar x\in\bar U_{n}$.  Then there is a
$y\in\units$ such that $xy^p\in U_{n}$.  But then $\hat x=\hat y^{-p}$, hence
$\hat x\in(\ogoth/\pgoth^{n})^{\times p}$.  Conversely, suppose that $\hat
x\in(\ogoth/\pgoth^{n})^{\times p}$, and write $\hat x=\hat y^{p}$ for some
$y\in\units$.  Then $xy^{-p}$ belongs to $U_{n}$ and hence $\bar x\in\bar
U_{n}$.

Prop.~{45} says that $\bar U_n=\Ker(\units\!/\ogoth^{\times p}\rightarrow
G_n/G_n^p)$, where $G_n=(\ogoth/\pgoth^n)^\times$.  In the light of the case
$p=2$ of this prop., th.~{14} is equivalent to ``$d_{L|K}\in\bar U_{2e}$'',
whereas prop.~{15} says, writing $\bar U_{2e}=\{1,u\}$, that $d_{L|K}=1$ for
$[L:K]$ odd and $d_{L|K}=u$ for $[L:K]$ even.  The difference is already clear
for $K=\Q_2$, where Hensel's th.~{6} is more precise, in that it specifies when
$D\equiv1\pmod8$ and when $D\equiv5\pmod8$, than Schur's th.~{11}
($D\equiv1\pmod4$).

Let us make the criterion of cor.~{35} explicit for $l=p$, in terms of
reduction modulo $\pgoth^{p{e_1}}$.

\th COROLLARY {46}
\enonce 
A unit\/ $x\in\units$ becomes a\/ $p$-th power in the maximal unramified
extension, or in the maximal tamely ramified extension, if and only if it is
a\/ $p$-th power modulo\/ $\pgoth^{p{e_1}}$.
\endth

\th EXAMPLE {47} (\citer\kraus(p.~376))
\enonce
Let $K=\Q_2(\pi)$, $\pi^3=2$, so that $e=3$ and $\pgoth^{p{e_1}}=4\ogoth$.  A 
unit becomes a square in the maximal unramified extension of\/ $K$ if
and only if it is congruent to one of 
$$
1,\ 1+\pi^2+\pi^4,\ 1+\pi^2+\pi^5,\ 1+\pi^4+\pi^5\pmod4.
$$
\endth
By cor.~{46}, we have to determine the squares in
$(\ogoth/4\ogoth)^\times=U_1/U_6$ or (prop.~{29}) the image of the map $(\
)^2: U_1/U_3\rightarrow U_1/U_6$.  As the above list consists of the squares
of $1$, $1+\pi$, $1+\pi^2$, $1+\pi+\pi^2$, we are done.

\th EXAMPLE {48}
\enonce
A unit in\/ $K=\Q_2(\pi)$, $\pi^3=2$, is the discriminant of an
odd-degree (resp.~even-degree) unramified extension if and only if, up to
squares in $(\ogoth/{\pi^7}\ogoth)^\times$, it is 
$$
\equiv1\ \ (\hbox{resp.}\equiv 1+\pi^6)\ \ \pmod{\pi^7}.
$$
\endth
Indeed, by the second definition of the filtration on
$\units\!/\ogoth^{\times2}$, the group $\bar U_6$ can be identified with
$\{1,1+\pi^6\}$, which, in the notation of that definition, is the kernel
$W_6=\{1,1+\pi^6\}$ of $U_1/U_7\rightarrow U_1/U_6$ modulo its intersection
$\one$ with the subgroup of squares in $U_1/U_7$.

\th EXAMPLE {49}
\enonce
Let\/ $\zeta\in\bar\Qp^\times$ be an element of order\/ $p$.  A unit of\/
$K=\Qp(\zeta)$ becomes a $p$-th power in the maximal tamely ramified extension
of\/ $K$ if, and only if, it~is\/ $\equiv1\pmod{\pgoth^p}$, where\/ $\pgoth$
is the maximal ideal of the ring of integers of\/ $K$.
\endth
Here $e=p-1$ and ${e_1}=1$.  By cor.~{35}, we have to determine the $p$-th
powers in $U_1/U_p$.  But the raising-to-the-exponent-$p$ map $(\ )^p$ takes
$U_1$ to $U_p$ (prop.~{27}), so the only $p$-th power in $U_1/U_p$ is $1$.

\th EXAMPLE {50} \citer\hasse(Kap.~15)
\enonce
$K(\sqrt{-1})$ is the unramified quadratic extension of\/ $K=\Q_2(\sqrt3)$. 
\endth
It suffices (prop.~{16}, prop.~{45}) to show that $-1$ is a square
$\pmod{\pi^4}$, where $\pi=\sqrt3-1$ is a uniformiser of $K$.  Indeed,
we have $-1=(\sqrt3)^2-4$, which implies that
$-1\equiv(1+\pi)^2\pmod{\pi^4}$. 

More generally, the compositum\/ $L_1L_2$ of two linearly disjoint
\hbox{degree-$p$} ramified kummerian extensions $L_1=K(\root p\of {D_1})$,
$L_2=K(\root p\of {D_2})$ is unramified over $L_1$, $L_2$ precisely when the
$\Fp$-plane $D_1D_2$ contains the line $\bar U_{pe_1}$~:

\th EXAMPLE {51}
\enonce
Let\/ $K$ be a finite extension of\/ $\Qp(\zeta)$ ($\zeta^p=1$ but\/
$\zeta\neq1$) and let\/ $D_1$, $D_2$ be two distinct\/ $\Fp$-lines in\/
$K^\times\!/K^{\times p}$, distinct from\/ $\bar U_{pe_1}$, such that the
plane\/ $D_1D_2$ contains\/ $\bar U_{pe_1}$.  Then the compositum\/ $L_1L_2$ is
the unramified degree-$p$ extension of\/ $L_1=K(\root p\of {D_1})$ and of\/ 
$L_2=K(\root p\of {D_2})$.
\endth
This results from the computation of the relative ramification index and
residual degree of the degree-$p^2$ extension $L_1L_2|K$ by multiplicativity
in three different ways, the intermediate extensions being respectively $L_1$,
$L_2$, and the degree-$p$ unramified extension $L$ of $K$, which is contained
in $L_1L_2$ because the plane $D_1D_2$ contains the line $\bar U_{pe_1}$.  The
only possibility for the said relative ramification index and residual degree
is $p$ and $p$, which forces the conclusion.

Specifically, for $K=\Qp(\zeta)$, take $D_1$ to be the $\Fp$-line generated by
$\bar\zeta$ in $K^\times\!/K^{\times p}$, and $D_2$ to be any line in the
plane $D_1\bar U_p$ distinct from the lines $D_1$, $\bar U_p$.  Then
$L_2(\root p\of\zeta)$ is the unramified degree-$p$ extension of $L_2=K(\root
p\of {D_2})$.  Ex.~{50} illustrates this observation for $p=2$.

(Abhyankar's lemma asserts that if $K|\Qp$, $L_1|K$, $L_2|K$ are finite
extensions with $L_1|K$ {\it tame\/} of ramification index dividing the degree
of $L_2|K$, then the extension $L_1L_2|L_2$ is unramified,
cf.~\citer\narkiewicz(p.~236).  One can see that the requirement of tameness
cannot be dispensed with~: take $K$ containing $\zeta$, and $L_1=K(\root p\of
{D_1})$, $L_2=K(\root p\of {D_2})$, where the $\Fp$-lines $D_1$, $D_2$ in
$K^\times\!/K^{\times p}$ are chosen to be different and such that the plane
$D_1D_2$ does {\it not\/} contain the line $\bar U_{pe_1}$, where $(p-1)e_1$
is the absolute ramification index of $K$.)

\bigbreak
\centerline{\bf VI.  A brief history of primary numbers}
\bigskip

We pass in review the problem of characterising, among prime-degree kummerian
extensions of local fields, the unramified one, and of computing the valuation
of the discriminant for the ramified ones.

As in Part~I, there is no pretence at exhaustiveness.  These two Parts have
been added merely to share with the reader some of the wonderful theorems
which we have discovered in the classical literature.  They are presented
roughly in the order in which we came across them~; this explains the
chronological zigzag.

\medbreak 

{\bf Fr{\"o}hlich (1960).}  After the proof of prop.~{16}, we wanted to
 generalise, from quadratic extensions to prime-degree kummerian extensions,
 results of Fr{\"o}hlich on the valuation of the discriminant.  We quote the
 result without explaining his peculiar notation~; a translation into our
 notation is therefore provided.

\th THEOREM {52} \citer\frohlich(p.~24)
\enonce
Let\/ $p$ be a prime lying above\/ $2$.

(i) If $\alpha\equiv0\pmod p$ then $\delta_p=4_p\alpha_p$.

(ii) If $\alpha_p\in U_p$ and if $s$ is the greatest positive integer such
that firstly
$$
4\equiv0\pmod {p^{2s}}
$$
and in the second place
$$
\exists\;\gamma_p\in U_p\quad\hbox{ with }\quad
\alpha_p\equiv\gamma_p^2\pmod{p^{2s}},
$$
then
$$
\delta_p=(2_p\pi^{-1})^2\alpha_p,\quad\hbox { where }
\pi\in\ogoth_p^\times,\ \ (\pi)=p^s.
$$
\endth
(Let $K$ be a finite extension of $\Q_2$ and let $L=K(\sqrt\alpha)$ for some
$\alpha\in K^\times$.  If $\bar\alpha\in\bar U_m$ but $\bar\alpha\not\in\bar
U_{m+1}$ for some $m<2e$ (with the convention that $\bar
U_0=K^\times\!/K^{\times2}$), then $v_K(d_{L|K})=1+2e-m$.  If
$\bar\alpha\in\bar U_{2e}$ but $\bar\alpha\not\in\bar U_{2e+1}$, then
$v_K(d_{L|K})=0$.)




\medbreak 

{\bf Hecke (1923).}  A search in the literature revealed {\it S{\"a}tze\/}
118--120 in Hecke's {\it Vorlesungen}, which crown his treatment of {\it
Allgemeine Arithmetik der Zahlk{\"o}rper\/}.  He is studying the behaviour of
primes in a degree-$l$ kummerian extension $K=k(\root l\of\mu)$ of a number
field $k$ ($l$ being a prime for which $k^\times$ has an element $\zeta$ of
order~$l$\/)~: when does a prime $\pgoth$ of $k$ remain a prime, or become the
$l$-th power of a prime, or split as a product of $l$ (distinct) primes in
$K$~?  {\it Satz\/}~118 says~:

\th THEOREM {53} (\citer\hecke(p.~150))
\enonce 
Es gehe das Primideal\/ $\pgoth$ in der Zahl $\mu$ genau in der Potenz
$\pgoth^a$ auf.  Wenn dann $a$ nicht durch $l$ teilbar ist, so wird $\pgoth$
die $l^{\rm te}$ Potenz eines Primideals in $K$~: $\pgoth={\goth P}^l$.  Wenn
aber $a=0$ ist und\/ $\pgoth$ nicht in $l$ aufgeht, so wird\/ $\pgoth$ in $K$
das Produkt von $l$ vershiedenen Primidealen, falls die Kongruenz
$$
\mu\equiv\xi^l\pmod\pgoth
$$
durch eine ganze Zahl $\xi$ in $k$ l{\"o}sbar ist, dagegen bleibt $\pgoth$
ein Primideal in $K$, falls diese Kongruenz unl{\"o}sbar ist.
\endth
(Let $a$ be the $\pgoth$-adic valuation of $\mu$.  If $l$ does not divide $a$,
then $\pgoth$ becomes the $l$-th power of a prime $\goth P$ in $K$.
When $a=0$ but $\pgoth$ does not divide $l$, then $\pgoth$ splits in $K$ as a
product of $l$ distinct prime ideals if the displayed congruence is solvable
by an integer $\xi$ of $k$, otherwise $\pgoth$ remains a prime in $K$.)

The first part of the corresponding local statement --- for $k|\Qp$ a finite
extension and $l$ a prime for which $k^\times$ has an element of order~$l$ ---
would say that if $\mu\in k^\times\!/k^{\times l}$ but
$\mu\notin\units\!/\ogoth^{\times l}$, then the extension $k(\root l\of\mu)$
is totally ramified.  This is easy, because $\mu$ can be taken to have
valuation~$1$, and the polynomial $T^l-\mu$ is then Eisenstein.

The second part would say that if $l\neq p$ and if $\mu\in\units$, then
$\mu\in\ogoth^{\times l}$ if $\mu$ is an $l$-th power in the residue
field~: this follows from Hensel's lemma.  However, if
$\mu\notin\ogoth^{\times l}$, then the extension $k(\root l\of\mu)$ is
unramified~: this is precisely the content of the case $l\neq p$ of prop.~{16}.

Hecke had the global result for $p=l$ as well~; {\it Satz\/}~119 says~:

\th THEOREM {54} (\citer\hecke(p.~152))
\enonce 
Es sei $\goth l$ ein primfaktor von $1-\zeta$, die darin genau zur $a^{\rm
  ten}$ Potenz aufgeht~: $1-\zeta={\goth l}^a{\goth l}_1$~; es gehe $\goth l$
nicht in $\mu$ auf.  Dann zerf{\"a}llt $\goth l$ in $l$ voneinander
verschiedene Faktoren in $K(\root l\of\mu;k)$, falls die Kongruenz
$$
\mu\equiv\xi^l\pmod{{\goth l}^{al+1}}\leqno{(82)}
$$
durch eine Zahl $\xi$ in $k$ l{\"o}sbar ist.  Es bleibt $\goth l$ auch in $K$
Primideal, wenn zwar die Kongruenz
$$
\mu\equiv\xi^l\pmod{{\goth l}^{al}}\leqno{(83)}
$$
aber nicht $(82)$ l{\"o}sbar ist.  Endlich wird $\goth l$ die $l^{\rm te}$
Potenz eines Primideals in $K$, wenn auch diese Kongruenz $(83)$ unl{\"o}sbar
ist. 
\endth
(Let $\lgoth$ be a prime of $k$ dividing $1-\zeta$ and let $a$ be the
$\lgoth$-adic valuation of $1-\zeta$~; suppose that $\lgoth$ does not divide
$\mu$.  Then $\lgoth$ splits as a product of $l$ distinct primes in $K=k(\root
l\of\mu)$ if the congruence $(82)$ can be solved by an integer $\xi$ of $k$.
If the congruence $(83)$ can be solved, but not $(82)$, then $\lgoth$ remains
a prime in $K$.  Finally, $\lgoth$ becomes the $l$-th power of a prime
in $K$ if the congruence $(83)$ cannot be solved.)


Even Hasse's formulation ({\it Ist\/ $\zeta_0$ eine primitive $l$-te
Einheitswurzel, $\lgoth_0=\lambda_0=1-\zeta_0$ der Primteiler von\/ $l$ im
K{\"o}rper der\/ $l$-ten Einheitswurzeln, und sind\/ $\lgoth^{e_0}$ und\/
$\lgoth^e=\lgoth^{e_0(l-1)}$ die Beitr{\"a}ge von\/ $\lgoth$ zu\/ $\lgoth_0$
und\/ $l=\lgoth_0^{l-1}$, so ist die Bedingung
$$
\alpha\equiv\alpha_0^l\pmod{\lgoth^{e+e_0}}
$$ 
notwendig und hinreichend daf{\"u}r, da\ss\/ $\lgoth$ nicht im F{\"u}hrer
von\/ $k(\root l\of\alpha)$ aufgeht.}) \citer\hassepotenz(p.~232) is global,
as is clear from the context.  

The local version (for $k|\Q_l$ finite, containing an $l$-th root of~$1$)
would say that for $\mu\neq1$ in $\units\!/\ogoth^{\times l}$, the extension
$k(\root l\of\mu)$ is totally ramified, unless $\mu\in\bar U_{l{e_1}}$, in
which case it is unramified~: this is the case $l=p$ of prop.~{16}.  The local
version of the first part of th.~{54} would say that for $\mu\in\units$, the
condition ``\thinspace$\mu\in\ogoth^{\times l}$\thinspace'' is equivalent to
``\thinspace$\bar\mu\in(\ogoth/\pgoth^{l{e_1}+1})^{\times l}$\thinspace''~:
this follows from props.~{29} and~{44}, in view of the fact that $a=e_1$
(prop.~{23}).  

In view of prop.~{30} and~{45}, we see that prop.~{16} is the precise local
counterpart of the global th.~{53} and~{54}, which are its immediate
corollaries.  Let us mention that th.~{53} and~{54} were first proved by
Furtw{\"a}ngler \citer\furtwangler() respectively as his {\it S{\"a}tze\/}~3
and~4~; the results go back in part to Kummer.  Hecke singles out the
following consequence ({\it Satz\/}~120)~: 

\th THEOREM {55} (\citer\hecke(p.~154))
\enonce
Die Relativdiskriminante von\/ $K(\root l\of\mu;k)$ in bezug auf\/ $k$ ist
dann und nur dann gleich\/ $1$, wenn $\mu$ die $l^{\rm te}$ Potenz eines
Ideals in $k$ ist, und gleichzeitig, sofern dann\/ $\mu$ zu\/ $l$ teilerfremd
gew{\"a}hlt wird, die Kongruenz\/ $\mu\equiv\xi^l\pmod{(1-\xi)^l}$ durch eine
Zahl $\xi$ in $k$ l{\"o}sbar ist.
\endth
(The relative discriminant of\/ $K=k(\root l\of\mu)$ over\/ $k$ equals\/ $1$
precisely when\/ $\mu$ is the $l$-th power of an ideal in\/ $k$ and moreover,
when\/ $\mu$ is prime to\/ $l$, the congruence
$\mu\equiv\xi^l\pmod{(1-\zeta)^l}$ admits a solution\/ $\xi$ in\/ $k$.)


\medbreak

{\bf Hilbert (1897).}  Somewhat later we found that the case $l\neq2$ of
Hecke's theorems is a generalisation of a part of {\it Satz\/}~148 in
Hilbert's {\it Zahlbericht\/}, which treats the case $k=\Q(\zeta)$~; the other
part computes the valuation of the discriminant for ramifid prime-degree
kummerian extensions, again for this special base field.  Hilbert's notation
for $\Q(\zeta)$ is $k(\zeta)$~; he takes an integer $\mu$ in this field which
is not an $l$-th power.  {\it Satz\/}~148 says~:

\th THEOREM {56} (\citer\zahlbericht(p.~251))
\enonce
Es werde $\lambda=1-\zeta$ und\/ $\lgoth=(\lambda)$ gesetzt.  Geht ein von\/
$\lgoth$ verschiedenes Primideal\/ $\pgoth$ des Kreisk{\"o}rpers\/ $k(\zeta)$
in der Zahl\/ $\mu$ genau zur\/ $e$-ten Potenz auf, so enth{\"a}hlt, wenn der
Exponent\/ $e$ zu\/ $l$ prim ist, die Relativdiskriminante des durch\/
$M=\root l\of\mu$ und\/ $\zeta$ bestimmten Kummerschen K{\"o}rpers in bezug
auf\/ $k(\zeta)$ genau die Potenz\/ $\pgoth^{l-1}$ von\/ $\pgoth$ als Faktor.
Ist dagegen der Exponent\/ $e$ ein vielfaches von\/ $l$, so f{\"a}llt diese
Relativdiskriminante prim zu\/ $\pgoth$ aus.

Was das Primideal\/ $\lgoth$ betrifft, so k{\"o}nnen wir zun{\"a}chst den
Umstand ausschli{\ss}en, da{\ss} die Zahl\/ $\mu$ durch\/ $\lgoth$ teilbar ist
und dabei\/ $\lgoth$ genau in einer solchen Potenz enth{\"a}lt, deren Exponent
ein Vielfaches von\/ $l$ ist~; denn alsdann k{\"o}nnte der Zahl\/ $\mu$ sofort
durch eine zu\/ $\lgoth$ prime Zahl\/ $\mu^*$ ersetzt werden, so da{\ss\/}
$k(\root l\of{\mu^*},\zeta)$ derselbe K{\"o}rper wie $k(\root
l\of{\mu},\zeta)$ ist.  Unter Ausschlu{\ss\/} des genannten Umstandes haben
wir die zwei m{\"o}glichen F{\"a}lle, da{\ss\/} $\mu$ genau eine Potez von\/
$\lgoth$ enth{\"a}lt, deren Exponent zu\/ $l$ prim ist, oder da{\ss\/} $\mu$
nicht durch\/ $\lgoth$ teilbar ist.  Im\/ {\rm ersteren} Falle ist die
Relativdiskriminante von\/ $k(\root l\of{\mu},\zeta)$ in bezug auf\/
$k(\zeta)$ genau durch die Potenz\/ $\lgoth^{l^2-1}$ teilbar.  Im\/ {\rm
zweiten} Falle sei $m$ der h{\"o}chste Exponent $\le l$, f{\"u}r den es eine
Zahl\/ $\alpha$ in\/ $k(\zeta)$ gibt, so da{\ss\/} $\mu\equiv\alpha^l$ nach\/
$\lgoth^m$ ausf{\"a}llt.  Jene Relativdiskriminante ist dann im Falle $m=l$
zu\/ $\lgoth$ prim~; sie ist dagengen im Falle $m<l$ genau durch die Potenz\/
$\lgoth^{(l-1)(l-m+1)}$ von\/ $\lgoth$ teilbar.
\endth
(Put $\lgoth=(1-\zeta)$.  If for some prime $\pgoth\neq\lgoth$ of $k(\zeta)$,
the number $\mu$ is divisible precisely by the $e$-th power of $\pgoth$, and
if $e$ is prime to $l$, then the relative discriminant of $k(\root
l\of\mu,\zeta)$ over $k(\zeta)$ is divisible precisely by $\pgoth^{l-1}$.  If,
however, $l$ divides $e$, then the relative discriminant is prime to $\pgoth$.

As for the prime $\lgoth$, we may exclude the case in which $\mu$ is divisible
by a power of $\lgoth$ whose exponent is a multiple of $l$, for in this case
we can replace $\mu$ by $\mu^*$ which is prime to $\lgoth$ and such that
$k(\root l\of{\mu^*},\zeta)=k(\root l\of{\mu},\zeta)$.

Leaving aside this case, there are two possibilities~: either $\mu$ is
divisible by a power of $\lgoth$ whose exponent is prime to $l$, or $\mu$ is
prime to $\lgoth$.  In the {\it first\/} case, the relative discriminant of
$k(\root l\of\mu,\zeta)$ over $k(\zeta)$ is divisible precisely by
$\lgoth^{l^2-1}$.  In the {\it second\/} case, let $m\le l$ be the highest
exponent for which there is an integer $\alpha$ in $k(\zeta)$ such that
$\mu\equiv\alpha^l\pmod{\lgoth^m}$.  The relative discriminant is then prime
to $\lgoth$ if $m=l$, and divisible precisely by $\lgoth^{(l-1)(1+l-m)}$ if
$m<l$).  

The proof in the {\it Zahlbericht\/} needs {\it tief{}liegende S{\"a}tze mit
schwierigen Beweisen}, in Hensel's words \citer\henselsolvable(p.~200).
Observe that in the first case ($v_{\lgoth}(\mu)$ is prime to~$l$) Hilbert
could have defined $m=0$.

\bigskip\centerline{***}\bigbreak

But Stickelberger's th.~{1} is more closely related to {\it Satz\/}~96, which 
treats the prime $2$ missing from {\it Satz\/}~148 (th.~{56}).  Hilbert is
working with the quadratic field $k=\Q(\sqrt m)$ for some squarefree integer
$m\neq1$.  Following Dedekind (cf.~th.~{10}), he first determines the ring of
integers of $k$ and its discriminant.  {\it Satz\/}~95 says~:

\th THEOREM {57} (\citer\zahlbericht(p.~157))
\enonce
Eine Basis der quadratischen K{\"o}rpers\/ $k$ bilden die Zahlen $1$,
$\omega$, wenn
$$
\omega={1+\sqrt m\over2},\ \hbox{ bzw. }\ \omega=\sqrt m
$$
genommen wird, je nachdem die Zahl $m\equiv1$ nach $4$ oder nicht.  Die
Diskriminante von $k$ ist, entsprechend diesen zwei F{\"a}llen,
$$
d=m,\ \hbox{ bzw. }\ d=4m.
$$
\endth
($\omega$ being defined as above according as $m\equiv1\pmod4$ or not, $1$,
$\omega$ is a $\Z$-basis of the ring of integers of $k$, and the discriminant 
of $k$ is $m$ or $4m$ respectively.)

The splitting in $k$ of rational primes is treated in {\it Satz\/}~96, which
says~:

\th THEOREM {58} (\citer\zahlbericht(p.~158))
\enonce
Jede in\/ $d$ aufgehende rationale Primzahl\/ $l$ ist gleich dem Quadrat eines
Primideals in\/ $k$.  Jede ungerade, in\/ $d$ nicht aufgehende rationale
Primzahl\/ $p$ zerf{\"a}llt in\/ $k$ entweder in das Produkt zweier
verschiedener, zu einander konjugierter Primideale ersten Grades\/ $\pgoth$
und\/ $\pgoth'$ oder stellt selbst ein Primideal zweiten Grades vor, je
nachdem\/ $d$ quadratischer Rest oder Nichtrest f{\"u}r\/ $p$ ist.  Die
Primzahl\/ $2$ ist im Falle\/ $m\equiv1$ nach\/ $4$ in\/ $k$ in ein Produkt
zweier voneinander verschiedener konjugierter Primideale zerlegbar oder selber
Primideal, je nachdem\/ $m\equiv1$ oder\/ $\equiv5$ nach\/ $8$ ausf{\"a}llt.
\endth
(If a prime number $l$ divides $d$, then it becomes the square of a prime
ideal in $k$.  An odd prime number $p$ which does not divide $d$ splits as the
product of two distinct degree-$1$ mutually conjugate prime ideals or becomes
a degree-$2$ prime ideal in $k$ according as $d$ is a quadratic residue or not
$\pmod p$.  In the case $m\equiv1\pmod4$, the prime $2$ splits into the
product of two distinct conjugate prime ideals or remains a prime in $k$
according as $m\equiv1$ or\/ $\equiv5\pmod8$.)

Notice that, except for the determination of the ring of integers, Hilbert
could have made these two theorems a part of his {\it Satz}~148 if he had
allowed the prime $l$ there to equal $2$.

\medbreak

{\bf Hasse (1927).}  Thus Hecke's {\it S{\"a}tze\/}~118 and~119 generalise a
part of Hilbert's {\it S{\"a}tze\/}~95, 96 and~148, but leave out the
computation of the valuation of the discriminant for ramified kummerian
extensions of degree equal to the residual characteristic.  Our conjectural
answer, which amounted to the equality $t+m=le_1$ involving the ``level'' of a
line $D\neq\bar U_{le_1}$ in $K^\times\!/K^{\times l}$ (the integer $m$ such
that $D\subset\bar U_m$ but $D\not\subset\bar U_{m+1}$) and the break $t$ in
the filtration of $\Gal(L|K)$, $L=K(\root l\of D)$, by higher ramification
groups as explained below, was found proved in {\it Satz\/}~10 (p.~266) of
Hasse's {\it Klassenk{\"o}rperbericht\/} \citer\klassenkorperbericht(), in the
form of the equality $u+v=e_0l$.  He is working with number fields, but the
question is purely local, so we formulate it for local fields.  See Part~VII.

\medbreak 


\bigskip\centerline{***}\bigbreak

{\bf Hensel (1921).} Our most recent find, and the one most relevant to our
proof of prop.~{16}, are three related papers \citer\henselkummer(),
\citer\henselsolvable(), \citer\henselmultii() of Hensel.  We merely reproduce
his arguments in the ($l,l)$-case, leaving the problem of interpretation to the
reader. In the first one (pp.~117--8), he denotes by $\zeta$ a primitive
$l$-th root of~$1$, where $l$ is a prime $\neq2$, and considers an irreducible
equation $x^l-A=0$ over the cylcotomic field $K(\zeta)$, which we would denote
by $\Q(\zeta)$, and in which $l\sim\lgoth^{l-1}$, with $\lgoth$ the prime
ideal generated by $\lambda=1-\zeta$.  In $K(\lgoth)=\Q_l(\zeta)$, or rather
in its multiplicative group modulo $l$-th powers, he writes
$$
A=\lambda^a(1-\lambda)^{c_1}(1-\lambda^2)^{c_2}\ldots(1-\lambda^l)^{c_l}
\quad(a, c_i<l),\leqno{(16.^{\rm a})}
$$
meaning $(a,c_1,c_2,\ldots,c_l)\in(\Z/l\Z)^{1+l}$. If
$(a,c_1,c_2,\ldots,c_l)=0$, then $A$ is an $l$-th power in $\Ql(\zeta)$ and 
therefore 
$$
\lgoth\sim{\goth L}_1{\goth L}_2\ldots{\goth L}_l\leqno{(18.)}
$$
is a product of $l$ distinct prime ideals in $K(x)=\Q(\zeta,x)$.  Secondly, if
$a>0$ (i.e.~if $a\neq0$), then one may take $a=1$, and it follows that
$$
\lgoth\sim{\goth L}_1^l.\leqno{(18.^{\rm a})}
$$
Thirdly, if $a=0$ and if $c_i\neq0$ for some $i\in[1,l]$, let $m$ be the
smallest integer such that $c_m\neq0$~; we may then assume that $c_m=1$.  If
$m<l$, choose $m'$ and $l'$ such that $mm'+ll'=1$~; then
$$
\Lambda=\lambda^{l'}(1-x)^{m'},\leqno{(20.)}
$$
which is the same as the second $\Omega$ in proof of {\it Satz\/}~148 in the
{\it Zahlbericht}, is a primitive element of $K(x)$ whose norm
$n(\Lambda)=\lambda(1+\cdots)$ has valuation~$1$ and therefore, as in the
second case ($a=1$),
$$
\lgoth\sim{\goth L}_1^l.\leqno{(22.)}
$$
Finally, if $m=l$, then the equation $x^l-A=0$ takes the simple form
$$
x^l=1-\lambda^l.\leqno{(23.)}
$$
Hensel puts $\displaystyle{1-x\over\lambda}=\xi$ (this too can be found in
Hilbert), so that the new primitive element $\xi$ satisfies 
$$
f(\xi)=\xi^l-{l_1\over\lambda}\xi^{l-1}+{l_2\over\lambda^2}\xi^{l-2}-\cdots
+{l\over\lambda^{l-1}}\xi-1=0.\leqno{(23.^{\rm a})}
$$
Reducing modulo $\lambda$, and noting that $-l=\lambda^{l-1}(1+\cdots)$, he
gets the congruence
$$
f(\xi)\equiv\xi^l-\xi-1\pmod\lgoth\leqno{(23.^{\rm b})}
$$
Let ${\goth L}_1$ be a prime factor of $\lgoth$ in $K(x)$ and $\xi_0$ a root
of $f(x)=0$ in the {\it Bereich\/} of ${\goth L}_1$ (not ${\goth P}_1$, which
is a {\it Druckfehler\/})~; it follows from $(23.^{\rm b})$ that
$$
\xi_0^l\equiv\xi_0+1\pmod{{\goth L}_1},
$$
and therefore, for $i=1,2,\ldots, l-1$, 
$$
\xi_0^{l^i}\equiv\xi_0+i\pmod{{\goth L}_1}.
$$
As the $l$ conjugates $\xi_0, \xi_0^l,\ldots,\xi_0^{l^{l-1}}$ of $\xi_0$ are
units incongurent modulo ${\goth L}_1$, congruent respectively to $\xi_0,
\xi_0+1,\ldots,\xi_0+(l-1)$, the relative degree of ${\goth L}_1$ is $l$, and
therefore
$$
\lgoth\sim{\goth L}_1.\leqno{(24.)}
$$ (Note that Artin and Schreier could have drawn the inspiration for their
theory from this proof~; they were instead inspired by the archimedean prime
$\R$, which is equally worthy of our contemplation.  Note also that, in giving
this local proof of a part of Hilbert's {\it Satz\/}~148, Hensel missed the
opportunity of proving the other part by computing the relative discriminant
and the ring of integers, even though he had the uniformisers $\root
l\of\lambda$ and $(20.)$ in the two ramified cases $a=1$, resp.~$a=0$, $m<l$.
Finally, by taking $l\neq2$, he cannot connect the unramified case $m=l$, when
$l=2$, with his earlier result (th.~{6}) on local discriminants).

In the second paper (\citer\henselsolvable()), although the formulation of the
theorem is still global, the above local proof is carried over (pp.~207--8) to
the case of any finite extension $K$ of $\Q(\zeta)$ and $\lgoth\,|\,l$ a prime
of $K$ of residual degree~$f$.  Hensel notes that $K(\lgoth)$ contains
(cf.~prop.~{24}) the number
$$
\Lambda=\root l-1\of{-l}.
$$ 
As he shows in \citer\henselmultii(), (cf.~prop.~{30} and~{33}, or
prop.~{42}), one is reduced to considering the equation
$x^l=1-\xi_0\Lambda^l$, with $\xi_0$ a unit of $K(\lgoth)$ whose trace
$$
s_0=\xi_0+\xi_0^l+\xi_0^{l^2}+\cdots+\xi_0^{l^{f-1}}\leqno{(11.^{\rm a})}
$$
is not divisible by $\lgoth$ and hence congruent to one of
$1,2,\ldots,l-1\pmod\lgoth$. (The meaning of $\xi_0$ is different in the two
papers.)  Then the $l$ conjugates 
$$
\xi,\;\xi^{l^f},\;\xi^{l^{2f}},\;\ldots,\;\xi^{l^{(l-1)f}}\leqno{(9)}
$$ 
of the new primitive element $\displaystyle\xi={1-x\over\Lambda}$ 
are units incongurent modulo ${\goth L}$, where ${\goth L}$ is a prime divisor
of $\lgoth$ in $K(x)$.  Indeed, dividing
$$
(\xi\Lambda-1)^l=\xi^l\Lambda^l
-l_1\xi^{l-1}\Lambda^{l-1}+l_2\xi^{l-2}\Lambda^{l-2}-\cdots+l\xi\Lambda-1
=\xi_0\Lambda^l-1
$$
throughout by $\Lambda^l$ and recalling that $\Lambda^{l-1}=-l$, he gets the
equation 
$$
f(\xi)=\xi^l-{l\over\Lambda}\xi^{l-1}+{l_2\over\Lambda^2}\xi^{l-2}-\cdots
-\xi-\xi_0=0.\leqno{(12.^{\rm a})}
$$
Reducing modulo ${\goth L}$ and noting that the coefficients of $\xi^{l-1}$,
$\xi^{l-2}$, $\ldots$, $\xi^{2}$ are divisible by $\Lambda$ and hence by
${\goth L}$, he gets the congruence
$$
\xi^l\equiv\xi+\xi_0\pmod{\goth L}.\leqno{(12.^{\rm b})}
$$
Successively raising this to the exponent $l$, he gets a series of
congruences 
$$
{\eqalign{
\xi^{l^2}&\equiv\xi+\xi_0+\xi_0^{l}\cr
\xi^{l^3}&\equiv\xi+\xi_0+\xi_0^{l}+\xi_0^{l^2}\cr
\multispan2\dotfill\qquad\ \cr
\xi^{l^f}&\equiv\xi+\xi_0+\xi_0^l+\cdots+\xi_0^{l^{f-1}}\cr
&\equiv\xi+s_0.\cr
}}\pmod\lgoth$$
This implies that the $l$ powers $(9)$ are congruent modulo ${\goth L}$
respectively to
$$
\xi,\;\xi+s_0,\;\xi+2s_0,\;\ldots,\;\xi+(l-1)s_0,
$$ 
({\it not\/} to $\xi,\;\xi+\xi_0,\;\xi+2\xi_0,\;\ldots,\;\xi+(l-1)\xi_0$, as
in the original).  In other words, the $l$ powers $(9)$ are distinct modulo
${\goth L}$.  Hence the relative discriminant is not divisible by $l$,
because the $l$ roots $\xi+is_0$ (the original has $\xi+i\xi_0$),
$i\in[0,l[$, are distinct modulo ${\goth L}$. 

This is Hensel's local proof of Furtw{\"a}ngler's {\it Satz\/}~4 ($=$ Hecke's
{\it Satz\/}~119)~; of course, he also had a local proof of {\it Satz\/}~3
(resp.~118), which treats the easier case $l\neq p$.  Curiously, Hensel does
not cite Furtw{\"a}ngler~; it is left to Hasse to do so in his review in the
{\it Jahrbuch\/} [JFM~48.1170.01].   Equally curiously, Hasse
\citer\klassenkorperbericht({\it Satz\/}~9) refers to Hecke for the theorem
and its proof, not to Furtw{\"a}ngler and Hensel.

\bigskip\centerline{***}\bigbreak


{\bf Eisenstein (1850)}.  The roots of Hensel's paper \citer\henselkummer()
can be traced more than 70~years back to Eisenstein's seminal work
\citer\eisenstein() on the $\lambda$-tic reciprocity law for a prime number
$\lambda\neq2$.  He takes a primitive $\lambda$-th root $\zeta$ of $1$ and
sets $\eta=1-\zeta$.  For two mutually prime ``complex integral numbers'' $A$,
$B$ ($\in\Z[\zeta]$), he defines a $\lambda$-th root $(A/B)$ of~1, ``the
$\lambda$-tic character of $A$ modulo $B$'', using Kummer's ideal prime
divisors. (For an {\it ideal\/} prime $\goth b$ other than $\eta\Z[\zeta]$,
define $(A/{\goth b})$ by the congruence $(A/{\goth b})\equiv A^{N{\goth
b}-1\over\lambda}\pmod{\goth b}$, where $N{\goth b}=\Card(\Z[\zeta]/{\goth
b})$ is the norm of ${\goth b}$~; extend the definition to {\it numbers\/} $B$
by multiplicativity).  He investigates, when $A$ and $B$ are not divisible by
$\eta$, the ratio
$$
\left({A\over B}\right):\left({B\over A}\right)=\zeta^{\varphi(A,B)}
\qquad\qquad\left(\varphi(A,B)\in\Z/\lambda\Z\right)
$$ 
(which, if $\lambda$ had been $2$ and consequently $\zeta=-1$ and $(A/B)$ the
quadratic character of $A$ modulo $B$, would have been given by
$\varphi(A,B)={A-1\over2}{B-1\over2}$ for $A,B\in\Z$ odd, mutually prime,
and positive~: ``quadratic reciprocity'').  First he shows that
$A\pmod{\eta^{\lambda+1}}$ can be taken to be of the form
$$
A\equiv g^\alpha(1-k_1\eta)^{\alpha_1}(1-k_2\eta^2)^{\alpha_2}\cdots
(1-k_\lambda\eta^\lambda)^{\alpha_\lambda}\pmod{\eta^{\lambda+1}}
$$ 
where, in our language, a generator $g$ of $\F_\lambda^\times$ and integers
$k_i\in\Z[\zeta]$ prime to $\eta$ are fixed, and the exponents
$\alpha\in\Z/(\lambda-1)\Z$, $\alpha_i\in\Z/\lambda\Z$ vary with $A$.  He
concludes from this that $\varphi(A,B)=\varphi(A',B')$ if $A\equiv A'$ and
$B\equiv B'$ $\pmod{\eta^{\lambda+1}}$ or, what comes to the same,
$\pmod{\lambda\eta^{2}}$.  He remarks that as far as the $\lambda$-tic
character is concerned, we may take $\alpha=0$ and $k_i=1$.  One is thus
reduced to considering the case of relatively prime $A=1-\eta^\mu$,
$B=1-\eta^\nu$ and to the determination of
$\varepsilon_{\mu,\nu}=\zeta^{\varphi(1-\eta^\mu,1-\eta^\nu)}$. 

Next, he considers $(\eta/A)$, which, by multiplicativity, reduces to the
cases $A=1-\eta^\mu$~; what is most relevant for us is the basic formula
$$
\left({\eta\over1-\eta^\mu}\right)=
\cases{1\quad\hbox{if }\mu\neq\lambda\phantom{.}\cr
\zeta\quad\hbox{if }\mu=\lambda.\cr}
\leqno{(8.)}
$$
It follows that, for $A\equiv(1-\eta)^{\alpha_1}(1-\eta^2)^{\alpha_2}\cdots
(1-\eta^\lambda)^{\alpha_\lambda}\pmod{\eta^{\lambda+1}}$, one has more
generally 
$$
\left({\eta\over A}\right)=\zeta^{\alpha_\lambda}.\leqno{(10.)}
$$ 
This law may be considered as a remote ancestor of our prop.~{38}, just as
Eisenstein's analysis $\pmod{\eta^{\lambda+1}}$ may be said to have prefigured
Hensel's equation $(16.^{\rm a})$.  Indeed, Eisenstein sums up the main
achievement of his memoir by saying that $\sigma=\lambda+1$ is the smallest
exponent --- if smallest exponent there is --- such that
$\varphi(A,B)=\varphi(A',B')$ whenever $A,B\equiv A',B'\pmod{\eta^{\sigma}}$.

Local arithmetic makes it possible for us not only to understand all this
wizardry but also {\it sozusagen\/} to anticipate it.

\medskip
{\it 
\rightline{Most of us know the parents or grandparents we come from.}
\rightline{ But we go back and back, forever~; we go back all of us} 
\rightline{to the very beginning~; in our  blood and bone and brain}
\rightline{we carry the memories of thousands of beings.}
\rightline{\rm --- V.\ S.\ Naipaul \citer\naipaul(p.~9).} 
}

\bigskip\centerline{***}\bigbreak

It is amusing to note that the shortest path between the global results of
Stickelberger (th.~{1} and~{2}) on the one hand, and of Hilbert (th.~{56}
and~{58}) and Furtw{\"a}ngler (Hecke's th.~{53} and~{54}) on the other, passes
through purely local results (th.~{5} and~{6}, prop.~{15} and~{16})~: the
statements are different globally, but they are the same locally.

{\it $l^n$-primary numbers}\pointir In a finite extension $K$ of $\Q_l$
containing a primitive $l^n$-th root of~$1$, an element $\alpha\in K^\times$
is called $l^n$-{\it primary\/} if the extension $K(\!\root l^n\of\alpha)|K$
is unramified.  These numbers have been characterised by Hasse~(1936) using
the theory of Witt vectors, and also by Shafarevich (1950) in his work on the
general reciprocity law~; see \citer\fesvost(VI,\S4) for a comprehensive
presentation.

\bigbreak
\centerline{\bf VII.  The valuation of the discriminant}
\bigskip

Let\/ $p$, $l$ be prime numbers and let\/ $K$ be a finite extension of\/
$\Qp$.  Unless otherwise stated (as in prop.~{63}), we assume that $K^\times$
has an element $\zeta$ of order\/ $l$.  Denote by $k$ the residue
field of\/ $K$.  For any extension $F$ of $\Qp$, we denote by $v_F$ the
normalised valuation of $F$, so that $v_F(\pi)=1$ if $\pi$ is a generator of
the maximal ideal of the ring of integers of $F$.

\th PROPOSITION {59}
\enonce
Suppose that $p\neq l$.  Let\/ $D\neq k^\times\!/k^{\times l}$ be an
$\F_l$-line in\/ $K^\times\!/K^{\times l}$ and let\/ $L=K(\root l\of D)$ be
the (ramified) degree-$l$ cyclic extension of\/ $K$ corresponding to\/ $D$.
Then\/ $v_K(d_{L|K})=l-1$.
\endth
\proof This is clear since $L|K$ is a totally ramified extension of degree~$l$
prime to $p$ (tame ramification).  What is needed is not so much that $l$ be
prime, but that it be prime to $p$.  Cf.~footnote to th.~{10}.  

Hasse uses Hilbert's theory of higher ramification groups (\citer\hilbert(),
\citer\zahlbericht(p.~140)), to compute the valuation of the discriminant in
the case $l=p$.  (To be consistent with Hilbert, Hecke and Hasse, we denote
the prime $p=l$ by $l$, not by $p$\/).  Let $L|K$ be a finite {\it
galoisian\/} extension and $G=\Gal(L|K)$, where $K$ is {\it any\/} finite
extension of $\Ql$.  The ring of integers $\ogoth_L$ of $L$ and its maximal
ideal $\lgoth_L$ are stable under the action of $G$~; there is thus an induced
action on $\ogoth_L/\lgoth_L^{n+1}$ for every integer $n\in[-1,+\infty[$.

Define $G_n$ to be the subgroup consisting of those $\sigma\in G$ which act
trivially on $\ogoth_L/\lgoth_L^{n+1}$~; we have $G_{-1}=G$, and $G_0$ is the
inertia subgroup~: the extension $L_0=L^{G_0}$ is unramified over $K$, whereas
$L$ is totally ramified over $L_0$.  Also, $G_1$ is the (unique) maximal
sub-$l$-group of $G_0$~: the extension $L_1=L^{G_1}$ of $L_0$ is (totally but)
tamely ramified (and hence cyclic of degree dividing $q-1$, where
$q=\Card\ogoth_L/\lgoth_L$), whereas $L$ is a (totally ramified) $l$-extension
of $L_1$.  The decreasing filtration $(G_n)_{n\in[-1,+\infty[}$ is exhaustive
and separated.  The valuation of the different of $L|K$ is given by
$$
v_L({\goth D}_{L|K})=\sum_{n\in[0,+\infty[}(\Card G_n-1),\leqno{(10)} 
$$ \citer\klassenkorperbericht(p.~249) which also equals $v_K(d_{L|K})$ when
$L|K$ is totally ramified.  This filtration in the lower numbering is
compatible with the passage to a subgroup, but {\it not\/} with the passage to
a quotient.  

The problem of computing the ramification filtration on a quotient of $G$ was
first solved by Herbrand \citer\herram()~; one needs to convert the filtration
$(G_n)_{n\in[-1,+\infty[}$ to the upper numbering $(G^t)_{t\in[-1,+\infty[}$,
defined for any {\it real\/} $t$, and then take the quotient.  See
\citer\serre(ch.~IV) for the details, where the upper-numbering filtration on
$\Gal(\Ql(\zeta_{l^n})\,|\,\Ql)=(\Z/l^n\Z)^\times$ is also shown, using
prop.~{23}, to be the quotient of the natural filtration of $\Z_l^\times$ by
units of various levels.

Let us now revert to a $K$ which contains a primitive $l$-th root of $1$ and
denote by $e=(l-1)e_1$ the ramification index of\/ $K$ over $\Ql$.

\th PROPOSITION {60}
\enonce
For an\/ $\F_l$-line\/ $D\neq\bar U_{le_1}$ in\/ $K^\times\!/K^{\times l}$,
let\/ $L=K(\root l\of D)$ be the (ramified) degree-$l$ cyclic extension of\/
$K$ corresponding to\/ $D$.  Let\/ $m$ be the integer such that\/
$D\subset\bar U_m$ but\/ $D\not\subset\bar U_{m+1}$, with the convention
that\/ $\bar U_0=K^\times\!/K^{\times l}$.  Then\/
$v_K(d_{L|K})=(l-1)(1+le_1-m)$.
\endth
\proof Suppose first that $m=0$.  Then $D$ is generated by the class modulo
$K^{\times l}$ of some uniformiser $\lambda$ of $K$, and $L$ is defined by
$g=T^l-\lambda$, which is Eisenstein.  As $g'=lT^{l-1}$, the exponent of the
different ${\goth D}_{L|K}$ is $v_L(l)+(l-1)=le+l-1=(l-1)(1+le_1)$, and, as
$L|K$ is totally ramified, this is also the valuation $v_K(d_{L|K})$ of the
discriminant $d_{L|K}$.

Suppose next that $m>0$~; then $m<le_1$ (because $D\neq\bar U_{le_1}$) and $m$
is prime to $l$ (prop.~{42}).  Let $(G_i)_{i\in[-1,+\infty[}$ be the
ramification filtration~; we have $G=G_{-1}=G_0=G_1=\cdots=G_t$ but
$G_{t+1}=\{\Id_L\}$ for some integer $t$ which is strictly positive because
$L|K$ is wildly ramified.  As $v_K(d_{L|K})=(l-1)(1+t)$
(cf. \citer\serre(prop.~IV.4)), it is enough to show that $t=le_1-m$.

The line $D$ is generated by the class modulo $K^{\times l}$ of some unit
$\mu$ of $K$ and there is a unit $\xi$ such that $v_K(\xi^l-\mu)=m$
(prop.~{45}).  Fix a root $\root l\of\mu$ of $T^l-\mu$ in $L$.  As
$N_{L|K}(\xi-\root l\of\mu)=\xi^l-\mu$ and as $L|K$ is totally ramified, we
also have $v_L(\xi-\root l\of\mu)=m$.

Let $x,y\in\Z$ be such that $mx+ly=1$~; the element $\Lambda=(\xi-\root
l\of\mu)^x\lambda^y$ is then a uniformiser of $L$ for any fixed but arbitrary
uniformiser $\lambda$ of $K$.  Taking  a generator $\sigma\in G$ 
and writing $\sigma(\root l\of\mu)=\zeta.\root l\of\mu$ for some order-$l$
element $\zeta$ of $K^\times$, we have
$$
{\sigma(\Lambda)\over\Lambda}
=\left({\xi-\zeta\root l\of\mu\over\xi-\root l\of\mu}\right)^{\!x}
=\left(1+{(1-\zeta)\root l\of\mu\over\xi-\root l\of\mu}\right)^{\!x}
=(1+\alpha)^x,
$$
defining $\alpha$.  Hence
$\sigma(\Lambda)/\Lambda\equiv1+x\alpha\pmod{\Lambda^{le_1-m+1}}$.  As
$v_L(\alpha)=le_1-m$, and as $x$ is prime to $l$, this means that $\sigma\in
G_{le_1-m}$ but $\sigma\notin G_{le_1-m+1}$, and hence $t=le_1-m$, as desired
(cf.~\citer\klassenkorperbericht(p.~266)).

Notice that when $m=0$, our direct computation of $v_K(d_{L|K})$ in this case
shows that $t=le_1$. Thus, $t+m=le_1$ for all ramified degree-$l$ kummerian
extensions.

(This allows one to compute the valuation of the discriminant of
$\Ql(\zeta_{l^n})\,|\,\Ql$ by induction, without invoking the fact that the
ramificaton filtration is a quotient of the filtration on $\Z_l^\times$ by
units of various levels.)

\th COROLLARY {61}
\enonce
For an\/ $\F_l$-line\/ $D\subset K^\times\!/K^{\times l}$, let\/ $L=K(\root
l\of D)$ and let\/ $\ogoth_L$ be the ring of integers of\/ $L$.  If\/
$D\not\subset\bar U_1$, then there is a uniformiser $\mu$ such that\/
$\bar\mu\in D$~; in this case, $\root l\of\mu$ is a uniformiser of\/ $L$ and
$\ogoth_L=\ogoth[\root l\of\mu]$.  Suppose next that\/ $D\subset\bar U_1$ and
let $\mu\in\units$ be a unit such that\/ $\bar\mu$ generates $D$.  If there is
an\/ $m<le_1$ such that\/ $\bar\mu\in\bar U_m$ but\/ $\bar\mu\notin\bar
U_{m+1}$, then\/ $m$ is prime to\/ $l$~; for any\/ $x,y\in\Z$ such that\/
$mx+ly=1$, the element\/ $\Lambda=(\xi-\root l\of\mu)^x\lambda^y$ is a
uniformiser of\/ $L$ and\/ $\ogoth_L=\ogoth[\Lambda]$.  Finally, if\/ $D=\bar
U_{le_1}$, then\/ $\lambda$ remains a uniformiser of\/ $L$, the class of\/
$\mu=1-\eta l(1-\zeta)$ generates\/ $D$ for any\/ $\eta\in\ogoth$ whose image
generates\/ $k/\wp(k)$, and\/ $\ogoth_L=\ogoth[(\root l\of\mu-1)/(1-\zeta)]$.
\endth

Only the last part (the one about $\bar U_{le_1}$) does not follow from the
proof of the previous prop.~; this part was already dealt with in cor.~{37}.

\th COROLLARY {62}
\enonce
Let\/ $K$ be a finite extension of\/ $\Ql(\zeta)$, $L$ a degree-$l$ cyclic
extension of\/ $K$, $(G_i)_{i\in[-1,+\infty[}$ the ramification filtration
of\/ $G=\Gal(L|K)$, $t$ the integer such that\/ $G_t=G$, $G_{t+1}=\{\Id_L\}$.
Then\/ $t=-1$ if $L|K$ is unramified~; otherwise\/ $t\in\;[1,le_1]$ and\/ $t$
is prime to\/~$l$, unless $t=le_1$.  Each such\/ $t$ occurs for some\/ $L$~;
$t=le_1$ occurs only when\/ $L$ is the splitting field of\/ $T^l-\lambda$ for
some uniformiser\/ $\lambda$ of\/ $K$.
\endth
This is clear if $L|K$ is unramified, for then the inertia subgroup $G_0$ is
trivial.  Otherwise, let $D=\Ker(K^\times\!/K^{\times l}\rightarrow
L^\times\!/L^{\times l})$ be the $\F_l$-line in $\bar U_0=K^\times\!/K^{\times
l}$ which corresponds to $L$~; we have $D\neq\bar U_{le_1}$ (prop.~{16}).  Let
$m$ be the integer such that $D\subset\bar U_m$ but $D\not\subset\bar
U_{m+1}$~; we have $m<le_1$ and $m$ is prime to $l$ (prop.~{42}) unless $m=0$,
which can happen only when $D$ is generated by the image of some uniformiser
of $K$.  As $m+t=le_1$ (prop.~{60}), the statement follows.
\medbreak
Hasse \citer\klassenkorperbericht(p.~251) first proves that in the ramified
case $t\in[1,le_1]$ is prime to\/~$l$ unless $t=le_1$ and uses it to conclude
from the equality $m+t=le_1$ (\citer\klassenkorperbericht(p.~266)) that
$m\in[0,le_1[$ is prime to $l$ unless $m=0$, a fact which we had seen directly
in prop.~{42}.  Indeed, cor.~{62} can be used to recover
\citer\klassenkorperbericht(p.~251), which specifies the break in the
ramification filtration of a $(\Z/l\Z)$-extensions $L|K$, and the
possibilities for $v_K(d_{L|K})$, even when the finite extension $K$ of $\Ql$
does {\it not\/} contain an element of order~$l$~:


\th PROPOSITION {63}
\enonce
Suppose that\/ $K^\times$ does\/ {\rm not} contain an element of order~$l$,
let\/ $L|K$ be a degree-$l$ cyclic extension, $G=\Gal(L|K)$, and\/ $t$ the 
integer such that\/ $G_t=G$ but\/ $G_{t+1}=\{\Id_L\}$. Then\/ $t=-1$ if\/
$L|K$ is unramified~; otherwise\/ $t\in\;[1,le_1[$ and\/ $t$ is prime
to\/~$l$. 
\endth
Clearly $t=-1$ if the extension $L|K$ is unramified, so assume that it is
(totally) ramified.  Let's first give two proofs of the fact that
$t\in\;[1,le_1]$ and that $t$ is prime to $l$ if $t\neq le_1$~; then we will
show that $t\neq le_1$.

Put $K'=K(\zeta)$ (where $\zeta^l=1$ but $\zeta\neq1$) and $L'=LK'$.  The
degree of the extension $K'|K$, and {\it a fortiori\/} the relative
ramification index $s$, divide the degree $l-1$ (cf.~prop.~{23}) of the
extension $\Ql(\zeta)|\Ql$.  The extension $L'|K'$ is a totally ramified
kummerian extension of degree~$l$.  Denoting by $t'$ the break in the
ramification filtration of $\Gal(L'|K')$, it is sufficient, by cor.~{62}, to
show that $t'=st$.  (Note that the absolute ramification index of $K'$ is
$s.(l-1)e_1$, where $(l-1)e_1$ is the absolute ramification index of $K$~; the
number $e_1$ need not be an integer, cf.~prop.~{25}.)

Now, $L'|L$ is a tamely ramified extension whose relative ramification index
is also $s$~; choose a uniformiser $\varpi$ of $L'$ such that $\varpi^s$ is a
uniformiser of $L$.  Also, the restriction map
$\Gal(L'|K')\rightarrow\Gal(L|K)$ is an isomorphism~; choose a generator
$\sigma$ of these groups.  By the definition of $t$,
$$\eqalign{
s.(t+1)&=s.v_L(\sigma(\varpi^s)-\varpi^s)\cr
&=v_{L'}(\sigma(\varpi)^s-\varpi^s).\cr
}$$
Let $\xi$ be an element of order $s$ in $\Ql^\times$ (recall that $s$
divides $l-1$) and write
$$
\sigma(\varpi)^s-\varpi^s
=(\sigma(\varpi)-\varpi)
(\sigma(\varpi)-\xi\varpi)\ldots(\sigma(\varpi)-\xi^{s-1}\varpi). 
$$ 
By the definition of $t'$, the $L'$-valuation of the first factor is $t'+1$~;
let us show that it is $1$ for the other $s-1$ factors.  For $0<i<s$, the
element $1-\xi^i$ is a unit of $L'$~; writing
$\sigma(\varpi)-\xi^i\varpi=(\sigma(\varpi)-\varpi)+(1-\xi^i)\varpi$ and
noting that $v_{L'}(\sigma(\varpi)-\varpi)=t'+1>1$ because $t'>0$ (cor.~{62}),
we have $v_{L'}(\sigma(\varpi)-\xi^i\varpi)=1$.  Therefore
$v_{L'}(\sigma(\varpi)^s-\varpi^s)=(t'+1)+(s-1)=t'+s$, and thus
$$
s(t+1)=v_{L'}(\sigma(\varpi)^s-\varpi^s)=t'+s.
$$
Hence $t'=st$, as desired.  Cf.~\citer\childs(p.~127).

(Notice, before moving on, that the claim about $K'|K$ being totally ramified
made in \citer\childs(p.~127) is not correct when $p\neq2$.  Take, for
example, $K=\Qp(\root p-1\of p)$, which is totally ramified of degree $p-1$
and admits $\pi=\root p-1\of p$ as a uniformiser.  Writing
$-p=\varepsilon\pi^{p-1}$, we deduce $\varepsilon=-1$~; in particular
$\bar\varepsilon\notin\F_p^{\times(p-1)}$, therefore $\zeta\notin K$
(prop.~{25}).  On the other hand, the extension $K(\root p-1\of{-1})$ is
unramified over $K$ and contains $\zeta$, for it contains $\root p-1\of{-p}$
(prop.~{24}).  Therefore $K'=K(\zeta)$ is unramified over $K$~; as the degree
of $K'|K$ is $>1$, it cannot be totally ramified. More generally, if the
ramification index of a finite extension $K|\Qp$ is a multiple of the
ramification index $p-1$ of $\Qp(\zeta)|\Qp$, then $K(\zeta)|K$ is unramified,
as follows from Abhyankar's lemma \citer\narkiewicz(p.~236).)

Remark that, with the above notation, the equality $t'=st$ can also be derived
by computing the different ${\goth D}_{L'|K}$ in two different ways by
transitivity in the square
$$
\diagram{
L&\hfl{s-1}{s}{10mm}&L'\cr
\ufl{(l-1)(1+t)}{l}{5mm}&&\ufl{l}{(l-1)(1+t')}{5mm}\cr
K&\hfl{s}{s-1}{10mm}&K'\cr
}
$$
in which the internal letters indicate the relative ramification index and the
external expressions the valuation of the different of the extension in
question.  Using the tower of extensions $L'|L|K$, we get 
$$
v_{L'}({\goth D}_{L'|K})
=v_{L'}({\goth D}_{L'|L})+s.v_{L}({\goth D}_{L|K})
=(s-1)+s.(1-l)(1+t)
$$
whereas using the tower $L'|K'|K$, we get
$$
v_{L'}({\goth D}_{L'|K})
=v_{L'}({\goth D}_{L'|K'})+l.v_{K'}({\goth D}_{K'|K})
=(1-l)(1+t')+l.(s-1).
$$
Comparing the two expressions, we deduce $t'=st$, as claimed.  The hypothesis
${}_lK^\times=\one$ has not been used so far.

Let us show finally that if $L|K$ is a degree-$l$ cyclic extension whose
ramification break occurs at $le_1$, then $K^\times$ has an element of
order~$l$ (and therefore $L|K$ is kummerian).  Let $\lambda$ (resp.~$\Lambda$)
be a uniformiser of $K$ (resp.~$L$)~; $\Lambda^l/\lambda$ is a unit of $L$.
By assumption, for any generator $\sigma\in\Gal(L|K)$, we have
$$
{\sigma(\Lambda)\over\Lambda}
\equiv1+\theta\lambda^{e_1}\pmod{\Lambda^{le_1+1}}
$$
for some $\theta\in k^\times$ invertible in the common residue field $k$ of
$K$ and $L$.  Applying $N_{L|K}$, we get
$1\equiv(1+\theta\lambda^{e_1})^l\pmod{\lambda^{le_1+1}}$, which we take to
mean that the map $(\ )^l:U_{e_1}/U_{e_1+1}\to U_{le_1}/U_{le_1+1}$ is not
injective.  This is possible only if $K^\times$ has an element of order~$l$
(prop.~{29}).  Cf.~\citer\fesvost(p.~75).

(Note in passing that a finite extension $K$ of $\Ql$ ($l$ being an odd prime)
may very well not contain a primitive $l$-th root of~$1$ and still admit
ramified $(\Z/l\Z)$-extensions, contrary to the inadvertent exercise~30 in
\citer\cohen(p.~280).  The simplest example would be the unique degree-$l$
extension $L$ of $K=\Ql$ contained in $M=\Ql(\root l\of\zeta)$, where
$\zeta^l=1$ but $\zeta\neq1$.  Recall that $M|\Ql$ is totally ramified
galoisian with $\Gal(M|\Ql)=(\Z/l^2\Z)^\times$ (cf.~prop.~{23}), a group which
has a unique order-$l$ quotient (prop.~{34}).  Another aside~: for {\it
  every\/} ramified $(\Z/l\Z)$-extension of $\Ql$ ($l\neq2$), the break in the
ramification filtration occurs at $t=1$, because $le_1<2$.  Consequently, for
every $(\Z/l\Z)$-extension of $K'=\Ql(\zeta)$ which comes from a
$(\Z/l\Z)$-extension of $K=\Ql$ (such as $M|K'$; indeed $M=LK'$), the
ramification break occurs at $l-1$ (and hence $\bar\zeta\notin\bar U_2$).
However, $\Q_2$ has two quadratic extensions with $t=1$ and four with $t=2$,
because there are two $\F_2$-lines (generated respectively by $\overline{-1}$
and $\overline{-u}$, where $u=1+2^2$) in $\Q_2^\times\!/\Q_2^{\times2}$ of
``level'' $m=1$ and four lines with $m=0$~; the remaining line $\{\bar1,\bar
u\}$ is of level $m=2=2e_1$ and gives the unramified quadratic extension, for
which of course $t=-1$.)

\th COROLLARY {64}
\enonce
As\/ $L$ runs through the degree-$l$ cyclic extensions of\/ $K$, the
possibilities for\/ $v_K(d_{L|K})$ are\/ $0$, $(l-1)(1+t)$ (for\/ $0<t<le_1$
prime to\/ $l$) and\/ $(l-1)(1+le_1)$.
\endth

Let us note finally that the number of degree-$l$ kummerian extensions of $K$
(when $K^\times$ has an element of order~$l$) which have a given ramification
break $t$ (equivalently, a given valuation of the discriminant or a given
``level'' $m$) can be determined using cor.~{43}.  For every positive integer
$n$, denote by 
$$
\delta_l(n)={l^{1+n}-1\over l-1}=1+l+l^2+\cdots+l^n
$$
the number of lines in the vector $\F_l$-space of dimension $1+n$~;
equivalently, $\delta_l(n)$ is the number of points in $\P_n(\F_l)$.
According to our current convention, $\bar U_0=K^\times\!/K^{\times l}$, a
vector space of dimension $2+d$.  Also, $\bar U_1$ is of dimension $1+d$, so
the number of $\F_l$-points of $\P(\bar U_0)$ which are not in $\P(\bar U_1)$
is $\delta_l(1+d)-\delta_l(d)$.  Thus $\delta_l(1+d)-\delta_l(d)=l^{1+d}$ is
the number of degree-$l$ cyclic extension of $K$ whose ramification break
occurs at $le_1$, or, equivalently, the normalised valuation of whose
discriminant is $(l-1)(1+le_1)$.  For example, $\Ql(\zeta)$ has $l^l$
extensions whose ramification break occurs at $l$.  In other words, there are
$l^l$ extensions of $\Ql(\zeta)$ the valuation of whose discriminant is
$l^2-1$, where the valuation of $1-\zeta$ is~$1$.

Let us come to the other possibilities for the ramification break (cor.~{62}).
Define
$$\mu(t)=\left(t-\left[t\over l\right]\right)\!f\quad(t\in[0,le_1[\,),\qquad 
\mu(le_1)=1+d.
$$ 
We have seen (cor.~{43}) that $\mu(t)$ is the
$\F_l$-dimension of the projective space $\P(\bar U_m)$, with $m=le_1-t$ and
$\bar U_0=K^\times\!/K^{\times l}$, for every $t\in[0,le_1]$.

For $t\in[1,le_1[$, the number of degree-$l$ cyclic extensions of $K$ having a
ramification break at $t$ is the number of $\F_l$-points in $\P(\bar U_m)$
($m=le_1-t$) which are not in $\P(\bar U_{m+1})$~; this number equals
$\delta_l(\mu(t))-\delta_l(\mu(t-1))$.  Notice that $\mu(le_1-1)=d$, so the
number of degree-$l$ cyclic extensions of $K$ whose ramification break occurs
at $le_1$ can also be written as $\delta_l(\mu(le_1))-\delta_l(\mu(le_1-1))$.
When $t=1$, the number of extensions is
$$
\delta_l(\mu(1))-\delta_l(\mu(0))
=\delta_l(f)-\delta_l(0)=\delta_l(f)-1=l+l^2+\cdots+l^f.
$$
For example, $\Ql(\zeta)$ has exactly $l^t$ degree-$l$
cyclic extensions whose ramification break occurs at~$t$, for every
$t\in[1,l]$. 

In general, we have proved 

\th COROLLARY {65}
\enonce
Suppose that\/ $K^\times$ has an element of order\/ $l$.  For\/
$t\in[1,le_1]$, the number of degree-$l$\/ kummerian extensions of\/ $K$ whose
ramification break occurs at\/ $t$ is\/
$\delta_l(\mu(t))-\delta_l(\mu(t-1))$~; this number vanishes when\/ $t$ is a
multiple of\/ $l$, except when\/ $t=le_1$. 
\endth
More concretely, write $l+l^2+\ldots+l^d=n_1+n_2+\ldots+n_e$, where $n_1$ is
the sum of the first $f$ terms on the left, $n_2$ the sum of the next $f$
terms, $\ldots$, $n_e$ the sum of the last $f$ terms.  There are exactly $e$
numbers $1=t_1<t_2<\ldots<t_e=le_1-1$ in $[1,le_1[$ which are prime to~$l$~;
also put $t_{e+1}=le_1$ and $n_{e+1}=l^{1+d}$.  We have seen (cor.~{62}) that
the ramification break of a degree-$l$ cyclic extension $L|K$ occurs at
(precisely) one of $-1$, $t_1$, $\ldots$, $t_e$, $t_{e+1}$.  Cor.~{65} says
that the number of degree-$l$ cyclic extensions $L|K$ having these breaks is
respectively $1$, $n_1$, $\ldots$, $n_e$, $n_{e+1}$.

Let us rewrite the number of (totally) ramified degree-$l$ cyclic extensions
of $K$.  The number in question is the sum over the number of degree-$l$
cyclic extensions whose ramification break occurs at $t_1$, $\ldots$, $t_e$,
$t_{e+1}$ respectively, as indicated~:
$$
\overbrace{\strut n_1}^{t_1}+
\overbrace{\strut \cdots}^{\cdots}+
\overbrace{\strut n_e}^{t_e}+
\overbrace{\strut n_{e+1}}^{t_{e+1}}=
\overbrace{\strut l+l^2+\cdots+l^d}^{t_1,\ldots,t_e}+
\overbrace{\strut l^{1+d}}^{t_{e+1}}
=\delta_l(1+d)-1.
$$

\th COROLLARY {66}
\enonce
Suppose that\/ $K^\times$ has an element of order\/ $l$.  The number of
degree-$l$\/ kummerian extensions\/ $L|K$ with\/ $v_K(d_{L|K})=(l-1)(1+t_i)$
is\/ $n_i=p^{(i-1)f+1}+\ldots+p^{if}$ for\/ $i\in[1,e]$ but\/
$n_{e+1}=l^{1+d}$ for\/ $i=e+1$.
\endth

This allows us to compute the contribution of degree-$l$ kummerian extensions
to the degree-$l$ ``mass formula'' of Serre.  Recall that this formula states
that 
$$
\sum_{L}{1\over l^{c(L)f}}=l,
$$ 
where $L$ runs through all (totally) ramified degree-$l$ extensions of
$K$ and $c(L)=v_K(d_{L|K})-(l-1)$ \citer\serremass().  On the other hand,
cor.~{65} allows us (when $\zeta\in K$) to compute the contribution to this
sum coming from {\it kummerian\/} $L|K$.  If the ramification break occurs at
$t_i$ ($i\in[1,e+1]$), then $c(L)=(l-1)t_i$, and the number of such extensions
is $n_i$.  Thus
$$
\sum_{i\in[1,e+1]}{n_i\over l^{(l-1)t_if}}\leqno{(11)}
$$
is the contribution of degree-$l$ kummerian extensions to Serre's
degree-$l$ mass formula.  This gives the proportion of ramified degree-$l$
extensions which are kummerian.

When $l=2$, we are dealing with quadratic extensions, which are all
kummerian, so the sum $(11)$ must equal~$2$.  This can be verified directly~:

\th LEMMA {67}
\enonce
Let\/ $e>0$ and $f>0$ be integers.  We have the identity
$$
{2+2^2+\cdots+2^f\over 2^{(2-1)f}}+
\cdots+
{2^{(e-1)f+1}+2^{(e-1)f+2}+\cdots+2^{ef}\over 2^{(2e-1)f}}+
{2^{1+ef}\over 2^{2ef}}
=2.
$$
\endth\noindent
For $i\in[1,e]$, the numerator of the $i$-th term on the left-hand
side is 
$$\eqalign{
2^{(i-1)f+1}+2^{(i-1)f+2}+\cdots+2^{if}
&=2^{(i-1)f}(2+2^2+\cdots+2^{f})\cr
&=2^{(i-1)f+1}(2^f-1),\cr
}$$
so the $i$-th term is $2(2^f-1)/2^{if}$.  Thus the sum of the first
$e$ terms is
$$\eqalign{
2(2^f-1)\left({1\over2^f}+{1\over2^{2f}}+\cdots+{1\over2^{ef}}\right)
&={2^{ef}-1\over2^{ef-1}}\cr
&=2-{1\over2^{ef-1}},\cr
}$$
which, when added to the $(e+1)$-th term, gives $2$, proving the statement.

{\it Final remark.}  The perceptive reader must have noticed that the problem
of computing the relative discriminant of a finite extension $L|K$ of {\it
number fields\/} has been reduced to the computation of the relative
discriminant of a local kummerian extension of degree equal to the residual
characteristic, and that prop.~{60} and prop.~{45} solve this local problem.
Let us briefly indicate the argument for reducing the global problem to the
local problem, and the local problem to the prime-degree kummerian case.

The relative discriminant ${\goth d}_{L|K}$ is an integral ideal of $K$, so it
is determined by the knowledge of the exponents $v_\pgoth({\goth d}_{L|K})$
with which various primes $\pgoth$ of $K$ appear in its prime decomposition.
This exponent is a purely local invariant~: it is the valuation of the
discriminant of the {\'e}tale $K_\pgoth$-algebra $L\otimes_KK_\pgoth$.  This
valuation is the sum of the valuations of the discriminants of the various
finite extensions of $K_\pgoth$ into which $L\otimes_KK_\pgoth$ splits.  In
other words, it is sufficient to know how to compute the discriminant of an
extension of local fields.

So let $L|K$ now be a finite extension of $\Qp$ ($p$ prime).  The first
reduction is to the galoisian case~: if $M|K$ is a galoisian extension
containing $L$, then $M|L$ is also galoisian, and if we know how to compute
the differents ${\goth D}_{M|K}$ and ${\goth D}_{M|L}$ of these galoisian
extensions, then we can compute ${\goth D}_{L|K}$, since
$$
{\goth D}_{M|K}={\goth D}_{M|L}.i_{M|L}({\goth D}_{L|K}),
$$ 
of which the {\it Schachtelungssatz\/} (8) is but an avatar, and where
$i_{M|L}$ takes ideals of $L$ to ideals of $M$.  We have already used this
formula in the more conceptual proof of prop.~{63}.

So assume now that $L|K$ is galoisian.  The next reduction is to the case of a
totally ramified $p$-extension.  Indeed, there are intermediate extensions
$L|L_1|L_0|K$ such that $L_0$ is the maximal unramified extension of $K$ in
$L$ and $L_1$ is the maximal tamely ramified extension of $L_0$ in $L$~; the
extension $L|L_1$ has degree a power of $p$.  As we know how to compute the
valuations $v_K(d_{L_0|K})=0$ and $v_{L_0}(d_{L_1|L_0})=[L_1:L_0]-1$
(prop.~{59}), it suffices to know how to compute $v_{L_1}(d_{L|L_1})$.

So assume now that $L|K$ is a (finite) $p$-extension, and let $G=\Gal(L|K)$.
The $p$-group $G$ admits a finite decreasing sequence of subgroups $(G_n)_n$
(we are {\it not\/} talking about higher ramification groups here) such that
$G_{n+1}$ is normal in $G_n$ and each quotient $G_{n}/G_{n+1}$ is a group of
order $p$.  We are thus reduced to the case of (cyclic) degree-$p$ extensions.

So assume now that $L|K$ is a cyclic degree-$p$ extension.  The valuation of
the discriminant $v_K(d_{L|K})$ can be computed if we know the break $t$ in
the ramification filtration of $\Gal(L|K)$, which is given by the formula
$t'=st$, where $s$ is the ramification index of $K'=K(\zeta)$ ($\zeta^p=1$,
$\zeta\neq1$) over $K$, and $t'$ is the break in the ramification filtration
of the degree-$p$ kummerian extension $L'|K'$, with $L'=LK'$~; cf.~the proof
of prop.~{63}.

So assume now that $\zeta\in K$ and that $L|K$ is a degree-$p$ cyclic
extension~; it corresponds to an $\Fp$-line $D\subset K^\times\!/K^{\times
p}$.  As we saw (prop.~{16}), $L|K$ is unramified ({\it le cas non
ramifi{\'e}}) precisely when $D=\bar U_{pe_1}$ (where $(p-1)e_1$ is the
absolute ramification index of $K$).  In this case $v_K(d_{L|K})=0$, of
course. 

Assume finally that $L|K$ is a ramified kummerian degree-$p$ extension, so it
corresponds to a line $D\neq\bar U_{pe_1}$, and denote by $m$ the integer such
that $D\subset\bar U_m$ but $D\not\subset\bar U_{m+1}$~; we have
$m\in[0,pe_1[$, and $m$ is prime to $p$ unless $m=0$.  The break in the
ramification filtration of $L|K$ occurs at $pe_1-m$ and therefore
$v_K(d_{L|K})=(p-1)(1+pe_1-m)$ (prop.~{60}).  As for computing the integer
$m$, there are two cases.  Either $D$ can be generated by the image of a
uniformiser ({\it le cas tr{\`e}s ramifi{\'e}\/}), in which case $m=0$, or it
can be generated by the image $\bar u$ of a unit ({\it le cas peu
  ramifi{\'e}\/}), in which case $m$ is the exponent in the highest power
$\pgoth^m$ of the maximal ideal $\pgoth$ of the ring of integers $\ogoth$ of
$K$ modulo which $u$ is a $p$-th power~: $\hat u\in(\ogoth/\pgoth^{m})^{\times
  p}$ but $\hat u\notin(\ogoth/\pgoth^{m+1})^{\times p}$ (prop.~{45}).

This solves the problem in the local kummerian case (of degree equal to the
residual characteristic), and hence the global problem of determining the
relative discriminant of an extension of number fields.  Taking the base field
to be $\Q$, this also helps decide if a given order in a number field is the
maximal order, which happens precisely when the discriminant of the order
equals the absolute discriminant of the number field in question.

\bigbreak
\centerline{\bf VIII. Discriminants of elliptic curves over local fields}
\bigskip

Let $p$ be a prime number and let $K$ be a finite extension of $\Qp$.  The
basic reason why the discriminant $d_{L|K}$ of a finite extension $L|K$ is an
element of $K^\times\!/\ogoth^{\times 2}$ is that the discriminant $d_B\in
K^\times$ of an $\ogoth$-basis $B$ of $\ogoth_L$ changes by an element of
$\ogoth^{\times 2}$ if we change $B$ (to another $\ogoth$-basis of
$\ogoth_L$).  Moreover, $\ogoth$-bases always exist, and any two
$\ogoth$-bases of $\ogoth_L$ differ by an $\ogoth$-automorphism of the module
$\ogoth_L$.

Now let $E$ be an elliptic curve over $K$.  It can be defined by a minimal
weierstrassian cubic $f$, whose discriminant $d_f$ --- the result of
elliminating the indeterminates $x,y$ from $f,f'_x,f'_y$ --- belongs to
$K^\times$.  If we replace $f$ by another minimal weierstrassian cubic $g$
defining $E$, then $d_f$ gets multiplied by an element of $\ogoth^{\times12}$.
Moreover, any two minimal weierstrassian cubics defining $E$ can be changed
into each other.  This suggests the definition of the discriminant of $E$ as
the class in $K^\times\!/\ogoth^{\times12}$ of the discriminant of any minimal
weierstrassian cubic defining $E$ (def.~68).

More precisely, if $ f=y^2+a_1xy+a_3y-x^3-a_2x^2-a_4x-a_6=0 $ is a minimal
weierstrassian cubic defining $E$, then, according to \citer\tate(p.~180), 
$$
d_f=-b_2^2b_8-2^3b_4^3-3^3b_6^2+3^2b_2b_4b_6
$$
where
$$
b_2=a_1^2+2^2a_2,\quad b_4=a_1a_3+2a_4,\quad b_6=a_3^2+2^2a_6
$$
and
$$
b_8=b_2a_6-a_1a_3a_4+a_2a_3^2-a_4^2.
$$
If we replace $f$ by another minimal weierstrassian cubic $g=0$, there is
change of variables $x=u^2x'+r$, $y=u^3y'+u^2sx'+t$ with $u\in\ogoth^\times$
and $r,s,t\in\ogoth$.  One has $d_{f}=u^{12}d_{g}$ \citer\tate(p.~181), which
leads to the following definition.

\th DEFINITION {68}
\enonce
Let\/ $E$ be an elliptic curve over\/ $K$.  The discriminant\/ $d_{E|K}\in
K^\times\!/\ogoth^{\times12}$ of\/ $E$ is the class of the discriminant
$d_f\in K^\times$ of any of the minimal weierstrassian cubics $f$ defining\/
$E$.
\footnote{$({}^4)$}{\rm It turns out that the global version of this
definition, which we also arrived at \citer\minimal() --- the discriminant of
an elliptic curve over a number field as an
id{\`e}le-modulo-twelfth-powers-of-unit-id{\`e}les, in perfect analogy with
Fr{\"o}hlich's definition of the id{\'e}lic discriminant of an extension of
number fields as an id{\`e}le-modulo-squares-of-unit-id{\`e}les --- has
already been considered by J.~Silverman \citer\silverman().  He was inspired
by the same paper \citer\frohlich() as us, and used it practically for the
same purpose --- a criterion for the existence of a global minimal
weierstrassian cubic. }
\endth

If $E$ has good reduction, then $d_{E|K}\in\ogoth^\times\!/\ogoth^{\times12}$.
We might ask which elements of this finite group occur as $d_{E|K}$ for some
(good-reduction elliptic curve) $E$ over $K$.  We shall see that they all do~;
this should be contrasted with th.~{6} and prop.~{15} which exclude, when
$p=2$, certain elements of $\ogoth^\times\!/\ogoth^{\times2}$ from being
discriminants of (unramified) extensions of $K$.

\th PROPOSITION {69}
\enonce
Suppose that $p\neq2,3$ and let $\delta\in\units$.  There exists a minimal
weierstrassian cubic $f$ such that $d_f=\delta$.
\endth

The invariants 
$c_4=b_2^2-2^3.3.b_4$, $c_6=-b_2^3+2^2.3^2b_2b_4-2^3.3^3.b_6$, 
of a weierstrassian cubic $f$ satisfy $c_6^2=c_4^3-2^6.3^3.d_f$.  In
imitation, consider the cubic $\Gamma:\eta^2=\xi^3-2^6.3^3.\delta$.  Its
discriminant is $-2^4.3^3.(2^6.3^3)^2\delta^2$, so $\Gamma$ is an elliptic
curve over~$K$~; it even has good reduction~$\tilde\Gamma$.  Let $(\xi,\eta)$
be a point in $\Gamma(K)$ with $\xi,\eta\in\ogoth$, for example a point whose
reduction is $\neq0$ in $\tilde\Gamma(k)$, a group which is not reduced to
$\zero$ because it has at least $q+1-2\sqrt q$ points (Hasse) and $q$ is at
least~$5$. Then the weierstrassian cubic
$$
y^2=x^3-(\xi/2^3.3)x-(\eta/2^5.3^3),
$$ 
has discriminant $\delta$ (and $c_4=\xi$, $c_6=\eta$).  It is minimal because
its coefficients are in $\ogoth$ and discriminant in
$\units$. Cf.~\citer\krausinv(p.~76). 

\th PROPOSITION {70}
\enonce
Suppose that $p=2$~or~$3$ and let $\delta\in\units$.  There exists a minimal
weierstrassian cubic $f$ such that $d_f=\delta$.  
\endth

The following proof was suggested by Joseph Oesterl{\'e}.  Consider the
weierstrassian cubic $f=y^2+xy-x^3-a_6=0$.  It has 
$$
b_2=1,\ \ \ b_4=0,\ \ \ b_6=4a_6,\ \  \ b_8=a_6,\ \ \
\hbox{ and }\ \ d_f=-a_6-2^4.3^3.a_6^2.
$$
The last equation can be solved, by Hensel's lemma, for $a_6$ when
$d_f=\delta$, both for $p=2$ and for $p=3$.  Explicitly, we write
$$
\eqalign{
a_6&=-\delta-2^4.3^3.a_6^2\cr
&= -\delta-2^4.3^3.(-\delta-2^4.3^3.a_6^2)^2\cr
}
$$
and so on, which also shows that $a_6$ is a unit and hence $f$ minimal.  Thus
every unit is the discriminant of some minimal weierstrassian cubic.    

\th COROLLARY {71}
\enonce
Every element $\delta\in\ogoth^\times\!/\ogoth^{\times12}$ is the
discriminant of a (good-reduction) elliptic curve over~$K$.
\endth

\th COROLLARY {72}
\enonce
Let $k$ be a finite field.  Every element $\delta\in k^\times\!/k^{\times12}$
is the discriminant of an elliptic $k$-curve.
\endth
Let $K|\Q_p$ be the unramified extension whose residue field is $k$, denote by
$\ogoth$ its ring of integers, and let
$\eta\in\ogoth^\times\!/\ogoth^{\times12}$ an element whose image is $\delta$.
There is a good-reduction elliptic curve $E|K$ whose discriminant is $\eta$~;
the discriminant of its reduction is $\delta$.

\bigbreak
\centerline{\bf IX. Starting from $D\equiv0,1\pmod4$}
\bigskip

Before coming to the genesis of these notes, let us summarise their main
features.  First, we have done everything intrinsically, without choosing any
bases for the spaces which appear.  The invariant language of points and lines
and $\Fp$-spaces is also an aid to the imagination~; compare, for example, the
statement of th.~{52} with its translation into our language.  Secondly, we
have emphasised that questions about discriminants are purely local, and that
locally they have analogues for $p$-primary numbers.  Thirdly, a central role
is played by the exact sequence
$$
1\rightarrow
{}_p\mu\rightarrow
U_{{e_1}}/U_{{{e_1}}+1}\hfl{(\ )^p}{}{8mm}
U_{p{{e_1}}}/U_{p{{e_1}}+1}\rightarrow
\bar U_{p{e_1}}\rightarrow1.
\leqno{(6)}$$ 
which becomes isomorphic, upon choosing a primitive $p$-th root of~$1$ or a
$(p-1)$-th root of $-p$, to the exact sequence
$
0\rightarrow\F_p\rightarrow k\rightarrow
k\rightarrow\F_p\rightarrow0
$
involving $\wp$~; the isomorphism has been made explicit.  Finally, we have
systematised a certain number of results from the literature, and seen that
many of them continue lines of enquiry which can be traced back to the
beginning of the $20$-th century and even earlier.

This investigation was begun in response to a question by a student
(K.~Srilakshmi) as to why the absolute discriminant of a number field is
$\equiv0,1\pmod4$.  I told her that it is a purely local matter at the
prime~$2$ and perhaps Cassels \citer\cassels() had a proof.  He had it
indeed~; more importantly, he had a reference to Fr{\"o}hlich
\citer\frohlich() for the id{\'e}lic notion of the discriminant.  It is the
dissatisfaction with the formulation of the absolute version in Cassels and of
the relative version in Fr{\"o}hlich --- after such a beautiful definition ---
which led to prop.~{15} and, somewhat later, generalising from the prime~$2$
to any prime~$p$, to prop.~{16}.  The analogue for good-reduction elliptic
curves (cor.~{71}) was just an amusing afterthought.

At this point, seeking to generalise Fr{\"o}hlich's theorem (th.~{52}) about
rings of integers from dyadic quadratic extensions to kummerian extensions of
degree equal to the residual characteristic, we discovered Hecke's theorems
and, somewhat later, Hilbert's theorems from the {\it Zahlbericht\/} which
Hecke generalises.

The question arose as to what the generalisation of the other part --- the one
about the valuation of the discriminant --- of Hilbert's theorems is.  Luckily
we found the key to this question in Hasse's {\it Klassenk{\"o}rperbericht\/}
in the form of the relationship between the ``level'' of a line and the break
in the ramification filtration of the corresponding kummerian extension.  By
providing a uniformiser in the ramified case, this relationship determines the
ring of integers (prop.~{61}) and the ramification break (prop.~{62}) at one
stroke. The fact that we could recently acquire copies of Hecke \citer\hecke()
and Hilbert \citer\zahlbericht(), and that Hasse (at least Teil~I and Teil~Ia
\citer\klassenkorperbericht()) is available online, was crucial.  Very
recently we learnt about the local proof by Hensel and the papers by
Furtw{\"a}ngler and Eisenstein.  Ironically, the sought-after uniformiser can
be guessed at from the second $\Omega$ in the {\it proof\/} of {\it
  Satz\/}~148 in the {\it Zahlbericht}.

Let it be mentioned that a local version, and a generalisation to all
$p$-primary numbers (for any prime $p$), of Martinet's congruence on the
absolute norm of the relative discriminant of an extension of number fields
\citer\martinet() has recently been obtained by S.~Pisolkar \citer\supriya(). 

Prop.~{16} can be taken to mean that the inertia subgroup $G^0$ of
$G=\Gal(M|K)$ --- where $M=K(\root p\of{K^\times})$ and $K|\Qp$ is a finite
extension containing a primitive $p$-th root of~$1$ --- which happens to be
the same as the higher ramification subgroup $G^1$, equals $\bar
U_{pe_1}^\perp$ (orthogonal for the Kummer pairing), $(p-1)e_1$ being the
absolute ramification index of $K$.  The question arose --- and we had put it
to a student (S.~Das) --- as to the orthogonality $G^n=\bar
U_{pe_1-n+1}^\perp$ (for $n\in[1, pe_1+1]$) of the ramification filtration in
the upper numbering on $G$ with respect to the natural filtration on $\bar
U_0=K^\times\!/K^{\times p}$.  The basic idea was that the filtration on $G$
is uniquely determined by the filtrations on the order-$p$ quotients of $G$,
for which see cor.~{62}.  We have discovered that this orthogonality relation
has recently been established by I.~Del~Corso and R.~Dvornicich
\citer\delcorso(p.~286), albeit for somewhat special $K$, namely those of the
form $K=F(\root p-1\of{F^\times})$ for some finite extension $F$ of $\Qp$.

Let us summarise their proof in the general situation where $K$ is any finite
extension of $\Qp$ containing a primitive $p$-th root of~$1$.

Think of $G$ and $\bar U_0=K^\times\!/K^{\times p}$ as finite-dimensional
$\Fp$-spaces dual to each other under the pairing 
$$
\langle\ ,\ \rangle\,:\,
G\,\times\, K^\times\!/K^{\times p}\longrightarrow{}_p\mu,
\qquad
\langle\sigma,\bar\eta\rangle={\sigma(\xi)\over\xi} 
\quad (\xi^p=\eta),
$$ 
where $\eta$ runs through $K^\times$, with image $\bar\eta$ in $
K^\times\!/K^{\times p}$, and $\xi\in M^\times$ is any $p$-th root --- it
doesn't matter which --- of $\eta$ (``Kummer theory''~; see
\citer\bourbaki(p.~V.84--87)).  For every subspace $H\subset G$,
the subspace $D=H^\perp$ of $\bar U_0$ satisfies
$$
K(\root p\of D)=M^H,\quad \Gal(K(\root p\of D)|K)=G/H.
$$
Conversely, for every subspace $D\subset \bar U_0$, these relations hold with
$H$ the subspace $D^\perp\subset G$. 

Now let $n>0$ be an integer and consider the subspace $G^n\subset G$.  As we
know that $G^1=\bar U_{pe_1}^\perp$ (prop.~{16}), assume that $n>1$.  Let us
determine the hyperplanes of $G$ which contain $G^n$~; clearly $G^n\subset
\bar U_{pe_1}^\perp$.

Every hyperplane $H\neq\bar U_{pe_1}^\perp$ is of the form $H=D^\perp$ for
some $\Fp$-line $D\neq \bar U_{pe_1}$ in $\bar U_0$.  Let $m\in[0,pe_1[$ be
the integer such that $D\subset\bar U_m$ but $D\not\subset\bar U_{m+1}$. (As
we saw in prop.~{42}, $m$ is prime to $p$ unless $m=0$~; this fact plays no
role in the present proof).  We know that the upper as well as the lower
ramification break of $\Gal(K(\root p\of D)|K)=G/H$ occurs at $pe_1-m$ (cf.\
cor.~{62}).  In other words, $(G/H)^{pe_1-m}\neq0$ but $(G/H)^{pe_1+1-m}=0$.
By the compatibility of the upper-numbering filtraton under passage to the
quotient (see \citer\serre(Chapter~IV, Proposition~14), for example), we also
have $(G/H)^n=G^n/(G^n\cap H)$.  This implies that $G^n\subset H$ if and only
if $n>pe_1-m$ (equivalently, $m>pe_1-n$ or still $D\subset\bar
U_{pe_1-n+1}$). 

We have seen that for every $n>0$ and every hyperplane $H=D^\perp$ in $G$, 
$$
G^n\subset H
\ \Longleftrightarrow\ 
D\subset\bar U_{pe_1-n+1}.
$$
As $G^n$ is the intersection of all hyperplanes containing it,
$(G^n)^\perp$ is generated in $\bar U_0$ by the union of all lines contained
in $\bar U_{pe_1-n+1}$.  Therefore $(G^n)^\perp=\bar U_{pe_1-n+1}$ (for
$n\in[1, pe_1+1]$, with the convention that $\bar U_0=K^\times\!/K^{\times
  p}$), as claimed.  In particular $G^{pe_1+1}=\{1\}$ and the breaks in the
ramification filtration in the upper numbering occur at $-1$ ($M|K$ is not
totally ramified), at the $e$ integers in $[1,pe_1]$ which are prime to $p$,
and at $pe_1$~; no break occurs at~$0$ because $M|K$ is a $p$-extension.

The argument can be carried out for every real $t\in[-1,+\infty[$ to determine
the filtration $(G^t)_{t\in[-1,+\infty[}$ in terms of the filtration $(\bar
U_m)_{m\in[0,+\infty[}$.  Combined with our knowledge of $\dim_\Fp\bar U_m$
(cor.~{43}), this allows one to determine the breaks in the ramification
filtration and to compute $v_K(d_{M|K})$ by formula $(10)$, which becomes
applicable after converting the filtration to the lower numbering.

An interesting application can be made by taking $K=F(\root
p-1\of{F^\times})$, where $F$ is {\it any\/} finite extension of $\Qp$~; note
that $K^\times$ has an element of order~$p$ (cf.~prop.~{24}).  Recall that
$F^\times$ has an element of order $p-1$ (cf.~prop.~{22}), and that the
$(\Z/(p-1)\Z)$-module $F^\times/F^{\times (p-1)}$ is free of rank~$2$,
containing the free rank-$1$ submodule
$\ogoth_F^\times/\ogoth_F^{\times(p-1)}$.  As the extension $K_0|F$ obtained
by adjoining the $(p-1)$-th roots of $\ogoth_F^\times$ is unramified of degree
$p-1$ (cf.~the proof of prop.~{16} in the case $l\neq p$), and $K|K_0$ is
totally but tamely ramified of the same degree, we have
$v_F(d_{K|F})=(p-1)(p-2)$, by the {\it Schachtelungssatz\/} $(8)$ and
prop.~{59}.  Equivalently, $v_K({\goth D}_{K|F})=p-2$, by an application of
prop.~{59} or of formula~$(10)$ to the filtered group $\Gal(K|F)$~:
$$
\Gal(K|F)_0=(\ogoth_F^\times\!/\ogoth_F^{\times(p-1)})^\perp,\quad
\Gal(K|F)_n=\{\Id_L\} \ \ (n>0),
$$
where the orthogonal is taken with respect to the Kummer pairing between
$\Gal(K|F)$ and $F^\times\!/F^{\times(p-1)}$, with values in ${}_{p-1}\mu$.  

By another application of the {\it Schachtelungssatz\/}, one can compute
$v_F(d_{M|F})$ in terms of $v_K(d_{M|K})$, which was computed above.
Alternately, one can compute the filtration $(\Gal(M|F)_n)_{n\in\N}$ from our
knowledge of the filtrations $(\Gal(M|K)_n)_{n\in\N}$ and
$(\Gal(K|F)_n)_{n\in\N}$, and thereby recover the value of $v_F(d_{M|F})$ by
applying formula $(10)$.

The interest in $v_F(d_{M|F})$ comes from the fact that $M$ is the compositum
of {\it all\/} \hbox{degree-$p$} extensions of $F$ \citer\delcorso(p.~273).

Note that the orthogonality relation $G^n=\bar U_{pe_1-n+1}^\perp$ is closely
related to explicit formulas for the Hilbert symbol, for which see
\citer\vostokov().

\medskip

It is entirely fitting that these ramblings, which were prompted by a question
from a student as to why $D\equiv0,1\pmod4$, should close with what was in
effect a question to a student, as to whether $G^n=\bar U_{pe_1-n+1}^\perp$
($n\in[1,pe_1]$).

\medskip\centerline{***}\medbreak

{\bf Acknowledgements.}  I am grateful to Joseph Oesterl{\'e} for his timely
encouragement and advice, to Dipendra Prasad for his interest, to Peter
Roquette for his indulgence and his letter, to Victor Abrashkin and Ivan
Fesenko for their emails, to the G{\"o}ttinger Digitalisierungszentrum and the
Cellule Mathdoc in Grenoble for making many of the bibliographic items freely
available online, and to the American Mathematical Society for awarding points
which can be redeemed for \citer\zahlbericht() and \citer\hecke().

I thank the organisers of the 2007 Annual Meeting of the Ramanujan
Mathematical Society in Mangalore, Microsoft Research in Bangalore, and the
universities of Ratisbon and Heidelberg for their invitations to present some
of these results.

\bigbreak
\unvbox\bibbox 

\bye